
\pretolerance=500 \tolerance=1000  \brokenpenalty=5000

\catcode`\@=11


\font\twelvebf=cmbx10 at 12pt
\font\eightrm=cmr8
\font\eighti=cmmi8
\font\eightsy=cmsy8
\font\eightbf=cmbx8
\font\eighttt=cmtt8
\font\eightit=cmti8
\font\eightsl=cmsl8
\font\sevenrm=cmr7
\font\seveni=cmmi7
\font\sevensy=cmsy7
\font\sevenbf=cmbx7

\font\sixrm=cmr6
\font\sixi=cmmi6
\font\sixsy=cmsy6
\font\sixbf=cmbx6

\skewchar\eighti='177 \skewchar\sixi='177
\skewchar\eightsy='60 \skewchar\sixsy='60

\def\tenpoint{%
  \textfont0=\tenrm \scriptfont0=\sevenrm
  \scriptscriptfont0=\fiverm
  \def\rm{\fam\z@\tenrm}%
  \textfont1=\teni  \scriptfont1=\seveni
  \scriptscriptfont1=\fivei
  \def\oldstyle{\fam\@ne\teni}\let\old=\oldstyle
  \textfont2=\tensy \scriptfont2=\sevensy
  \scriptscriptfont2=\fivesy
  \textfont\itfam=\tenit
  \def\it{\fam\itfam\tenit}%
  \textfont\slfam=\tensl
  \def\sl{\fam\slfam\tensl}%
  \textfont\bffam=\tenbf
  \scriptfont\bffam=\sevenbf
  \scriptscriptfont\bffam=\fivebf
  \def\bf{\fam\bffam\tenbf}%
  \textfont\ttfam=\tentt
  \def\tt{\fam\ttfam\tentt}%
  \abovedisplayskip=12pt plus 3pt minus 9pt
  \belowdisplayskip=\abovedisplayskip
  \abovedisplayshortskip=0pt plus 3pt
  \belowdisplayshortskip=4pt plus 3pt 
  \smallskipamount=3pt plus 1pt minus 1pt
  \medskipamount=6pt plus 2pt minus 2pt
  \bigskipamount=12pt plus 4pt minus 4pt
  \normalbaselineskip=12pt
  \setbox\strutbox=\hbox{\vrule height8.5pt depth3.5pt width0pt}%
  \let\bigf@nt=\tenrm
  \let\smallf@nt=\sevenrm
  \normalbaselines\rm}
  
\def\eightpoint{%
  \textfont0=\eightrm \scriptfont0=\sixrm
  \scriptscriptfont0=\fiverm
  \def\rm{\fam\z@\eightrm}%
  \textfont1=\eighti  \scriptfont1=\sixi
  \scriptscriptfont1=\fivei
  \def\oldstyle{\fam\@ne\eighti}\let\old=\oldstyle
  \textfont2=\eightsy \scriptfont2=\sixsy
  \scriptscriptfont2=\fivesy
  \textfont\itfam=\eightit
  \def\it{\fam\itfam\eightit}%
  \textfont\slfam=\eightsl
  \def\sl{\fam\slfam\eightsl}%
  \textfont\bffam=\eightbf
  \scriptfont\bffam=\sixbf
  \scriptscriptfont\bffam=\fivebf
  \def\bf{\fam\bffam\eightbf}%
  \textfont\ttfam=\eighttt
  \def\tt{\fam\ttfam\eighttt}%
  \abovedisplayskip=9pt plus 3pt minus 9pt
  \belowdisplayskip=\abovedisplayskip
  \abovedisplayshortskip=0pt plus 3pt
  \belowdisplayshortskip=3pt plus 3pt 
  \smallskipamount=2pt plus 1pt minus 1pt
  \medskipamount=4pt plus 2pt minus 1pt
  \bigskipamount=9pt plus 3pt minus 3pt
  \normalbaselineskip=9pt
  \setbox\strutbox=\hbox{\vrule height7pt depth2pt width0pt}%
  \let\bigf@nt=\eightrm
  \let\smallf@nt=\sixrm
  \normalbaselines\rm}
\tenpoint

\font\tencal=eusm10

\font\sevencal=eusm7

\font\fivecal=eusm5
\newfam\calfam
\textfont\calfam=\tencal
\scriptfont\calfam=\sevencal
\scriptscriptfont\calfam=\fivecal
\def\cal#1{{\fam\calfam\relax#1}}

\def\pc#1{\bigf@nt#1\smallf@nt}

\catcode`\;=\active
\def;{\relax\ifhmode\ifdim\lastskip>\z@\unskip\fi \kern\fontdimen2 
 -1.2 \fontdimen3 \string;}

\catcode`\:=\active
\def:{\relax\ifhmode\ifdim\lastskip>\z@\unskip\fi\penalty\@M\ 
\fi\string:}

\catcode`\!=\active
\def!{\relax\ifhmode\ifdim\lastskip>\z@ \unskip\fi\kern\fontdimen2 
 -1.1 \fontdimen3 \string!}

\catcode`\?=\active
\def?{\relax\ifhmode\ifdim\lastskip>\z@ \unskip\fi\kern\fontdimen2 
 -1.1 \fontdimen3 \string?}

\frenchspacing

\def\og{\leavevmode\raise.3ex\hbox{$\scriptscriptstyle 
\langle\!\langle\,$}}
\def\fg{\leavevmode\raise.3ex\hbox{$\scriptscriptstyle 
\,\rangle\!\rangle$}}

\def\pointir{\unskip . --- \ignorespaces}


\def\Medbreak{\vskip-\lastskip\medbreak}

\def\rem#1\endrem{%
\Medbreak {\it#1\unskip} : }

\long\def\thm#1 #2\enonce#3\endthm{%
\Medbreak {\pc#1} {#2\unskip}\pointir{\it #3}\medskip}

\def\decale#1{\smallbreak\hskip 28pt\llap{#1}\kern 5pt}
\def\decaledecale#1{\smallbreak\hskip 34pt\llap{#1}\kern 5pt}

\let\@ldmessage=\message

\def\message#1{{\def\pc{\string\pc\space}%
\def\'{\string'}\def\`{\string`}%
\def\^{\string^}\def\"{\string"}%
\@ldmessage{#1}}}


\def\up#1{\raise 1ex\hbox{\smallf@nt#1}}

\def\diagram#1{\def\normalbaselines{\baselineskip=0pt 
\lineskip=5pt}\matrix{#1}}

\def\longmapsto#1{\mapstochar\mathrel{\joinrel \kern-0.2mm\hbox to
#1mm{\rightarrowfill}}}

\catcode`\@=12

\showboxbreadth=-1  \showboxdepth=-1


\message{`lline' & `vector' macros from LaTeX}

\def\Grille{\setbox13=\vbox to 5\unitlength{\hrule width 109mm \vfill}
\setbox13=\vbox to 65mm
{\offinterlineskip\leaders\copy13\vfill\kern-1pt\hrule} 
\setbox14=\hbox to 5\unitlength{\vrule height 65mm\hfill} 
\setbox14=\hbox to 109mm{\leaders\copy14\hfill\kern-2mm \vrule height 
65mm}
\ht14=0pt\dp14=0pt\wd14=0pt \setbox13=\vbox to 0pt
{\vss\box13\offinterlineskip\box14} \wd13=0pt\box13}

\def\rule(#1,#2)\dir(#3,#4)\long#5{%
\noalign{\leftput(#1,#2){\lline(#3,#4){#5}}}}
\def\arrow(#1,#2)\dir(#3,#4)\length#5{%
\noalign{\leftput(#1,#2){\vector(#3,#4){#5}}}}
\def\put(#1,#2)#3{\noalign{\setbox1=\hbox{%
\kern #1\unitlength \raise #2\unitlength\hbox{$#3$}}%
\ht1=0pt \wd1=0pt \dp1=0pt\box1}}

\catcode`@=11

\def\{{\relax\ifmmode\lbrace\else$\lbrace$\fi}
\def\}{\relax\ifmmode\rbrace\else$\rbrace$\fi}
\def\newcount{\alloc@0\count\countdef\insc@unt}
\def\newdimen{\alloc@1\dimen\dimendef\insc@unt}
\def\newwrite{\alloc@7\write\chardef\sixt@@n}

\newwrite\@unused
\newcount\@tempcnta
\newcount\@tempcntb
\newdimen\@tempdima
\newdimen\@tempdimb
\newbox\@tempboxa

\def\@spaces{\space\space\space\space}
\def\@whilenoop#1{}
\def\@whiledim#1\do #2{\ifdim #1\relax#2\@iwhiledim{#1\relax#2}\fi}
\def\@iwhiledim#1{\ifdim #1\let\@nextwhile=\@iwhiledim 
\else\let\@nextwhile=\@whilenoop\fi\@nextwhile{#1}}
\def\@badlinearg{\@latexerr{Bad \string\line\space or \string\vector 
\space argument}}
\def\@latexerr#1#2{\begingroup 
\edef\@tempc{#2}\expandafter\errhelp\expandafter{\@tempc}%

\def\@eha{Your command was ignored.^^JType \space I <command> <return> 
\space to replace it with another command,^^Jor \space <return> \space 
to continue without it.}
\def\@ehb{You've lost some text. \space \@ehc}
\def\@ehc{Try typing \space <return> \space to proceed.^^JIf that 
doesn't work, type \space X <return> \space to quit.}
\def\@ehd{You're in trouble here.  \space\@ehc}
\typeout{LaTeX error.  \space See LaTeX manual for explanation.^^J 
\space\@spaces\@spaces\@spaces Type \space H <return> \space for 
immediate help.}\errmessage{#1}\endgroup}
\def\typeout#1{{\let\protect\string\immediate\write\@unused{#1}}}

\font\tenln = line10
\font\tenlnw = linew10

\newdimen\@wholewidth
\newdimen\@halfwidth
\newdimen\unitlength 

\unitlength =1pt

\def\thinlines{\let\@linefnt\tenln \let\@circlefnt\tencirc 
\@wholewidth\fontdimen8\tenln \@halfwidth .5\@wholewidth}
\def\thicklines{\let\@linefnt\tenlnw \let\@circlefnt\tencircw 
\@wholewidth\fontdimen8\tenlnw \@halfwidth .5\@wholewidth}
\def\linethickness#1{\@wholewidth #1\relax \@halfwidth .5
\@wholewidth}

\newif\if@negarg

\def\lline(#1,#2)#3{\@xarg #1\relax \@yarg #2\relax 
\@linelen=#3\unitlength \ifnum\@xarg =0 \@vline \else \ifnum\@yarg =0 
\@hline \else \@sline\fi \fi}
\def\@sline{\ifnum\@xarg< 0 \@negargtrue \@xarg -\@xarg \@yyarg 
-\@yarg \else \@negargfalse \@yyarg \@yarg \fi
\ifnum \@yyarg >0 \@tempcnta\@yyarg \else \@tempcnta - \@yyarg \fi
\ifnum\@tempcnta>6 \@badlinearg\@tempcnta0 \fi
\setbox\@linechar\hbox{\@linefnt\@getlinechar(\@xarg,\@yyarg)}%
\ifnum \@yarg >0 \let\@upordown\raise \@clnht\z@ 
\else\let\@upordown\lower \@clnht \ht\@linechar\fi
\@clnwd=\wd\@linechar
\if@negarg \hskip -\wd\@linechar \def\@tempa{\hskip -2\wd \@linechar} 
\else \let\@tempa\relax \fi
\@whiledim \@clnwd <\@linelen \do {\@upordown\@clnht\copy\@linechar 
\@tempa \advance\@clnht \ht\@linechar \advance\@clnwd \wd\@linechar}%
\advance\@clnht -\ht\@linechar \advance\@clnwd -\wd\@linechar 
\@tempdima\@linelen\advance\@tempdima -\@clnwd 
\@tempdimb\@tempdima\advance\@tempdimb -\wd\@linechar
\if@negarg \hskip -\@tempdimb \else \hskip \@tempdimb \fi
\multiply\@tempdima \@m\@tempcnta \@tempdima \@tempdima \wd\@linechar 
\divide\@tempcnta \@tempdima \@tempdima \ht\@linechar 
\multiply\@tempdima \@tempcnta \divide\@tempdima \@m \advance\@clnht 
\@tempdima
\ifdim \@linelen <\wd\@linechar \hskip \wd\@linechar 
\else\@upordown\@clnht\copy\@linechar\fi}
\def\@hline{\ifnum \@xarg <0 \hskip -\@linelen \fi
\vrule height \@halfwidth depth \@halfwidth width \@linelen
\ifnum \@xarg <0 \hskip -\@linelen \fi}
\def\@getlinechar(#1,#2){\@tempcnta#1\relax
\multiply\@tempcnta 8\advance\@tempcnta -9
\ifnum #2>0 \advance\@tempcnta #2\relax
 \else\advance\@tempcnta -#2\relax\advance\@tempcnta 64 \fi
\char\@tempcnta}
\def\vector(#1,#2)#3{\@xarg #1\relax \@yarg #2\relax 
\@linelen=#3\unitlength
\ifnum\@xarg =0 \@vvector \else \ifnum\@yarg =0 \@hvector \else 
\@svector\fi \fi}
\def\@hvector{\@hline\hbox to 0pt{\@linefnt \ifnum \@xarg <0 
\@getlarrow(1,0)\hss\else \hss\@getrarrow(1,0)\fi}}
\def\@vvector{\ifnum \@yarg <0 \@downvector \else \@upvector \fi}
\def\@svector{\@sline\@tempcnta\@yarg \ifnum\@tempcnta <0 
\@tempcnta=-\@tempcnta\fi \ifnum\@tempcnta <5 \hskip -\wd\@linechar 
\@upordown\@clnht \hbox{\@linefnt \if@negarg 
\@getlarrow(\@xarg,\@yyarg) \else \@getrarrow(\@xarg,\@yyarg) 
\fi}\else\@badlinearg\fi}
\def\@getlarrow(#1,#2){\ifnum #2 =\z@ \@tempcnta='33\else 
\@tempcnta=#1\relax\multiply\@tempcnta \sixt@@n \advance\@tempcnta -9 
\@tempcntb=#2\relax \multiply\@tempcntb \tw@ \ifnum \@tempcntb >0 
\advance\@tempcnta \@tempcntb\relax \else\advance\@tempcnta 
-\@tempcntb\advance\@tempcnta 64 \fi\fi \char\@tempcnta}
\def\@getrarrow(#1,#2){\@tempcntb=#2\relax \ifnum\@tempcntb < 0 
\@tempcntb=-\@tempcntb\relax\fi \ifcase \@tempcntb\relax 
\@tempcnta='55 \or \ifnum #1<3 \@tempcnta=#1\relax\multiply\@tempcnta 
24 \advance\@tempcnta -6 \else \ifnum #1=3 \@tempcnta=49 
\else\@tempcnta=58 \fi\fi\or \ifnum #1<3 
\@tempcnta=#1\relax\multiply\@tempcnta 24 \advance\@tempcnta -3 \else 
\@tempcnta=51\fi\or \@tempcnta=#1\relax\multiply\@tempcnta \sixt@@n 
\advance\@tempcnta -\tw@ \else \@tempcnta=#1\relax\multiply\@tempcnta 
\sixt@@n \advance\@tempcnta 7 \fi \ifnum #2<0 \advance\@tempcnta 64 
\fi \char\@tempcnta}
\def\@vline{\ifnum \@yarg <0 \@downline \else \@upline\fi}
\def\@upline{\hbox to \z@{\hskip -\@halfwidth \vrule width 
\@wholewidth height \@linelen depth \z@\hss}}
\def\@downline{\hbox to \z@{\hskip -\@halfwidth \vrule width 
\@wholewidth height \z@ depth \@linelen \hss}}
\def\@upvector{\@upline\setbox\@tempboxa 
\hbox{\@linefnt\char'66}\raise \@linelen \hbox to\z@{\lower 
\ht\@tempboxa \box\@tempboxa\hss}}
\def\@downvector{\@downline\lower \@linelen \hbox to 
\z@{\@linefnt\char'77\hss}}

\thinlines

\newcount\@xarg
\newcount\@yarg
\newcount\@yyarg
\newcount\@multicnt
\newdimen\@xdim
\newdimen\@ydim
\newbox\@linechar
\newdimen\@linelen
\newdimen\@clnwd
\newdimen\@clnht
\newdimen\@dashdim
\newbox\@dashbox
\newcount\@dashcnt
\catcode`@=12

\newbox\tbox
\newbox\tboxa

\def\leftzer#1{\setbox\tbox=\hbox to 0pt{#1\hss}%
\ht\tbox=0pt \dp\tbox=0pt \box\tbox}
\def\rightzer#1{\setbox\tbox=\hbox to 0pt{\hss #1}%
\ht\tbox=0pt \dp\tbox=0pt \box\tbox}
\def\centerzer#1{\setbox\tbox=\hbox to 0pt{\hss #1\hss}%
\ht\tbox=0pt \dp\tbox=0pt \box\tbox}

\def\leftput(#1,#2)#3{\setbox\tboxa=\hbox{%
\kern #1\unitlength \raise #2\unitlength\hbox{\leftzer{#3}}}%
\ht\tboxa=0pt \wd\tboxa=0pt \dp\tboxa=0pt\box\tboxa}

\def\rightput(#1,#2)#3{\setbox\tboxa=\hbox{%
\kern #1\unitlength \raise #2\unitlength\hbox{\rightzer{#3}}}%
\ht\tboxa=0pt \wd\tboxa=0pt \dp\tboxa=0pt\box\tboxa}

\def\centerput(#1,#2)#3{\setbox\tboxa=\hbox{%
\kern #1\unitlength \raise #2\unitlength\hbox{\centerzer{#3}}}%
\ht\tboxa=0pt \wd\tboxa=0pt \dp\tboxa=0pt\box\tboxa}

\unitlength=1mm

\expandafter\ifx\csname amssym.def\endcsname\relax \else
\endinput\fi
%
\expandafter\edef\csname amssym.def\endcsname{%
       \catcode`\noexpand\@=\the\catcode`\@\space}
\catcode`\@=11
%

\def\undefine#1{\let#1\undefined}
\def\newsymbol#1#2#3#4#5{\let\next@\relax
 \ifnum#2=\@ne\let\next@\msafam@\else
 \ifnum#2=\tw@\let\next@\msbfam@\fi\fi
 \mathchardef#1="#3\next@#4#5}
\def\mathhexbox@#1#2#3{\relax
 \ifmmode\mathpalette{}{\m@th\mathchar"#1#2#3}%
 \else\leavevmode\hbox{$\m@th\mathchar"#1#2#3$}\fi}
\def\hexnumber@#1{\ifcase#1 0\or 1\or 2\or 3
\or 4\or 5\or 6\or 7\or 8\or
 9\or A\or B\or C\or D\or E\or F\fi}

\font\tenmsa=msam10
\font\sevenmsa=msam7
\font\fivemsa=msam5
\newfam\msafam
\textfont\msafam=\tenmsa
\scriptfont\msafam=\sevenmsa
\scriptscriptfont\msafam=\fivemsa
\edef\msafam@{\hexnumber@\msafam}
\mathchardef\dabar@"0\msafam@39
\def\dashrightarrow{\mathrel{\dabar@\dabar@\mathchar"0
\msafam@4B}}
\def\dashleftarrow{\mathrel{\mathchar"0\msafam@4C
\dabar@\dabar@}}

\def\ulcorner{\delimiter"4\msafam@70\msafam@70 }
\def\urcorner{\delimiter"5\msafam@71\msafam@71 }
\def\llcorner{\delimiter"4\msafam@78\msafam@78 }
\def\lrcorner{\delimiter"5\msafam@79\msafam@79 }
\def\yen{{\mathhexbox@\msafam@55}}
\def\checkmark{{\mathhexbox@\msafam@58}}
\def\circledR{{\mathhexbox@\msafam@72}}
\def\maltese{{\mathhexbox@\msafam@7A}}

\font\tenmsb=msbm10
\font\sevenmsb=msbm7
\font\fivemsb=msbm5
\newfam\msbfam
\textfont\msbfam=\tenmsb
\scriptfont\msbfam=\sevenmsb
\scriptscriptfont\msbfam=\fivemsb
\edef\msbfam@{\hexnumber@\msbfam}
\def\Bbb#1{{\fam\msbfam\relax#1}}
\def\widehat#1{\setbox\z@\hbox{$\m@th#1$}%
 \ifdim\wd\z@>\tw@ em\mathaccent"0\msbfam@5B{#1}%
 \else\mathaccent"0362{#1}\fi}
\def\widetilde#1{\setbox\z@\hbox{$\m@th#1$}%
 \ifdim\wd\z@>\tw@ em\mathaccent"0\msbfam@5D{#1}%
 \else\mathaccent"0365{#1}\fi}
\font\teneufm=eufm10
\font\seveneufm=eufm7
\font\fiveeufm=eufm5
\newfam\eufmfam
\textfont\eufmfam=\teneufm
\scriptfont\eufmfam=\seveneufm
\scriptscriptfont\eufmfam=\fiveeufm
\def\frak#1{{\fam\eufmfam\relax#1}}

\csname amssym.def\endcsname

\expandafter\ifx\csname pre amssym.tex at\endcsname\relax \else 
\endinput\fi
\expandafter\chardef\csname pre amssym.tex at\endcsname=\the
\catcode`\@
\catcode`\@=11
\begingroup\ifx\undefined\newsymbol \else\def\input#1
{\endgroup}\fi
\input amssym.def \relax
\newsymbol\boxdot 1200
\newsymbol\boxplus 1201
\newsymbol\boxtimes 1202
\newsymbol\square 1003
\newsymbol\blacksquare 1004
\newsymbol\centerdot 1205
\newsymbol\lozenge 1006
\newsymbol\blacklozenge 1007
\newsymbol\circlearrowright 1308
\newsymbol\circlearrowleft 1309
\undefine\rightleftharpoons
\newsymbol\rightleftharpoons 130A
\newsymbol\leftrightharpoons 130B
\newsymbol\boxminus 120C
\newsymbol\Vdash 130D
\newsymbol\Vvdash 130E
\newsymbol\vDash 130F
\newsymbol\twoheadrightarrow 1310
\newsymbol\twoheadleftarrow 1311
\newsymbol\leftleftarrows 1312
\newsymbol\rightrightarrows 1313
\newsymbol\upuparrows 1314
\newsymbol\downdownarrows 1315
\newsymbol\upharpoonright 1316
 
\newsymbol\downharpoonright 1317
\newsymbol\upharpoonleft 1318
\newsymbol\downharpoonleft 1319
\newsymbol\rightarrowtail 131A
\newsymbol\leftarrowtail 131B
\newsymbol\leftrightarrows 131C
\newsymbol\rightleftarrows 131D
\newsymbol\Lsh 131E
\newsymbol\Rsh 131F
\newsymbol\rightsquigarrow 1320
\newsymbol\leftrightsquigarrow 1321
\newsymbol\looparrowleft 1322
\newsymbol\looparrowright 1323
\newsymbol\circeq 1324
\newsymbol\succsim 1325
\newsymbol\gtrsim 1326
\newsymbol\gtrapprox 1327
\newsymbol\multimap 1328
\newsymbol\therefore 1329
\newsymbol\because 132A
\newsymbol\doteqdot 132B
 
\newsymbol\triangleq 132C
\newsymbol\precsim 132D
\newsymbol\lesssim 132E
\newsymbol\lessapprox 132F
\newsymbol\eqslantless 1330
\newsymbol\eqslantgtr 1331
\newsymbol\curlyeqprec 1332
\newsymbol\curlyeqsucc 1333
\newsymbol\preccurlyeq 1334
\newsymbol\leqq 1335
\newsymbol\leqslant 1336
\newsymbol\lessgtr 1337
\newsymbol\backprime 1038
\newsymbol\risingdotseq 133A
\newsymbol\fallingdotseq 133B
\newsymbol\succcurlyeq 133C
\newsymbol\geqq 133D
\newsymbol\geqslant 133E
\newsymbol\gtrless 133F
\newsymbol\sqsubset 1340
\newsymbol\sqsupset 1341
\newsymbol\vartriangleright 1342
\newsymbol\vartriangleleft 1343
\newsymbol\trianglerighteq 1344
\newsymbol\trianglelefteq 1345
\newsymbol\bigstar 1046
\newsymbol\between 1347
\newsymbol\blacktriangledown 1048
\newsymbol\blacktriangleright 1349
\newsymbol\blacktriangleleft 134A
\newsymbol\vartriangle 134D
\newsymbol\blacktriangle 104E
\newsymbol\triangledown 104F
\newsymbol\eqcirc 1350
\newsymbol\lesseqgtr 1351
\newsymbol\gtreqless 1352
\newsymbol\lesseqqgtr 1353
\newsymbol\gtreqqless 1354
\newsymbol\Rrightarrow 1356
\newsymbol\Lleftarrow 1357
\newsymbol\veebar 1259
\newsymbol\barwedge 125A
\newsymbol\doublebarwedge 125B
\undefine\angle
\newsymbol\angle 105C
\newsymbol\measuredangle 105D
\newsymbol\sphericalangle 105E
\newsymbol\varpropto 135F
\newsymbol\smallsmile 1360
\newsymbol\smallfrown 1361
\newsymbol\Subset 1362
\newsymbol\Supset 1363
\newsymbol\Cup 1264
 
\newsymbol\Cap 1265
 
\newsymbol\curlywedge 1266
\newsymbol\curlyvee 1267
\newsymbol\leftthreetimes 1268
\newsymbol\rightthreetimes 1269
\newsymbol\subseteqq 136A
\newsymbol\supseteqq 136B
\newsymbol\bumpeq 136C
\newsymbol\Bumpeq 136D
\newsymbol\lll 136E
 
\newsymbol\ggg 136F
 
\newsymbol\circledS 1073
\newsymbol\pitchfork 1374
\newsymbol\dotplus 1275
\newsymbol\backsim 1376
\newsymbol\backsimeq 1377
\newsymbol\complement 107B
\newsymbol\intercal 127C
\newsymbol\circledcirc 127D
\newsymbol\circledast 127E
\newsymbol\circleddash 127F
\newsymbol\lvertneqq 2300
\newsymbol\gvertneqq 2301
\newsymbol\nleq 2302
\newsymbol\ngeq 2303
\newsymbol\nless 2304
\newsymbol\ngtr 2305
\newsymbol\nprec 2306
\newsymbol\nsucc 2307
\newsymbol\lneqq 2308
\newsymbol\gneqq 2309
\newsymbol\nleqslant 230A
\newsymbol\ngeqslant 230B
\newsymbol\lneq 230C
\newsymbol\gneq 230D
\newsymbol\npreceq 230E
\newsymbol\nsucceq 230F
\newsymbol\precnsim 2310
\newsymbol\succnsim 2311
\newsymbol\lnsim 2312
\newsymbol\gnsim 2313
\newsymbol\nleqq 2314
\newsymbol\ngeqq 2315
\newsymbol\precneqq 2316
\newsymbol\succneqq 2317
\newsymbol\precnapprox 2318
\newsymbol\succnapprox 2319
\newsymbol\lnapprox 231A
\newsymbol\gnapprox 231B
\newsymbol\nsim 231C
\newsymbol\ncong 231D
\newsymbol\diagup 201E
\newsymbol\diagdown 201F
\newsymbol\varsubsetneq 2320
\newsymbol\varsupsetneq 2321
\newsymbol\nsubseteqq 2322
\newsymbol\nsupseteqq 2323
\newsymbol\subsetneqq 2324
\newsymbol\supsetneqq 2325
\newsymbol\varsubsetneqq 2326
\newsymbol\varsupsetneqq 2327
\newsymbol\subsetneq 2328
\newsymbol\supsetneq 2329
\newsymbol\nsubseteq 232A
\newsymbol\nsupseteq 232B
\newsymbol\nparallel 232C
\newsymbol\nmid 232D
\newsymbol\nshortmid 232E
\newsymbol\nshortparallel 232F
\newsymbol\nvdash 2330
\newsymbol\nVdash 2331
\newsymbol\nvDash 2332
\newsymbol\nVDash 2333
\newsymbol\ntrianglerighteq 2334
\newsymbol\ntrianglelefteq 2335
\newsymbol\ntriangleleft 2336
\newsymbol\ntriangleright 2337
\newsymbol\nleftarrow 2338
\newsymbol\nrightarrow 2339
\newsymbol\nLeftarrow 233A
\newsymbol\nRightarrow 233B
\newsymbol\nLeftrightarrow 233C
\newsymbol\nleftrightarrow 233D
\newsymbol\divideontimes 223E
\newsymbol\varnothing 203F
\newsymbol\nexists 2040
\newsymbol\Finv 2060
\newsymbol\Game 2061
\newsymbol\mho 2066
\newsymbol\eth 2067
\newsymbol\eqsim 2368
\newsymbol\beth 2069
\newsymbol\gimel 206A
\newsymbol\daleth 206B
\newsymbol\lessdot 236C
\newsymbol\gtrdot 236D
\newsymbol\ltimes 226E
\newsymbol\rtimes 226F
\newsymbol\shortmid 2370
\newsymbol\shortparallel 2371
\newsymbol\smallsetminus 2272
\newsymbol\thicksim 2373
\newsymbol\thickapprox 2374
\newsymbol\approxeq 2375
\newsymbol\succapprox 2376
\newsymbol\precapprox 2377
\newsymbol\curvearrowleft 2378
\newsymbol\curvearrowright 2379
\newsymbol\digamma 207A
\newsymbol\varkappa 207B
\newsymbol\Bbbk 207C
\newsymbol\hslash 207D
\undefine\hbar
\newsymbol\hbar 207E
\newsymbol\backepsilon 237F
\catcode`\@=\csname pre amssym.tex at\endcsname

\magnification=1200
\hsize=160 true mm
\vsize=240 true mm
\hoffset=-2mm
\voffset=8mm
\parindent=12pt   \parskip=2pt     
\hfuzz=1pt
\pageno=1
\headline={\ifnum\pageno=1\else\tenrm\hfil 873-\number\folio\hfil\fi}
\footline={\hfil}

\noindent\vtop to 65 true mm{%
\hbox to 160 true mm{S\'eminaire Bourbaki\hfil Mars 2000}
\hbox to 160 true mm{52\`{e}me ann\'{e}e, 1999-2000,
n\raise4pt\hbox{o} 873\hss}
\vfil
\centerline{\bf LA CORRESPONDANCE DE LANGLANDS}
\centerline{\bf SUR LES CORPS DE FONCTIONS}
\bigskip
\centerline{\bf [d'apr\`{e}s Laurent LAFFORGUE]}
\vfil
\centerline{\bf {\rm par} G\'{e}rard LAUMON}
\vfil}\par

\topskip=1 true cm
\baselineskip=13pt

Lafforgue a r\'{e}cemment \'{e}tabli la correspondance de Langlands
pour $\mathop{\rm GL}\nolimits_{r}$ sur un corps de fonctions. Sa
preuve suit la strat\'{e}gie introduite, il y a plus de 25 ans, par
Drinfeld pour traiter le cas $r=2$. Une des innovations principales
est la construction d'une compactification toro\"{\i}dale du
sch\'{e}ma simplicial $\mathop{\rm PGL} \nolimits_{r}^{\bullet
+1}/\mathop{\rm PGL}\nolimits_{r}$ classifiant de $\mathop{\rm
PGL}\nolimits_{r}$ qui prolonge la compactification de De Concini et
Procesi de $\mathop{\rm PGL}\nolimits_{r}^{2}/\mathop{\rm
PGL}\nolimits_{r}$.

Apr\`{e}s avoir \'{e}nonc\'{e} le th\'{e}or\`{e}me principal et ses
cons\'{e}quences, nous pr\'{e}senterons les grandes lignes de la
d\'{e}monstration, en renvoyant pour les d\'{e}tails aux publications
de Lafforgue ([La 1] \`{a} [La 10]). Les sections 4 (Compactification
du classifiant de $\mathop{\rm PGL}\nolimits_{r}$) et 7 (Une variante
d'un th\'{e}or\`{e}me de Pink) sont de nature g\'{e}n\'{e}rale et
peuvent \^{e}tre lues ind\'{e}pendamment du reste du texte.
\vskip 3mm

Je remercie chaleureusement L. Lafforgue pour son aide dans
la pr\'{e}paration de cet expos\'{e}.

\vskip 7mm

\centerline{\bf 1. \'{E}NONC\'{E} DU TH\'{E}OR\`{E}ME PRINCIPAL}
\vskip 10mm

On fixe dans tout cet expos\'{e} une courbe $X$ projective, lisse et
g\'{e}om\'{e}triquement connexe sur un corps fini ${\Bbb F}_{q}$ \`{a}
$q$ \'{e}l\'{e}ments. On note $|X|$ l'ensemble des points ferm\'{e}s
de $X$.

Soit $F$ le corps des fonctions de $X$. On identifie les places de $F$
aux \'{e}l\'{e}ments de $|X|$. Pour chaque $x\in |X|$ on peut former
le compl\'{e}t\'{e} $F_{x}$ de $F$ en la place $x$. C'est un corps de
valuation discr\`{e}te complet. On notera encore $x:F_{x}^{\times}
\rightarrow {\Bbb Z}$ sa valuation. L'anneau des entiers de $F_{x}$
est l'anneau de valuation discr\`{e}te ${\cal O}_{x}=\{a_{x}\in
F_{x}^{\times}\mid x(a_{x})\geq 0\}\cup\{0\}$, ${\frak
p}_{x}=\{a_{x}\in F_{x}^{\times}\mid x(a_{x})>0\}\cup\{0\}$ est
l'unique id\'{e}al maximal de ${\cal O}_{x}$ et le corps r\'{e}siduel
$\kappa (x)={\cal O}_{x}/{\frak p}_{x}$ est une extension finie de
${\Bbb F}_{q}$ dont on notera $\mathop{\rm deg}(x)$ le degr\'{e}.

L'anneau (topologique) des ad\`{e}les de $F$ est le produit restreint
$$
{\Bbb A}=\{a=(a_{x})_{x\in |X|}\mid a_{x}\in {\cal O}_{x}\hbox{ pour
presque tout }x\}\subset\prod_{x\in |X|}^{}F_{x}
$$
o\`{u} l'expression {\og}{pour presque tout $x$}{\fg} signifie
{\og}{pour tous les $x$ sauf un nombre fini}{\fg}. Le corps $F$ se
plonge diagonalement dans ${\Bbb A}$ et l'anneau topologique compact
${\cal O}:=\prod_{x\in |X|}^{}{\cal O}_{x}$ est un sous-anneau de
${\Bbb A}$. 

On dispose d'un homomorphisme de groupes $\mathop{\rm deg}:{\Bbb
A}^{\times}\rightarrow {\Bbb Z}$ d\'{e}fini par
$$
\mathop{\rm deg}(a)=\sum_{x\in |X|}^{}\mathop{\rm deg}(x)x(a_{x}).
$$
Cet homomorphisme est surjectif et est identiquement nul sur
$F^{\times}$ et aussi sur ${\cal O}^{\times}$. Son noyau est compact
modulo $F^{\times}$: en d'autres termes, pour tout $a\in {\Bbb
A}^{\times}$ de degr\'{e} non nul, le groupe quotient $F^{\times}
\backslash {\Bbb A}^{\times}/{\cal O}^{\times}a^{{\Bbb Z}}$ est fini.
\vskip 5mm

{\bf 1.1. Repr\'{e}sentations automorphes cuspidales}
\vskip 5mm

Soit $r$ un entier $\geq 1$. On consid\`{e}re le groupe ad\'{e}lique
$\mathop{\rm GL}\nolimits_{r}({\Bbb A})$ des matrices inversibles de
taille $r\times r$ et \`{a} coefficients dans l'anneau ${\Bbb A}$, et
ses sous-groupes $\mathop{\rm GL}\nolimits_{r}(F)$ et
$$
K=\prod_{x\in |X|}K_{x}=\mathop{\rm GL}\nolimits_{r}({\cal
O})=\prod_{x\in |X|}\mathop{\rm GL}\nolimits_{r}({\cal O}_{x}).
$$
On fixe une fois pour toute une mesure de Haar $dg$ sur $\mathop{\rm
GL}\nolimits_{r}({\Bbb A})$, et une d\'{e}composition $dg=\prod_{x\in
|X|}^{}dg_{x}$ de cette mesure en produit de mesures de Haar locales.
Il est commode de normaliser $dg$ et les $dg_{x}$ par les conditions
$\mathop{\rm vol}(K,dg)=1$ et $\mathop{\rm vol}(K_{x},dg_{x})=1$.

Une {\it forme automorphe cuspidale} (pour $\mathop{\rm
GL}\nolimits_{r}^{}$ sur $F$) est une fonction $\varphi :\mathop{\rm
GL}\nolimits_{r}({\Bbb A})\rightarrow {\Bbb C}$ ayant les
propri\'{e}t\'{e}s suivantes:
\vskip 2mm

\itemitem{1)} $\varphi (\gamma g)=\varphi (g)$, $\forall\gamma\in
\mathop{\rm GL}\nolimits_{r}(F)$, $\forall g\in\mathop{\rm
GL}\nolimits_{r}({\Bbb A})$,
\vskip 2mm

\itemitem{2)} il existe un sous-groupe $K_{\varphi}\subset K$ d'indice
fini tel que $\varphi (gk)= \varphi (g)$, $\forall g\in\mathop{\rm
GL}\nolimits_{r}({\Bbb A})$, $\forall k\in K_{\varphi}$,
\vskip 2mm

\itemitem{3)} il existe $a\in {\Bbb A}^{\times}$ tel que $\mathop{\rm
deg}(a)\not=0$ et $\varphi (ga)=\varphi (g)$, $\forall g\in
\mathop{\rm GL}\nolimits_{r}({\Bbb A})$,
\vskip 2mm

\itemitem{4)} pour toute d\'{e}composition non triviale
$r=r_{1}+\cdots +r_{s}$ de $r$ en entiers strictement positifs, qui
d\'{e}finit un sous-groupe parabolique standard
$P=MU\subsetneq\mathop{\rm GL}\nolimits_{r}$ de radical unipotent $U$
et de composante de Levi $M\cong\mathop{\rm
GL}\nolimits_{r_{1}}\times\cdots\times\mathop{\rm
GL}\nolimits_{r_{s}}$, le {\it terme constant}
$$
\mathop{\rm GL}\nolimits_{r}({\Bbb A})\rightarrow {\Bbb
C},~g\mapsto\int_{U(F)\backslash U({\Bbb A})}\varphi (ug)du,
$$
est identiquement nul (ici $du$ est n'importe quelle mesure de Haar
sur le quotient compact $U(F)\backslash U({\Bbb A})$).
\vskip 3mm

Si l'on fixe un \'{e}l\'{e}ment $a_{0}$ de degr\'{e} non nul dans
${\Bbb A}^{\times}$, on a
$$
L_{{\rm cusp}}=\bigcup_{n\geq 1}^{}L_{{\rm cusp}}(a_{0}^{n})
$$
o\`{u}, pour tout $a\in {\Bbb A}^{\times}$ de degr\'{e} non nul,
$L_{{\rm cusp}}(a)\subset L_{{\rm cusp}}$ est le sous-espace
d\'{e}fini en imposant ce $a$ particulier dans la propri\'{e}t\'{e}
3).

Toujours pour $a\in {\Bbb A}^{\times}$ de degr\'{e} non nul, on montre
que toute fonction $\varphi\in L_{{\rm cusp}}(a)$ est \`{a} support
compact sur $\mathop{\rm GL}\nolimits_{r}(F)\backslash \mathop{\rm
GL}\nolimits_{r}({\Bbb A})/a^{{\Bbb Z}}$, et on peut donc munir $L_{{\rm
cusp}}(a)$ du produit scalaire hermitien d\'{e}fini positif
$$
(\varphi_{1},\varphi_{2}):=\int_{\mathop{\rm GL}\nolimits_{r}(F)
\backslash \mathop{\rm GL}\nolimits_{r}({\Bbb A})/a^{{\Bbb Z}}}
\overline{\varphi_{1}(g)}\varphi_{2}(g)dg.
$$
\vskip 3mm

Le groupe $\mathop{\rm GL}\nolimits_{r}({\Bbb A})$ agit par
translation \`{a} droite sur l'espace vectoriel complexe
$$
L_{{\rm cusp}}=L_{{\rm cusp}}(\mathop{\rm GL}\nolimits_{r}(F)
\backslash \mathop{\rm GL}\nolimits_{r}({\Bbb A}))
$$
des formes automorphes cuspidales. Cette repr\'{e}sentation est {\it
lisse} (le fixateur de tout vecteur dans $L_{{\rm cusp}}$ est un
sous-groupe ouvert de $\mathop{\rm GL}\nolimits_{r}({\Bbb A})$). Comme
on a fix\'{e} une mesure de Haar $dg$ sur $\mathop{\rm GL}
\nolimits_{r}({\Bbb A})$, la donn\'{e}e de cette repr\'{e}sentation
\'{e}quivaut \`{a} celle d'une structure de ${\cal H}$-module sur
$L_{{\rm cusp}}$, o\`{u} l'{\it alg\`{e}bre de Hecke}
$$
{\cal H}={\cal C}_{{\rm c}}^{\infty}(\mathop{\rm GL}
\nolimits_{r}({\Bbb A}))
$$
est l'alg\`{e}bre de convolution des fonctions complexes, localement
constantes et \`{a} support compact, sur $\mathop{\rm GL}
\nolimits_{r}({\Bbb A})$.

Pour chaque $a\in {\Bbb A}^{\times}$ de degr\'{e} non nul, l'action de
$\mathop{\rm GL}\nolimits_{r}({\Bbb A})$ respecte le sous-espace
$L_{{\rm cusp}}(a)$. La repr\'{e}sentation induite sur ce sous-espace
est {\it admissible} (pour tout sous-groupe $K'\subset K$ d'indice
fini, l'espace vectoriel des invariants sous $K'$ dans $L_{{\rm
cusp}}(a)$ est de dimension finie). Elle est de plus {\it unitaire}
pour le produit scalaire d\'{e}fini ci-dessus. La repr\'{e}sentation
de $\mathop{\rm GL}\nolimits_{r}({\Bbb A})$ sur $L_{{\rm cusp}}$ est
donc semi-simple: elle admet une d\'{e}composition isotypique
$$
L_{{\rm cusp}}\cong\bigoplus_{\pi}^{}V_{\pi}^{\oplus m(\pi )}
$$
o\`{u} $\pi$ parcourt un syst\`{e}me de repr\'{e}sentants des classes
d'isomorphie de re\-pr\'{e}\-sen\-ta\-tions complexes
irr\'{e}ductibles admissibles de $\mathop{\rm GL}\nolimits_{r}({\Bbb
A})$, o\`{u} $V_{\pi}$ est l'espace de $\pi$ (en g\'{e}n\'{e}ral de
dimension infinie) et o\`{u} les multiplicit\'{e}s $m (\pi )$ sont des
entiers $\geq 0$.

\thm D\'{E}FINITION 
\enonce
Une {\rm repr\'{e}sentation automorphe cuspidale irr\'{e}ductible}
{\rm (}pour $\mathop{\rm GL}\nolimits_{r}$ sur $F${\rm )} est une
repr\'{e}sentation complexe admissible irr\'{e}ductible de
$\mathop{\rm GL}\nolimits_{r}({\Bbb A})$ qui est isomorphe \`{a} un
facteur direct de $L_{{\rm cusp}}$.
\endthm

On notera ${\cal A}_{r}$ un syst\`{e}me de repr\'{e}sentants des
classes d'isomorphie de ces re\-pr\'{e}\-sen\-ta\-tions automorphes
cuspidales irr\'{e}ductibles. Chaque $\pi\in {\cal A}_{r}$ admet un
caract\`{e}re central $\omega_{\pi}:F^{\times}\backslash {\Bbb
A}^{\times}\rightarrow {\Bbb C}^{\times}$ qui est d'ordre fini puisque
$\omega_{\pi}(a)=1$ pour au moins un $a\in {\Bbb A}^{\times}$ de
degr\'{e} non nul.

\thm TH\'{E}OR\`{E}ME (Th\'{e}or\`{e}me de multiplicit\'{e} un, [PS 
2], [Sh])
\enonce
On a 
$$
L_{{\rm cusp}}\cong\bigoplus_{\pi\in {\cal A}_{r}}^{}V_{\pi}.
$$
En d'autres termes, les multiplicit\'{e}s dans la d\'{e}composition
isotypique de la re\-pr\'{e}\-sen\-ta\-tion $L_{{\rm cusp}}$ de
$\mathop{\rm GL}\nolimits_{r}({\Bbb A})$ sont donn\'{e}es par $m(\pi
)=1$ si $\pi$ est automorphe cuspidale et $m(\pi )=0$ sinon.
\endthm

Pour tout $x\in |X|$, on d\'{e}finit aussi l'{\it alg\`{e}bre de Hecke
locale}
$$
{\cal H}_{x}={\cal C}_{{\rm c}}^{\infty}(\mathop{\rm 
GL}\nolimits_{r}^{}(F_{x})).
$$
C'est l'alg\`{e}bre de convolution pour la mesure de Haar $dg_{x}$ sur
$\mathop{\rm GL}\nolimits_{r}^{}(F_{x})$. L'espace de toute
repr\'{e}sentation lisse de $\mathop{\rm GL}\nolimits_{r}^{}(F_{x})$
est naturellement muni d'une structure de ${\cal H}_{x}$-module.

On notera $e_{K_{x}}\in {\cal H}_{x}$ la fonction caract\'{e}ristique
de $K_{x}=\mathop{\rm GL}\nolimits_{r}^{}({\cal O}_{x})$ dans
$\mathop{\rm GL}\nolimits_{r}^{}(F_{x})$ et on appellera alg\`{e}bre
de Hecke locale {\it non ramifi\'{e}e} la sous-alg\`{e}bre
$$
{\cal H}_{x}(K_{x})=e_{K_{x}}\ast {\cal H}_{x}\ast e_{K_{x}}\subset
{\cal H}_{x}.
$$
C'est une alg\`{e}bre commutative unitaire munie d'un isomorphisme,
construit par Satake,
$$
{\cal H}_{x}(K_{x})\cong {\Bbb C}[z_{1},z_{1}^{-1},\ldots
,z_{r},z_{r}^{-1}]^{{\frak S}_{r}},
$$
o\`{u} le groupe sym\'{e}trique ${\frak S}_{r}$ agit par permutation
des ind\'{e}termin\'{e}es $z_{1},\ldots ,z_{r}$. En particulier, si
$(\pi_{x},V_{\pi_{x}})$ est une repr\'{e}sentation admissible
irr\'{e}ductible de $\mathop{\rm GL}\nolimits_{r} (F_{x})$ {\it non
ramifi\'{e}e}, c'est-\`{a}-dire telle que
$V_{\pi_{x}}^{K_{x}}\not=(0)$, on a en fait
$$
\mathop{\rm dim}(V_{\pi_{x}}^{K_{x}})=1
$$
et le ${\cal H}_{x}(K_{x})$-module $V_{\pi_{x}}^{K_{x}}$ de rang $1$
et aussi la repr\'{e}sentation $(\pi_{x},V_{\pi_{x}})$ sont (\`{a}
isomorphisme pr\`{e}s) uniquement d\'{e}termin\'{e}s par la
donn\'{e}e d'un $r$-uplet non ordonn\'{e}
$$
(z_{1}(\pi_{x}),\ldots ,z_{r}(\pi_{x}))
$$
de nombres complexes non nuls, appel\'{e}s les {\it valeurs propres
de Hecke} de $\pi_{x}$.

L'alg\`{e}bre de Hecke globale ${\cal H}$ est le produit tensoriel
restreint
$$
{\cal H}=\lim_{{\scriptstyle\longrightarrow\atop\scriptstyle S}}
\left(\bigotimes_{x\in S}{\cal H}_{x}\right)\otimes
\left(\bigotimes_{x\notin S}e_{K_{x}}\right)
$$
des alg\`{e}bres de Hecke locales, o\`{u} $S$ parcourt l'ensemble des
parties finies de $|X|$. Une repr\'{e}sentation admissible
irr\'{e}ductible $(\pi ,V_{\pi})$ de $\mathop{\rm GL}
\nolimits_{r}^{}({\Bbb A})$ admet donc, pour chaque $x\in |X|$, une
composante locale $(\pi_{x},V_{\pi_{x}})$ qui est une
repr\'{e}sentation irr\'{e}ductible admissible de $\mathop{\rm GL}
\nolimits_{r}^{}(F_{x})$ bien d\'{e}finie \`{a} isomorphisme pr\`{e}s.
Pour presque tout $x$, $\pi$ est non ramifi\'{e}e en $x$,
c'est-\`{a}-dire admet une composante locale non ramifi\'{e}e en $x$,
et $\pi$ est le produit tensoriel restreint
$$
\pi\cong \bigotimes_{x\in |X|}\!\!{}^{\prime}\,\pi_{x}
$$
de ses composantes locales. (Une fois que l'on a fix\'{e} une base
$v_{x}$ de $V_{\pi_{x}}^{K_{x}}$ pour chaque place $x$ non
ramifi\'{e}e pour $\pi$, l'espace de la repr\'{e}sentation du membre
de droite est par d\'{e}finition la limite inductive
$$
\lim_{{\scriptstyle\longrightarrow\atop\scriptstyle S}}
\left(\bigotimes_{x\in S}V_{\pi_{x}}\right)\otimes
\left(\bigotimes_{x\notin S}v_{x}\right)
$$
pour $S$ parcourant l'ensemble des parties finies de $|X|$ qui
contiennent toutes les places ramifi\'{e}es pour $\pi$.)

On d\'{e}finit la {\it fonction $L$ partielle} de $\pi$ comme le
produit eul\'{e}rien formel
$$
L^{S_{\pi}}(\pi ,T)=\prod_{x\notin S_{\pi}}L(\pi_{x},T)
=\prod_{x\notin S_{\pi}}{1\over \prod_{i=1}^{r}
(1-z_{i}(\pi_{x})T^{\mathop{\rm deg}(x)})}\in {\Bbb C}[[T]]
$$
o\`{u} $S_{\pi}\subset |X|$ est l'ensemble fini des places
ramifi\'{e}es de $\pi$. Godement et Jacquet ([G-J]) ont
d\'{e}montr\'{e} que le produit eul\'{e}rien $L^{S_{\pi}}(\pi
,q^{-s})$ converge absolument dans un demi-plan $\mathop{\rm
Re}(s)>\sigma$ pour un nombre r\'{e}el $\sigma$ assez grand, et que la
s\'{e}rie formelle en $T$ d\'{e}finie par $L^{S_{\pi}}(\pi ,T)$ est en
fait le d\'{e}veloppement d'une fraction rationnelle dans ${\Bbb
C}(T)$, et m\^{e}me d'un polyn\^{o}me dans ${\Bbb C}[T]$ si $r\geq 2$.
\vskip 5mm

{\bf 1.2. Repr\'{e}sentations du groupe de Galois}
\vskip 5mm

On fixe une cl\^{o}ture s\'{e}parable $\overline{F}$ de $F$ et on note
$\Gamma_{F}$ le groupe de Galois de $\overline{F}$ sur $F$. Pour
chaque $x\in |X|$ on choisit arbitrairement un sous-groupe de
d\'{e}composition $D_{x}\subset\Gamma_{F}$ en $x$. Ce groupe de
d\'{e}composition est le groupe de Galois $\Gamma_{F_{x}}$ d'une
certaine cl\^{o}ture s\'{e}parable $\overline{F}_{x}$ de $F_{x}$. On
note $I_{x}\subset D_{x}$ le sous-groupe d'inertie de ce groupe de
d\'{e}composition. Le groupe quotient $D_{x}/I_{x}$ est isomorphe au
groupe de Galois $\Gamma_{\kappa (x)}$ de $\overline{\kappa (x)}$ sur
$\kappa (x)$ pour une certaine cl\^{o}ture alg\'{e}brique
$\overline{\kappa (x)}$ de $\kappa (x)$ et est donc isomorphe \`{a}
$\mathop{\rm Frob} \nolimits_{x}^{\widehat{{\Bbb Z}}}$ o\`{u}
$\mathop{\rm Frob} \nolimits_{x}\in\Gamma_{\kappa (x)}$ est
l'\'{e}l\'{e}ment de Frobenius g\'{e}om\'{e}trique (l'inverse de
l'\'{e}l\'{e}vation \`{a} la puissance $|\kappa (x)|=q^{\mathop{\rm
deg}(x)}$).

Muni de la topologie de Krull, $\Gamma_{F}$ est un groupe topologique
pro-fini, dont les repr\'{e}sentations complexes continues et de
dimension finie se factorisent n\'{e}cessairement par un quotient fini.
Suivant Serre et Grothendieck, on obtient une cat\'{e}gorie plus vaste
de repr\'{e}sentations de $\Gamma_{F}$ en rempla\c{c}ant le corps des
coefficients ${\Bbb C}$ par un corps $\ell$-adique. Fixons donc un
nombre premier auxiliaire $\ell$, distinct de la caract\'{e}ristique de
${\Bbb F}_{q}$, et fixons une cl\^{o}ture alg\'{e}brique
$\overline{{\Bbb Q}}_{\ell}$ de ${\Bbb Q}_{\ell}$.

Une {\it repr\'{e}sentation $\ell$-adique} $\sigma$ de $\Gamma_{F}$
est un $\overline{{\Bbb Q}}_{\ell}$-espace vectoriel $V_{\sigma}$ de
dimension finie, muni d'un homomorphisme de groupes $\sigma
:\Gamma_{F}\rightarrow\mathop{\rm Aut}\nolimits_{\overline{{\Bbb
Q}}_{\ell}}(V_{\sigma})$, ayant les propri\'{e}t\'{e}s suivantes:
\vskip 2mm

\itemitem{1)} il existe une base de $V_{\sigma}$ identifiant
$\mathop{\rm Aut}\nolimits_{\overline{{\Bbb Q}}_{\ell}}(V_{\sigma})$
\`{a} $\mathop{\rm GL}\nolimits_{r}(\overline{{\Bbb Q}}_{\ell})$,
o\`{u} $r$ est la dimension de $V_{\sigma}$, et une extension finie
$E_{\lambda}$ de ${\Bbb Q}_{\ell}$ contenue dans $\overline{{\Bbb
Q}}_{\ell}$ telles que $\sigma (\Gamma_{F})\subset \mathop{\rm
GL}\nolimits_{r}(E_{\lambda}) \subset\mathop{\rm
GL}\nolimits_{r}(\overline{{\Bbb Q}}_{\ell})$,
\vskip 2mm

\itemitem{2)} $\sigma :\Gamma_{F}\rightarrow \mathop{\rm GL}
\nolimits_{r}^{}(E_{\lambda})$ est continu pour la topologie de Krull
sur $\Gamma_{F}$ et la topologie $\ell$-adique sur $\mathop{\rm
GL}\nolimits_{r}^{}(E_{\lambda})$,
\vskip 2mm

\itemitem{3)} pour presque tout $x\in |X|$, $\sigma$ est {\it non
ramifi\'{e}e} en $x$, c'est-\`{a}-dire que la restriction $\sigma_{x}$
de $\sigma$ au sous-groupe de d\'{e}composition $D_{x}\subset
\Gamma_{F}$ est triviale sur le sous-groupe d'inertie $I_{x}\subset
D_{x}$.
\vskip 3mm

Ces repr\'{e}sentations forment une cat\'{e}gorie ab\'{e}lienne o\`{u}
tout objet est de longueur finie. 

Soit ${\cal G}_{r}$ un syst\`{e}me de repr\'{e}sentants des classes
d'isomorphie de repr\'{e}sentations $\ell$-adiques irr\'{e}ductibles
de rang $r$ de $\Gamma_{F}$ dont le d\'{e}terminant (la puissance
ext\'{e}rieure maximale) est d'ordre fini. Pour chaque $\sigma\in
{\cal G}_{r}$, on note $S_{\sigma}$ l'ensemble fini des places de $F$
o\`{u} $\sigma$ est ramifi\'{e}e. Pour chaque $x\notin S_{\sigma}$, on
dispose de l'automorphisme $\sigma_{x}(\mathop{\rm Frob}
\nolimits_{x})$ de $V_{\sigma}$, et on note
$$
(z_{1}(\sigma_{x}),\ldots ,z_{r}(\sigma_{x}))
$$
les {\it valeurs propres de Frobenius de $\sigma$ en $x$},
c'est-\`{a}-dire le $r$-uplet non ordonn\'{e} d'\'{e}l\'{e}ments non
nuls de $\overline{{\Bbb Q}}_{\ell}$ form\'{e} des valeurs propres de
$\sigma_{x}(\mathop{\rm Frob} \nolimits_{x})$. La {\it fonction $L$
partielle} de $\sigma$ est par d\'{e}finition le produit eul\'{e}rien
formel
$$
L^{S_{\sigma}}(\sigma ,T)=\prod_{x\in |X|-S_{\sigma}}
L(\sigma_{x},T)=\prod_{x\in |X|-S_{\sigma}}{1\over
\prod_{i=1}^{r}(1-z_{i}(\sigma_{x})T^{\mathop{\rm deg}(x)})}\in
\overline{{\Bbb Q}}_{\ell}[[T]].
$$
Il r\'{e}sulte de la formule des points fixes de
Grothendieck-Lefschetz ([Gr]) que cette s\'{e}rie formelle est le
d\'{e}veloppement d'une fraction rationnelle dans $\overline{{\Bbb
Q}}_{\ell}(T)$, et m\^{e}me d'un polyn\^{o}me dans $\overline{{\Bbb
Q}}_{\ell}[T]$ si $r\geq 2$.
\vskip 5mm

{\bf 1.3. La correspondance de Langlands et ses cons\'{e}quences}
\vskip 5mm

On fixe un isomorphisme de corps $\iota :\overline{{\Bbb Q}}_{\ell}
\buildrel\sim\over\longrightarrow {\Bbb C}$ (il en existe d'apr\`{e}s
l'axiome du choix). On notera encore $\iota$ les isomorphismes de
$\overline{{\Bbb Q}}_{\ell}(T)$ (resp. $\overline{{\Bbb
Q}}_{\ell}[[T]]$) sur ${\Bbb C}(T)$ (resp. ${\Bbb C}[[T]]$) induits
par $\iota$.

\thm TH\'{E}OR\`{E}ME {\pc PRINCIPAL} 
\enonce
{\rm (i)} {\rm (Correspondance de Langlands)} Il existe une unique
bijection
$$
{\cal A}_{r}\buildrel\sim\over\longrightarrow {\cal
G}_{r},~\pi\mapsto\sigma (\pi ),
$$
telle que, pour tout $\pi\in {\cal A}_{r}$ on ait l'\'{e}galit\'{e} de
facteurs $L$ locaux
$$
\iota (L(\sigma (\pi )_{x},T))=L(\pi_{x},T),
$$
pour presque tout $x\notin S_{\sigma (\pi )}\cup S_{\pi}$.

\decale{\rm (ii)} {\rm (Conjecture de Ramanujan-Petersson)} Pour tout
$\pi\in {\cal A}_{r}$ et toute place $x\notin S_{\pi}$ on a
$$
|z_{i}(\pi_{x})|=1,~\forall i=1,\ldots ,r.
$$
\endthm

Le cas $r=1$ de ce th\'{e}or\`{e}me est une reformulation de la
th\'{e}orie du corps de classes ab\'{e}lien pour les corps de
fonctions. Le cas $r=2$ a \'{e}t\'{e} d\'{e}montr\'{e} par Drinfeld
([Dr 6], [Dr 7]). Le cas g\'{e}n\'{e}ral est d\^{u} \`{a} Lafforgue
([La 10]).
\vskip 3mm

Les cons\'{e}quences de ce th\'{e}or\`{e}me qui sont formul\'{e}es 
ci-dessous \'{e}taient attendues, et leur d\'{e}duction du 
th\'{e}or\`{e}me principal est bien connue.

\thm TH\'{E}OR\`{E}ME (Compatibilit\'{e} avec la correspondance de 
Langlands locale)
\enonce
Pour tout $\pi\in {\cal A}_{r}$ et tout $x\in |X|$, la
repr\'{e}sentation de Galois locale $\sigma (\pi )_{x}$ de
$D_{x}=\Gamma_{F_{x}}$ est l'image $\sigma_{x}(\pi_{x})$ de la
repr\'{e}sentation locale $\pi_{x}$ de $\mathop{\rm
GL}\nolimits_{r}^{}(F_{x})$ par la correspondance de Langlands locale
{\rm (}cf. {\rm [L-R-S]}{\rm )}. En particulier, $\pi$ et $\sigma (\pi
)$ ont les m\^{e}mes facteurs $L$ locaux en toutes les places de $F$.
\endthm

\thm TH\'{E}OR\`{E}ME (Conjectur\'{e} par Deligne, [De](1.2.10))
\enonce
Soit $\sigma\in {\cal G}_{r}$.

\decale{\rm (i)} Le sous-corps $E=E(\sigma )$ de $\overline{{\Bbb
Q}}_{\ell}$ engendr\'{e} sur ${\Bbb Q}$ par les coefficients des
polyn\^{o}mes $\prod_{i=1}^{r}(1-z_{i}(\sigma_{x})T)$ pour tous les
$x\not\in S_{\sigma}$ est un corps de nombres {\rm (}une extension
finie de ${\Bbb Q}${\rm )}.

\decale{\rm (ii)} Pour tout nombre premier $\ell'$ distinct de la
caract\'{e}ristique de ${\Bbb F}_{q}$, fixons une cl\^{o}ture 
alg\'{e}brique $\overline{{\Bbb Q}}_{\ell'}$ de ${\Bbb Q}_{\ell'}$ et 
notons ${\cal G}_{r,\ell'}$ l'ensemble ${\cal G}_{r}$ correspondant. 

Il existe alors une unique famille de repr\'{e}sentations
$\sigma_{\ell',\lambda'}\in {\cal G}_{r,\ell'}$, index\'{e}e par les
couples $(\ell',\lambda')$ form\'{e}s d'un nombre premier $\ell'$
distinct de la caract\'{e}ristique de ${\Bbb F}_{q}$ et d'un
plongement $\lambda': E\hookrightarrow\overline{{\Bbb Q}}_{\ell'}$, et
ayant les propri\'{e}t\'{e}s suivantes:
\vskip 2mm

\itemitem{-} $\sigma =\sigma_{\ell ,\lambda}$ o\`{u} $\lambda$ est
l'inclusion de $E$ dans $\overline{{\Bbb Q}}_{\ell}$,
\vskip 2mm

\itemitem{-} chaque $\sigma_{\ell',\lambda'}$ a le m\^{e}me ensemble
de places ramifi\'{e}es que $\sigma$ et pour tout $x\notin S_{\sigma}$
on a l'\'{e}galit\'{e}
$$
\prod_{i=1}^{r}(1-z_{i}(\sigma_{\ell' ,\lambda' ,x})T)=
\lambda'\left(\prod_{i=1}^{r}(1-z_{i}(\sigma_{x})T)\right)\in
\overline{{\Bbb Q}}_{\ell'}[T].
$$
\hfill\hfill$\square$
\endthm

\thm TH\'{E}OR\`{E}ME (Conjectur\'{e} par Deligne, [De](1.2.10))
\enonce
Soient $Y$ un sch\'{e}ma normal de type fini sur le corps fini ${\Bbb
F}_{q}$ et ${\cal L}$ un syst\`{e}me local $\ell$-adique de rang $r$
sur $Y$. On suppose que ${\cal L}$ est irr\'{e}ductible et que sa
puissance ext\'{e}rieure maximale $\mathop{\rm det}({\cal L})$ est
d'ordre fini {\rm (}$\mathop{\rm det}({\cal L})^{\otimes N}$ est
isomorphe au faisceau $\ell$-adique constant de valeur
$\overline{{\Bbb Q}}_{\ell}$ sur $Y$ pour tout entier $N\geq 1$ assez
divisible{\rm )}. Alors ${\cal L}$ est pur de poids $0$. 
\hfill\hfill$\square$
\endthm

Une cons\'{e}quence plus indirecte du th\'{e}or\`{e}me principal est
l'\'{e}nonc\'{e} dit {\og}{de descente}{\fg} dans la premi\`{e}re
construction de Drinfeld (cf. [Dr 12], [F-G-K-V], [Lau 5], [Lau 6])
des faisceaux automorphes partout non ramifi\'{e}s pour $\mathop{\rm
GL}\nolimits_{r}^{}$ sur une courbe lisse, projective et
g\'{e}om\'{e}triquement connexe sur un corps arbitraire
(\'{e}ventuellement de caract\'{e}ristique nulle).
\vskip 5mm

{\bf 1.4. Historique}
\vskip 5mm

L'outil principal invent\'{e} par Drinfeld et utilis\'{e} par Drinfeld
et Lafforgue pour d\'{e}montrer le th\'{e}or\`{e}me principal, est le
champ modulaire des chtoucas. Avant de d\'{e}couvrir ce champ,
Drinfeld avait introduit des vari\'{e}t\'{e}s modulaires sur les corps
de fonctions tr\`{e}s semblables aux vari\'{e}t\'{e}s de Shimura sur
les corps de nombres: les vari\'{e}t\'{e}s de {\it modules
elliptiques} (ou de {\it faisceaux elliptiques}). L'utilisation de ces
vari\'{e}t\'{e}s avait permis de construire $\sigma (\pi )$ sous
certaines conditions sur la repr\'{e}sentation locale $\pi_{\infty}$
en une place donn\'{e}e $\infty\in |X|$. Ainsi, Drinfeld ([Dr 9] et
[Dr 10]) avait construit $\sigma (\pi )$ quand $r=2$ et $\pi_{\infty}$
est dans la s\'{e}rie discr\`{e}te, et Flicker et Kazhdan, puis
moi-m\^{e}me, avions g\'{e}n\'{e}ralis\'{e} cette construction de
Drinfeld, pour $r$ arbitraire, dans les cas o\`{u} $\pi_{\infty}$ est
soit supercuspidale ([F-K]), soit la repr\'{e}sentation de Steinberg
([Lau 1], [Lau 2]).

L'assertion d'unicit\'{e} du th\'{e}or\`{e}me principal \'{e}tait
connue depuis longtemps: elle est en effet une cons\'{e}quence
imm\'{e}diate du th\'{e}or\`{e}me de densit\'{e} de \v{C}ebotarev
([Se]). Il en \'{e}tait de m\^{e}me de l'injectivit\'{e} de
l'application $\pi\rightarrow\sigma (\pi )$ qui est automatique
d'apr\`{e}s le th\'{e}or\`{e}me de multiplicit\'{e} un fort de
Piatetski-Shapiro ([PS 2]).

Deligne avait remarqu\'{e} un principe de r\'{e}currence permettant
de d\'{e}duire la surjectivit\'{e} de l'application $\pi\rightarrow
\sigma (\pi )$ de son existence. Plus pr\'{e}cis\'{e}ment, si pour
tout $r'=1,\ldots ,r-1$, on a construit une application $\pi'
\rightarrow\sigma'(\pi')$ de ${\cal A}_{r'}\rightarrow {\cal G}_{r'}$
telle que, quel que soit $\pi'\in {\cal A}_{r'}$, on ait
l'\'{e}galit\'{e} de facteurs $L$ locaux
$$
\iota (L(\sigma'(\pi')_{x},T))=L(\pi_{x}',T)
$$
pour presque tout $x\notin S_{\sigma'(\pi')}\cup S_{\pi'}$, alors
l'\'{e}quation fonctionnelle de Grothendieck ([Gr]), la formule du
produit pour la constante de cette \'{e}quation fonctionnelle ([Lau
7]) et le th\'{e}or\`{e}me inverse de Hecke, Weil et Piatetski-Shapiro
([PS 1], [C-PS]) permettent de d\'{e}finir une application
$$
{\cal G}_{r}\rightarrow {\cal A}_{r},~\sigma\mapsto\pi (\sigma ),
$$
telle que, quel que soit $\sigma\in {\cal G}_{r}$, on ait
l'\'{e}galit\'{e} de facteurs $L$ locaux
$$
\iota (L(\pi (\sigma )_{x},T))=L(\sigma_{x},T)
$$
pour presque tout $x\notin S_{\sigma}\cup S_{\pi (\sigma )}$.

Lafforgue avait obtenu auparavant un \'{e}nonc\'{e} tr\`{e}s proche de
la conjecture de Ramanujan-Petersson ([La 1]) comme une
cons\'{e}quence de la description ad\'{e}lique de Drinfeld des
chtoucas sur les corps finis ([Dr 6]), de la formule des traces
d'Arthur-Selberg ([Ar]), de la formule des points fixes de
Grothendieck-Lefschetz ([Gr]), du th\'{e}or\`{e}me de puret\'{e} de
Deligne ([De]) et des estim\'{e}es de Jacquet-Shalika pour les valeurs
propres de Hecke des repr\'{e}sentations automorphes cuspidales
irr\'{e}ductibles ([J-S 1]).
\vskip 5mm

{\bf 1.5. La strat\'{e}gie}
\vskip 5mm

La strat\'{e}gie utilis\'{e}e par Drinfeld et Lafforgue pour
d\'{e}finir l'application $\pi\rightarrow\sigma (\pi )$ de ${\cal
A}_{r}$ dans ${\cal G}_{r}$ est inspir\'{e}e des travaux de Shimura,
Ihara, Deligne, Langlands, ... Sous sa forme la plus na\"{\i}ve elle
peut se d\'{e}crire comme suit:

On construit un sch\'{e}ma $V$ sur $F$ muni d'une action de
l'alg\`{e}bre de Hecke ${\cal H}$ de telle sorte que la cohomologie
$\ell$-adique \`{a} supports compacts
$$
H_{{\rm c}}^{\ast}(\overline{F}\otimes_{F}V,\overline{{\Bbb
Q}}_{\ell})
$$
soit une repr\'{e}sentation du produit de $\mathop{\rm
GL}\nolimits_{r}^{}({\Bbb A})$ et du groupe de Galois $\Gamma_{F}$.

Puis on calcule la trace de cette repr\'{e}sentation par la formule
des points fixes de Grothendieck-Lefschetz.

Enfin on compare cette formule des points fixes avec la formule des
traces d'Arthur-Selberg pour prouver que la repr\'{e}sentation
$$
\bigoplus_{\pi\in {\cal A}_{r}}\pi\otimes\sigma (\pi )
$$
de $\mathop{\rm GL}\nolimits_{r}^{}({\Bbb A})\times\Gamma_{F}$ que
l'on cherche est exactement la partie {\og}{cuspidale}{\fg} de
$H_{{\rm c}}^{\ast}(\overline{F}\otimes_{F}V, \overline{{\Bbb
Q}}_{\ell})$.
\vskip 3mm

En fait, comme la cohomologie $\ell$-adique ci-dessus est
automatiquement d\'{e}finie sur ${\Bbb Q}_{\ell}$, il y a une
obstruction (de rationalit\'{e}) \`{a} mener \`{a} bien un tel
programme (cf. [Ka]) et la strat\'{e}gie doit \^{e}tre
l\'{e}g\`{e}rement modifi\'{e}e. Comme l'a propos\'{e} Drinfeld, c'est
plut\^{o}t la repr\'{e}sentation
$$
\bigoplus_{\pi\in {\cal A}_{r}}\pi\otimes\sigma (\pi
)^{\vee}\otimes\sigma (\pi )
$$
du produit $\mathop{\rm GL}\nolimits_{r}^{}({\Bbb A})\times
\Gamma_{F}\times\Gamma_{F}$, o\`{u} $\sigma (\pi )^{\vee}$ est la
repr\'{e}sentation contragr\'{e}diente de $\sigma (\pi )$, que l'on
peut esp\'{e}rer obtenir comme la partie cuspidale de la cohomologie
$\ell$-adique \`{a} supports compacts d'un sch\'{e}ma sur $F\otimes
F$.

\vskip 7mm

\centerline{\bf 2. CHTOUCAS DE DRINFELD}
\vskip 10mm

{\bf 2.1. Le champ des chtoucas} 
\vskip 5mm

Tous les sch\'{e}mas (ou champs) consid\'{e}r\'{e}s seront sur ${\Bbb
F}_{q}$; on dira donc sch\'{e}ma (ou champ) au lieu de ${\Bbb
F}_{q}$-sch\'{e}ma (ou ${\Bbb F}_{q}$-champ) et on notera simplement
$U\times T$ le produit $U\times_{{\Bbb F}_{q}}T$ de deux tels
sch\'{e}mas (ou champs). Pour tout sch\'{e}ma (ou champ) $U$, on
notera $\mathop{\rm Frob}\nolimits_{U}:U\rightarrow U$ son
endomorphisme de Frobenius relativement \`{a} ${\Bbb F}_{q}$: si $U$
est un sch\'{e}ma, $\mathop{\rm Frob}\nolimits_{U}$ est donc
l'identit\'{e} sur l'espace topologique sous-jacent \`{a} $U$ et
l'\'{e}l\'{e}vation \`{a} la puissance $q$-i\`{e}me sur le faisceau
structural ${\cal O}_{U}$. Si ${\cal M}$ est un ${\cal O}_{U\times
X}$-Module, on notera ${}^{\tau}{\cal M}$ le ${\cal O}_{U\times
X}$-Module $(\mathop{\rm Frob}\nolimits_{U} \times\mathop{\rm
Id}\nolimits_{X})^{\ast}{\cal M}$.

\thm D\'{E}FINITION (Drinfeld)
\enonce
Un {\rm chtouca \`{a} droite} {\rm (}resp. {\rm \`{a} gauche}{\rm )}
$\widetilde{{\cal E}}$ de rang $r\geq 1$ sur un sch\'{e}ma $U$ est un
diagramme dans la cat\'{e}gorie ab\'{e}lienne des ${\cal O}_{U\times
X}$-Modules
$$
{\cal E}\,\smash{\mathop{ \lhook\joinrel\mathrel{\hbox to
6mm{\rightarrowfill}} }\limits^{\scriptstyle j}}\,{\cal E}'
\,\smash{\mathop{ {\hbox to 6mm{\leftarrowfill}}\joinrel\kern
-0.9mm\mathrel\rhook} \limits^{\scriptstyle t}}\,{}^{\tau}{\cal E}
\hbox{ {\rm (}resp. } {\cal E}\,\smash{\mathop{ {\hbox to
6mm{\leftarrowfill}}\joinrel\kern -0.9mm\mathrel\rhook}
\limits^{\scriptstyle t}}\,{\cal E}'\,\smash{\mathop{
\lhook\joinrel\mathrel{\hbox to 6mm{\rightarrowfill}}}
\limits^{\scriptstyle j}}\,{}^{\tau}{\cal E}\,)
$$
o\`{u}:

\itemitem{-} ${\cal E}$ et ${\cal E}'$ sont localement libres de rang
$r$,

\itemitem{-} $j$ et $t$ sont injectifs,

\itemitem{-} les conoyaux de $j$ et $t$ sont support\'{e}s par les
graphes $\Gamma_{\infty}\subset U\times X$ et $\Gamma_{o}\subset
U\times X$ de deux morphismes $\infty :U\rightarrow X$ et
$o:U\rightarrow X$, et sont localement libres de rang $1$ sur leurs
supports.

Le morphismes $\infty $ et $o$ sont appel\'{e}s respectivement le {\rm
p\^{o}le} et le {\rm z\'{e}ro} du chtouca.
\endthm

En d'autres termes, un chtouca \`{a} droite de rang $r$ est la
donn\'{e}e d'un fibr\'{e} vectoriel ${\cal E}$ de rang $r$, d'une
modification \'{e}l\'{e}mentaire sup\'{e}rieure $j:{\cal E}
\hookrightarrow {\cal E}'$ de ${\cal E}$, d'une modification
\'{e}l\'{e}mentaire inf\'{e}rieure $j':{\cal E}''\hookrightarrow {\cal
E}'$ de ${\cal E}'$ et d'un isomorphisme $u:{}^{\tau}{\cal E}
\buildrel\sim\over\longrightarrow {\cal E}''$. On a bien s\^{u}r une
description similaire pour les chtoucas \`{a} gauche.

En faisant varier $U$, on d\'{e}finit de mani\`{e}re \'{e}vidente le
champ $\mathop{\rm Cht}\nolimits^{r}$ des chtoucas \`{a} droite de
rang $r$, un morphisme $(\infty ,o):\mathop{\rm Cht}\nolimits^{r}
\rightarrow X\times X$ et un chtouca \`{a} droite de rang $r$
universel sur $\mathop{\rm Cht}\nolimits^{r}\times X$ de p\^{o}le et
de z\'{e}ro les deux composantes de ce morphisme. De m\^{e}me, on a le
champ ${}^{r}\!\mathop{\rm Cht}$ des chtoucas \`{a} gauche de rang
$r$, un morphisme $(\infty ,o):{}^{r}\!\mathop{\rm Cht} \rightarrow
X\times X$ et un chtouca \`{a} gauche de rang $r$ universel.

On a des morphismes de Frobenius partiels
$$
\mathop{\rm Frob}\nolimits_{\infty}:{}^{r}\!\mathop{\rm Cht}
\rightarrow\mathop{\rm Cht}\nolimits^{r}\hbox{ et }\mathop{\rm
Frob}\nolimits_{o}:\mathop{\rm Cht}\nolimits^{r} \rightarrow
{}^{r}\!\mathop{\rm Cht}
$$
au dessus de $(\infty ,o)\mapsto (\infty ,\mathop{\rm Frob}
\nolimits_{X}(o))$ et $(\infty ,o)\mapsto (\mathop{\rm Frob}
\nolimits_{X}(\infty ),o)$, qui envoient $\widetilde{{\cal E}}$ sur
$({\cal E}'\,\smash{\mathop{ \lhook\joinrel\mathrel{\hbox to
6mm{\rightarrowfill}} }\limits^{\scriptstyle j}}\,{}^{\tau}{\cal E}
\,\smash{\mathop{ {\hbox to 6mm{\leftarrowfill}}\joinrel\kern
-0.9mm\mathrel\rhook} \limits^{\scriptstyle {}^{\tau}\!t}}\,
{}^{\tau}{\cal E}')$ et $({\cal E}'\,\smash{\mathop{ {\hbox to
6mm{\leftarrowfill}}\joinrel\kern -0.9mm\mathrel\rhook}
\limits^{\scriptstyle t}}\,{}^{\tau}{\cal E}\,\smash{\mathop{
\lhook\joinrel\mathrel{\hbox to 6mm{\rightarrowfill}}}
\limits^{\scriptstyle {}^{\tau}\!j}}\,{}^{\tau}{\cal E}')$
respectivement.
\vskip 3mm

{\it Dans la suite et sauf mention explicite du contraire, les
chtoucas seront tous \`{a} droite, et on dira simplement
{\og}{chtouca}{\fg} pour {\og}{chtouca \`{a} droite}{\fg}.} Les
r\'{e}sultats d\'{e}montr\'{e}s pour les chtoucas \`{a} droite ont
bien entendu des analogues pour les chtoucas \`{a} gauche.

\thm PROPOSITION (Drinfeld)
\enonce
Le champ $\mathop{\rm Cht}\nolimits^{r}$ est alg\'{e}brique au sens de
Deligne-Mumford. Le morphisme $(\infty ,o):\mathop{\rm Cht}
\nolimits^{r}\rightarrow X\times X$ est lisse, purement de dimension
relative $2r-2$.
\endthm

Le champ $\mathop{\rm Cht}\nolimits^{r}$ s'\'{e}crit comme r\'{e}union
disjointe
$$
\mathop{\rm Cht}\nolimits^{r}=\coprod_{d\in {\Bbb Z}}^{}\mathop{\rm
Cht}\nolimits^{r,d}
$$
o\`{u} $\mathop{\rm Cht}\nolimits^{r,d}$ classifie les chtoucas de
degr\'{e}
$$
d=\mathop{\rm deg}{\cal E}=\mathop{\rm deg}{\cal E}'-1.
$$

\rem Exemple 
\endrem
Pour tout entier $d$, $\mathop{\rm Cht}\nolimits^{1,d}$ peut \^{e}tre
d\'{e}fini comme le produit fibr\'{e}
$$\diagram{
\mathop{\rm Cht}\nolimits^{1,d}&\kern -1mm\smash{\mathop{\hbox to
8mm{\rightarrowfill}}\limits^{\scriptstyle }}\kern -1mm&\mathop{\rm
Fib}\nolimits^{1,d}\cr
\llap{$\scriptstyle $}\left\downarrow\vbox to 4mm{}\right.\rlap{}
&\square&\llap{}\left\downarrow \vbox to 4mm{}\right.
\rlap{$\scriptstyle L$}\cr
X\times X&\kern -1mm\smash{\mathop{\hbox to 8mm{\rightarrowfill}}
\limits_{\scriptstyle J}}\kern -1mm&\mathop{\rm Fib}\nolimits^{1,0}\cr}
$$
o\`{u} $\mathop{\rm Fib}\nolimits^{1,d}$ est le champ alg\'{e}brique
des fibr\'{e}s en droites de degr\'{e} $d$ sur $X$, $J$ est
l'application d'Abel-Jacobi qui envoie $(\infty ,o)\in (X\times X)(U)$
sur le fibr\'{e} en droites ${\cal O}_{U\times X}(\infty -o)$, et $L$
est l'isog\'{e}nie de Lang qui envoie un fibr\'{e} en droites ${\cal
L}$ sur $U\times X$ sur le fibr\'{e} en droites ${\cal L}^{-1}
\otimes_{{\cal O}_{U\times X}}{}^{\tau}{\cal L}$.

En particulier, $\mathop{\rm Cht}\nolimits^{1,d}$ est de type fini et
admet un espace grossier qui est un rev\^{e}tement fini \'{e}tale
galoisien de $X\times X$.
\hfill\hfill$\square$
\vskip 3mm

{\it Mais, mis \`{a} part le cas $r=1$, aucun des champs
$\mathop{\rm Cht}\nolimits^{r,d}$ n'est de type fini.}
\vskip 5mm

{\bf 2.2. Troncatures}
\vskip 5mm

Soient $k$ un corps alg\'{e}briquement clos contenant ${\Bbb F}_{q}$
et $\widetilde{{\cal E}}$ un chtouca sur (le spectre de) $k$. On
appelle {\it sous-objet} de $\widetilde{{\cal E}}$, et on note
simplement
$$
\widetilde{{\cal F}}\subset\widetilde{{\cal E}},
$$
la donn\'{e}e de deux sous-${\cal O}_{k\otimes X}$-Modules ${\cal
F}\subset {\cal E}$ et ${\cal F}'\subset {\cal E}'$ tels que ${\cal
E}/{\cal F}$ et ${\cal E}'/{\cal F}'$ soient localement libres de
m\^{e}me rang et que $j({\cal F})\subset {\cal F}'$ et
$t({}^{\tau}{\cal F})\subset {\cal F}'$; \`{a} un tel sous-objet on
peut associer son {\it rang}
$$
\mathop{\rm rg}\widetilde{{\cal F}}=\mathop{\rm rg}{\cal 
F}=\mathop{\rm rg}{\cal F}'
$$
et, pour chaque nombre r\'{e}el $\alpha$, son degr\'{e} (d'indice
$\alpha$)
$$
\mathop{\rm deg}\nolimits_{\alpha}\widetilde{{\cal F}}=(1-\alpha 
)\mathop{\rm deg}{\cal F}+\alpha\mathop{\rm deg}{\cal F}'.
$$

Si ${0}=\widetilde{{\cal F}}_{0}\subsetneq\widetilde{{\cal F}}_{1}
\subsetneq\cdots\subsetneq\widetilde{{\cal F}}_{s}=\widetilde{{\cal
E}}$ est une filtration de $\widetilde{{\cal E}}$ par des sous-objets
comme ci-dessus, on peut lui associer son polygone (d'indice $\alpha$)
qui est la fonction affine par morceaux
$$
p:[0,r]\rightarrow {\Bbb R}
$$
avec $p(0)=p(r)=0$, dont les seules ruptures de pentes interviennent
en les entiers $\mathop{\rm rg}\widetilde{{\cal F}}_{\sigma}$, $\sigma
=1,\ldots ,s-1$, et qui prend en ces entiers-l\`{a} la valeur
$$
p(\mathop{\rm rg}\widetilde{{\cal F}}_{\sigma})=\mathop{\rm
deg}\nolimits_{\alpha}\widetilde{{\cal F}}_{\sigma}-{\mathop{\rm
rg}\widetilde{{\cal F}}_{\sigma}\over r}\mathop{\rm
deg}\nolimits_{\alpha}\widetilde{{\cal E}}.
$$

\thm PROPOSITION 
\enonce
Fixons un nombre r\'{e}el $\alpha\in [0,1]$. Alors, parmi tous les
polygones attach\'{e}s aux filtrations de $\widetilde{{\cal E}}$ comme
ci-dessus, il en existe un plus grand que tous les autres, et parmi
toutes les filtrations qui d\'{e}finissent ce polygone maximal il en
existe une moins fine que tous les autres.
\endthm

Le polygone et la filtration dont la proposition ci-dessus assurent
l'existence sont appel\'{e}s respectivement le {\it polygone de
Harder-Narasimhan} et la {\it filtration de Harder-Narasimhan}
d'indice $\alpha$ du chtouca $\widetilde{{\cal E}}$.

Nous appellerons {\it param\`{e}tre de troncature} toute fonction
continue, convexe, affine par morceaux $p:[0,r]\rightarrow {\Bbb
R}_{\geq 0}$ avec $p(0)=p(r)=0$ dont les points de ruptures de pente
ont des abscisses enti\`{e}res.

\thm PROPOSITION 
\enonce
\'{E}tant donn\'{e}s $\alpha\in [0,1]$ et un param\`{e}tre de
troncature $p$, il existe un unique ouvert $\mathop{\rm Cht}
\nolimits^{r;\,\leq_{\alpha} p}$ du champ $\mathop{\rm Cht}
\nolimits^{r}$ tel qu'un chtouca sur un corps alg\'{e}briquement clos
est dans cet ouvert si et seulement si son polygone de
Harder-Narasimhan d'indice $\alpha$ est major\'{e} par $p$.

De plus, pour chaque entier $d$, l'ouvert
$$
\mathop{\rm Cht}\nolimits^{r,d;\,\leq_{\alpha} p}:=\mathop{\rm
Cht}\nolimits^{r;\,\leq_{\alpha} p} \cap\mathop{\rm
Cht}\nolimits^{r,d}
$$
du champ $\mathop{\rm Cht}\nolimits^{r,d}$ est de type fini.
\endthm

On voit donc que, pour chaque entier $d$, le champ alg\'{e}brique
localement de type fini $\mathop{\rm Cht}\nolimits^{r,d}$ est
r\'{e}union filtrante des ouverts de type fini $\mathop{\rm
Cht}\nolimits_{\alpha}^{r,d;\,\leq p}$.
\vskip 5mm

{\bf 2.3. Structures de niveau sur les chtoucas} 
\vskip 5mm

Comme les courbes elliptiques, les chtoucas admettent des
structures de niveau. Un {\it niveau} est ici un sous-sch\'{e}ma
ferm\'{e} fini
$$
N=\mathop{\rm Spec}({\cal O}_{N})\subset X.
$$
Si ${\cal E}$ est un fibr\'{e} vectoriel sur $U\times X$ pour un
sch\'{e}ma (ou un champ) $U$, on notera simplement ${\cal E}_{N}$ la
restriction de ${\cal E}$ au ferm\'{e} $U\times N\subset U\times X$.
Si ${\cal M}$ est un ${\cal O}_{U\times N}$-Module on notera encore
${}^{\tau}{\cal M}$ le ${\cal O}_{U\times N}$-Module $(\mathop{\rm
Frob}\nolimits_{U}\times \mathop{\rm Id}_{N})^{\ast}{\cal M}$. Pour
tout fibr\'{e} vectoriel ${\cal E}$ sur $U\times X$ on a
\'{e}videmment ${}^{\tau}({\cal E}_{N})\cong ({}^{\tau}{\cal E})_{N}$.

On consid\`{e}re le champ alg\'{e}brique
$$
\mathop{\rm Tr}\nolimits_{N}^{r}
$$
qui associe au sch\'{e}ma $U$ la cat\'{e}gorie des couples $({\cal
F},u)$ form\'{e}s d'un fibr\'{e} vectoriel de rang $r$ sur $U\times N$
et d'un isomorphisme de fibr\'{e}s vectoriels $u:{}^{\tau}{\cal
F}\buildrel\sim\over \longrightarrow {\cal F}$. Si l'on note
$\mathop{\rm GL}\nolimits_{r}^{N}=\mathop{\rm Res}\nolimits_{N/
\mathop{\rm Spec}({\Bbb F}_{q})}\mathop{\rm GL}\nolimits_{r}$ le
sch\'{e}ma en groupes sur ${\Bbb F}_{q}$ restriction \`{a} la Weil du
sch\'{e}ma en groupes lin\'{e}aire $\mathop{\rm GL}\nolimits_{r}$ sur
$N$ et par $h\mapsto \tau (h)=\mathop{\rm Frob}\nolimits_{\mathop{\rm
GL} \nolimits_{r}}(h)$ son endomorphisme de Frobenius (relatif \`{a}
${\Bbb F}_{q}$), ce champ n'est autre que le quotient 
$$
\mathop{\rm Tr}\nolimits_{N}^{r}=[\mathop{\rm GL}\nolimits_{r}^{N}/
{}^{\tau}\!\!\mathop{\rm Ad}(\mathop{\rm GL}\nolimits_{r}^{N})]
$$
de $\mathop{\rm GL}\nolimits_{r}^{N}$ par l'action de $\mathop{\rm
GL}\nolimits_{r}^{N}$ sur lui-m\^{e}me par la conjugaison tordue
$$
{}^{\tau}\!\!\mathop{\rm Ad}(h)(g)=\tau (h)^{-1}gh.
$$
Comme l'{\it isog\'{e}nie de Lang}
$$
L:\mathop{\rm GL}\nolimits_{r}^{N}\rightarrow \mathop{\rm 
GL}\nolimits_{r}^{N},~h\mapsto L(h)=\tau (h)^{-1}h,
$$
est un torseur sous le groupe fini $\mathop{\rm GL}\nolimits_{r}^{}
({\cal O}_{N})=\mathop{\rm GL}\nolimits_{r}^{N}({\Bbb F}_{q})$ (pour
son action par translation \`{a} droite sur $\mathop{\rm GL}
\nolimits_{r}^{N}$), $\mathop{\rm Tr}\nolimits_{N}^{r}$ est
canoniquement isomorphe au champ classifiant $[\mathop{\rm Spec}
({\Bbb F}_{q})/\mathop{\rm GL}\nolimits_{r}^{}({\cal O}_{N})]$.

{\it On ne munira ici de structures de niveau que les chtoucas dont le
p\^{o}le et le z\'{e}ro ne rencontrent pas $N$}. Si
$$
\widetilde{{\cal E}}=({\cal E}\,\smash{\mathop{\lhook\joinrel
\mathrel{\hbox to 8mm {\rightarrowfill}}}\limits^{\scriptstyle
j}}\,{\cal E}'\, \smash{\mathop{ {\hbox to 8mm{\leftarrowfill}}
\joinrel\kern -0.9mm \mathrel\rhook} \limits^{\scriptstyle
t}}\,{}^{\tau}{\cal E})
$$
est un tel chtouca sur $U$ on remarque que les restrictions
$$
j_{N}:{\cal E}_{N}\rightarrow {\cal E}_{N}'\hbox{ et }t_{N}:
{}^{\tau}({\cal E}_{N})\cong ({}^{\tau}{\cal E})_{N}\rightarrow {\cal
E}_{N}'
$$
des fl\`{e}ches $j$ et $t$ sont des isomorphismes de ${\cal
O}_{U\times N}$-Modules. On a donc un morphisme de champs
alg\'{e}briques
$$
\mathop{\rm Cht}\nolimits^{r}\times_{X^{2}}^{}(X-N)^{2}\,
\smash{\mathop{\hbox to 8mm{\rightarrowfill}}\limits_{\scriptstyle
}}\,\mathop{\rm Tr}\nolimits_{N}^{r},~\widetilde{{\cal E}}\mapsto
({\cal E}_{N},j_{N}^{-1} \circ t_{N}).
$$

\thm D\'{E}FINITION (Drinfeld)
\enonce
Soit $\widetilde{{\cal E}}$ un chtouca de rang $r$ sur un sch\'{e}ma
$U$ dont le p\^{o}le et le z\'{e}ro se factorisent par l'immersion
ouverte $X-N\hookrightarrow X$. Une {\rm structure (principale) de
niveau} $N$ sur ce chtouca est la donn\'{e}e d'un isomorphisme
$$
\iota :{\cal O}_{U\times N}^{r}\buildrel\sim\over\longrightarrow {\cal
E}_{N}
$$
tel que
$$
\iota =(j_{N}^{-1}\circ t_{N})\circ {}^{\tau}\iota :{\cal O}_{U\times
N}^{r}\buildrel\sim\over\longrightarrow {\cal E}_{N}.
$$
\endthm

On d\'{e}finit de mani\`{e}re \'{e}vidente le champ $\mathop{\rm
Cht}\nolimits_{N}^{r}$ des chtoucas de rang $r$, dont le p\^{o}le et
le z\'{e}ro ne rencontrent pas $N$, et avec structures de niveau $N$.
On a un carr\'{e} cart\'{e}sien
$$\diagram{
\mathop{\rm Cht}\nolimits_{N}^{r}&\kern -12mm\smash{\mathop{\hbox to
19mm {\rightarrowfill}}\limits^{\scriptstyle }}\kern -1mm&\mathop{\rm 
Tr}\nolimits_{N}^{r,\tau}\cr
\llap{$\scriptstyle $}\left\downarrow
\vbox to 4mm{}\right.\rlap{}&\square&\llap{}\left\downarrow
\vbox to 4mm{}\right.\rlap{$\scriptstyle $}\cr
\mathop{\rm Cht}\nolimits^{r}\times_{X^{2}}^{}(X-N)^{2}&\kern
-1mm\smash{\mathop{\hbox to 8mm{\rightarrowfill}}
\limits_{\scriptstyle }}\kern -1mm&\mathop{\rm Tr}\nolimits_{N}^{r}\cr}
$$
o\`{u} la fl\`{e}che horizontale du bas est la fl\`{e}che de
restriction \`{a} $N$ d\'{e}finie ci-dessus, o\`{u} $\mathop{\rm Tr}
\nolimits_{N}^{r,\tau}$ est le champ qui associe au sch\'{e}ma $U$ la
cat\'{e}gorie des triplets $({\cal F},u,\iota )$ form\'{e}s d'un objet
$({\cal F},u)$ de $\mathop{\rm Tr}\nolimits_{N}^{r}(U)$ et d'un
isomorphisme de ${\cal O}_{U\times N}$-Modules $\iota :{\cal
O}_{U\times N}^{r}\buildrel\sim\over\longrightarrow {\cal F}$ tel que
$\iota =u\circ {}^{\tau}\iota$, et o\`{u} la fl\`{e}che verticale de
droite est l'oubli de la composante $\iota$. On remarquera que la
composante $\iota$ de $({\cal F},u,\iota )$ d\'{e}termine uniquement
la composante $u$ et donc que $\mathop{\rm Tr}\nolimits_{N}^{r,\tau}$
est canoniquement isomorphe au champ $[\mathop{\rm
GL}\nolimits_{r}^{N}/\mathop{\rm GL}\nolimits_{r}^{N}]$ quotient de
$\mathop{\rm GL}\nolimits_{r}^{N}$ par l'action de $\mathop{\rm GL}
\nolimits_{r}^{N}$ sur lui-m\^{e}me par translation \`{a} droite,
c'est-\`{a}-dire au sch\'{e}ma r\'{e}duit \`{a} un point $\mathop{\rm
Spec}({\Bbb F}_{q})$. Par cons\'{e}quent, la fl\`{e}che verticale de
droite du carr\'{e} ci-dessus se d\'{e}crit encore, soit comme la
fl\`{e}che $[\mathop{\rm GL}\nolimits_{r}^{N}/\mathop{\rm GL}
\nolimits_{r}^{N}] \rightarrow [\mathop{\rm GL}\nolimits_{r}^{N}
/{}^{\tau}\!\!\mathop{\rm Ad} (\mathop{\rm GL}\nolimits_{r}^{N})]$
induite par l'isog\'{e}nie de Lang $L$, soit comme le $\mathop{\rm GL}
\nolimits_{r}^{}({\cal O}_{N})$-torseur canonique $\mathop{\rm Spec}
({\Bbb F}_{q})\rightarrow [\mathop{\rm Spec}({\Bbb F}_{q})/
\mathop{\rm GL}\nolimits_{r}^{}({\cal O}_{N})]$.

La fl\`{e}che verticale de gauche est donc un $\mathop{\rm GL}
\nolimits_{r}^{}({\cal O}_{N})$-torseur et $\mathop{\rm Cht}
\nolimits_{N}^{r}$ est un champ alg\'{e}brique, lisse purement de
dimension $2r-2$ sur $(X-N)^{2}$.

\vskip 7mm

\centerline{\bf 3. HOMOMORPHISMES COMPLETS ET CHTOUCAS IT\'{E}R\'{E}S}
\vskip 10mm

{\bf 3.1. Homomorphismes complets}
\vskip 5mm

Soient $k$ un corps et $V$ et $W$ deux $k$-espaces vectoriels de
dimension $r$. On notera $H(V,W)$ le sch\'{e}ma des uplets
$$
(u_{1},\ldots ,u_{r};\lambda_{1},\ldots ,\lambda_{r-1})
$$
o\`{u}, pour $\rho =1,\ldots ,r$, $u_{\rho}$ est une application
lin\'{e}aire non nulle de $\wedge^{\rho}V$ dans $\wedge^{\rho}W$ et
o\`{u} $\lambda_{1},\ldots ,\lambda_{r-1}$ sont des scalaires. Les
conditions
\vskip 2mm

\itemitem{-} $u_{1}$ est un isomorphisme,
\vskip 2mm

\itemitem{-} $\lambda_{1},\ldots ,\lambda_{r-1}$ sont inversibles,
\vskip 2mm

\itemitem{-} $\wedge^{\rho}u_{1}=\lambda_{1}^{\rho -1}
\lambda_{2}^{\rho -2}\cdots \lambda_{\rho -1}u_{\rho}$ pour $\rho
=2,\ldots ,r$,
\vskip 2mm

\noindent d\'{e}finissent un sous-sch\'{e}ma localement ferm\'{e}
$\Omega^{\circ}(V,W)$ de ce sch\'{e}ma affine et la projection
$(u_{1}; \lambda_{1},\ldots ,\lambda_{r-1})$ identifie
$\Omega^{\circ}(V,W)$ au sch\'{e}ma $\mathop{\rm Isom}(V,W)\times_{k}
{\Bbb G}_{{\rm m},k}^{r-1}$. Nous noterons
$$
\Omega (V,W)\subset H(V,W)
$$
l'adh\'{e}rence sch\'{e}matique de $\Omega^{\circ}(V,W)$ dans
$H(V,W)$.

Le sch\'{e}ma ${\Bbb A}_{k}^{r-1}$ est naturellement muni d'un
diviseur \`{a} croisements normaux, r\'{e}union de $r-1$ diviseurs
lisses, \`{a} savoir les diviseurs $\{\lambda_{\rho}=0\}$ pour $\rho
=1,\ldots ,r-1$. Ce diviseur \`{a} croisement normaux d\'{e}finit de
la mani\`{e}re habituelle une stratification localement ferm\'{e}e de
${\Bbb A}_{k}^{r-1}$, index\'{e}e par les parties $R$ de
$[r-1]=\{1,\ldots ,r-1\}$ et de $R$-i\`{e}me strate
$$
{\Bbb A}_{k,R}^{r-1}=\{(\lambda_{1},\ldots ,\lambda_{r-1})\mid
\lambda_{\rho}=0\Leftrightarrow \rho\in R\}\cong {\Bbb G}_{{\rm
m},k}^{[r-1]-R}.
$$

Par construction, on dispose d'un morphisme
$$
(\lambda_{1},\ldots ,\lambda_{r-1}):\Omega (V,W)\rightarrow {\Bbb
A}_{k}^{r-1}.
$$
et on peut relever \`{a} $\Omega (V,W)$ le diviseur et la
stratification de ${\Bbb A}_{k}^{r-1}$ que l'on vient d'introduire. En
particulier, pour tout sous-ensemble $R\subset [r-1]$, on notera
$$
\Omega_{R}(V,W)\subset\Omega (V,W)
$$
l'image r\'{e}ciproque par ce morphisme de la strate localement
ferm\'{e}e ${\Bbb A}_{k,R}^{r-1}\subset {\Bbb A}_{k}^{r-1}$. On
v\'{e}rifie que $\Omega_{\emptyset}(V,W)=\Omega^{\circ}(V,W)$.

\thm PROPOSITION
\enonce
Pour chaque $R=\{r_{1},r_{2},\ldots ,r_{s-1}\}\subset [r-1]$ o\`{u}
$0=r_{0}<r_{1}<r_{2}<\cdots <r_{s-1}<r_{s}=r$, $\Omega_{R}(V,W)$ est
isomorphe au sch\'{e}ma des uplets
$$
(V^{\bullet},W_{\bullet},(v_{\sigma})_{\sigma =1,\ldots ,s-1};
(\lambda_{\rho})_{\rho\in [r-1]-R})
$$
o\`{u}:
\vskip 2mm

\itemitem{-} les $\lambda_{\rho}$ sont des scalaires inversibles,
\vskip 2mm

\itemitem{-} $V^{\bullet}=(V=V^{0}\supsetneq V^{1}\supsetneq\cdots
\supsetneq V^{s}=(0))$ est une filtration d\'{e}croissante par des
sous-espaces vectoriels de codimensions $0=r_{0},r_{1},\ldots
,r_{s}=r$,
\vskip 2mm

\itemitem{-} $W_{\bullet}=((0)=W_{0}\subsetneq W_{1}\subsetneq\cdots
\subsetneq W_{s}=W)$ est une filtration croissante par des
sous-espaces vectoriels de dimensions $0=r_{0},r_{1},\ldots ,r_{s}=r$,
\vskip 2mm

\itemitem{-} $v_{\sigma}:V^{\sigma -1}/V^{\sigma}\buildrel\sim\over
\longrightarrow W_{\sigma}/W_{\sigma -1}$ est un isomorphisme
d'espaces vectoriels.

A fortiori, le morphisme $(\lambda_{\rho})_{\rho\in [r-1]-R}:
\Omega_{R}(V,W)\rightarrow {\Bbb G}_{{\rm m},k}^{[r-1]-R}$ est lisse
de dimension relative $r^{2}$.
\hfill\hfill$\square$
\endthm

\rem Remarque
\endrem
L'\'{e}criture $R=\{r_{1},r_{2},\ldots ,r_{s-1}\}$ o\`{u}
$0=r_{0}<r_{1}<r_{2}<\cdots <r_{s-1}<r_{s}=r$ identifie les parties de
$[r-1]$ aux d\'{e}compositions $r=(r_{1}-r_{0})+(r_{2}-r_{1})+\cdots
+(r_{s}-r_{s-1})$ de $r$ en entiers strictement positifs.
\hfill\hfill$\square$
\vskip 3mm

\thm COROLLAIRE 
\enonce
Le morphisme $(\lambda_{1},\ldots ,\lambda_{r-1})$ est lisse
purement de dimension relative $r^{2}$ et $\Omega^{\circ}(V,W)\subset 
\Omega (V,W)$ est l'ouvert compl\'{e}mentaire d'un diviseur \`{a} 
croisements normaux, r\'{e}union des $r-1$ diviseurs lisses 
$\{\lambda_{\rho}=0\}$ pour $\rho =1,\ldots ,r-1$.
\endthm

Le tore ${\Bbb G}_{{\rm m},k}^{r-1}=\{(\mu_{1},\ldots ,\mu_{r-1})\}$
agit librement sur $\Omega (V,W)$ par
$$
u_{1}\mapsto u_{1},~u_{2}\mapsto \mu_{1}^{-1}u_{2},~u_{3}\mapsto
\mu_{1}^{-2}\mu_{2}u_{3},~\ldots ,~u_{r}\mapsto
\mu_{1}^{1-r}\mu_{2}^{2-r}\cdots\mu_{r-1}u_{r},
$$
et
$$
\lambda_{1}\mapsto\mu_{1}\lambda_{1},~
\lambda_{2}\mapsto\mu_{2}\lambda_{2},~\ldots ,~
\lambda_{r-1}\mapsto\mu_{r-1}\lambda_{r-1},
$$
et le quotient 
$$
\mathop{\rm H\widetilde{om}}(V,W):=\Omega (V,W)/{\Bbb G}_{{\rm
m},k}^{r-1}
$$
est un sch\'{e}ma quasi-projectif et lisse, qui contient comme ouvert
dense $\mathop{\rm Isom}(V,W)$ avec pour ferm\'{e} compl\'{e}mentaire
une r\'{e}union de $r-1$ diviseurs lisses \`{a} croisements normaux.
C'est par d\'{e}finition le {\it sch\'{e}ma des homomorphismes
complets de} $V$ dans $W$. Par construction, on a un morphisme
repr\'{e}sentable, lisse et purement de dimension relative $r^{2}$, de
champs alg\'{e}briques
$$
\mathop{\rm H\widetilde{om}}(V,W)\rightarrow [{\Bbb A}_{k}^{r-1}/{\Bbb
G}_{{\rm m},k}^{r-1}],
$$
et le diviseur \`{a} croisements normaux ci-dessus est l'image
r\'{e}ciproque par ce morphisme du diviseur \`{a} croisements normaux
\'{e}vident sur le champ alg\'{e}brique $[{\Bbb A}_{k}^{r-1}/{\Bbb
G}_{{\rm m},k}^{r-1}]$. 

Si $V=W=k^{r}$ on notera encore $\mathop{\widetilde{\rm gl}}
\nolimits_{r,k}$ le $k$-sch\'{e}ma $\mathop{\rm H\widetilde{om}}
(k^{r},k^{r})$ des {\it endomorphismes complets} de $k^{r}$. La strate
ouverte de $\mathop{\widetilde{\rm gl}}\nolimits_{r,k}$ est bien
entendu $\mathop{\rm GL}\nolimits_{r,k}$.

\rem Remarques 
\endrem
{\rm (i)} Le morphisme
$$
u_{1}:\mathop{\rm H\widetilde{om}}(V,W)\rightarrow \mathop{\rm 
Hom}(V,W)-\{0\}
$$
peut aussi \^{e}tre obtenu comme le compos\'{e} 
$$
\mathop{\rm H\widetilde{om}}(V,W)=H_{r-1}\rightarrow\cdots\rightarrow
H_{1}\rightarrow H_{0}=\mathop{\rm Hom}(V,W)-\{0\}
$$
o\`{u} $H_{1}\rightarrow H_{0}$ est l'\'{e}clatement de $H_{0}$ le
long du ferm\'{e} des homomorphismes de rang $1$ et o\`{u}, pour $\rho
=2,\ldots ,r-1$, $H_{\rho}\rightarrow H_{\rho -1}$ est
l'\'{e}clatement de $H_{\rho -1}$ le long du transform\'{e} strict du
ferm\'{e} de $H_{0}$ form\'{e} des homomorphismes de rang $\leq\rho$.

\decale{\rm (ii)} Le groupe multiplicatif agit par homoth\'{e}tie sur
les homomorphismes complets. Le quotient $\mathop{\widetilde{\rm
gl}}\nolimits_{r,k}/{\Bbb G}_{{\rm m},k}$ est un sch\'{e}ma projectif,
isomorphe \`{a} la compactification de De Concini et Procesi de
$\mathop{\rm PGL}\nolimits_{r,k}$.
\hfill\hfill$\square$
\vskip 3mm

Pour tout $k$-sch\'{e}ma $U$, la cat\'{e}gorie $[{\Bbb A}_{k}^{r-1}/
{\Bbb G}_{{\rm m},k}^{r-1}](U)$ a pour objets les uplets
$$
(({\cal L}_{1},\lambda_{1}),\ldots ,({\cal L}_{r-1},\lambda_{r-1}))
$$
o\`{u} ${\cal L}_{1},\ldots ,{\cal L}_{r-1}$ sont des ${\cal
O}_{U}$-Modules inversibles et o\`{u} $\lambda_{1},\ldots
,\lambda_{r-1}$ sont des sections globales de ces fibr\'{e}s en
droites. On peut donc voir un $U$-homomorphisme complet de $V$ dans
$W$, c'est-\`{a}-dire un $U$-point de $\mathop{\rm H\widetilde{om}}
(V,W)$, comme un uplet
$$
u=(u_{1},\ldots ,u_{r};({\cal L}_{1},\lambda_{1}),\ldots ,({\cal
L}_{r-1},\lambda_{r-1}))
$$
o\`{u} $(({\cal L}_{1},\lambda_{1}),\ldots ,({\cal L}_{r-1},
\lambda_{r-1}))\in\mathop{\rm ob}\,[{\Bbb A}_{k}^{r-1}/{\Bbb G}_{{\rm
m},k}^{r-1}](U)$ et o\`{u}
$$
u_{\rho}:\bigwedge^{\rho}V\otimes_{k}{\cal L}_{1}^{\otimes (\rho -1)}
\otimes {\cal L}_{2}^{\otimes (\rho -2)}\otimes\cdots\otimes {\cal
L}_{\rho -1}\rightarrow \bigwedge^{\rho}W
$$
est un homomorphisme partout non nul de ${\cal O}_{U}$-Modules pour
$\rho =1,\ldots ,r$. Un homomorphisme complet $u$ est donc bien un
homomorphisme $u_{1}$ {\it compl\'{e}t\'{e}} par des donn\'{e}es
suppl\'{e}mentaires et peut \^{e}tre not\'{e} commod\'{e}ment
$u:V\Rightarrow W$.

Soient maintenant $S$ un $k$-sch\'{e}ma et ${\cal V}$, ${\cal W}$ deux
${\cal O}_{S}$-modules localement libres de rang constant $r$. Si $U$
est un $S$-sch\'{e}ma, on peut consid\'{e}rer plus
g\'{e}n\'{e}ralement les uplets
$$
u=(u_{1},\ldots ,u_{r};({\cal L}_{1},\lambda_{1}),\ldots ,({\cal
L}_{r-1},\lambda_{r-1}))
$$
o\`{u} $(({\cal L}_{1},\lambda_{1}),\ldots ,({\cal L}_{r-1},
\lambda_{r-1}))\in\mathop{\rm ob}\,[{\Bbb A}_{k}^{r-1}/{\Bbb G}_{{\rm
m},k}^{r-1}](U)$ et o\`{u}
$$
u_{\rho}:\bigwedge^{\rho}{\cal V}_{U}\otimes {\cal L}_{1}^{\otimes (\rho -1)}
\otimes {\cal L}_{2}^{\otimes (\rho -2)}\otimes\cdots\otimes {\cal
L}_{\rho -1}\rightarrow \bigwedge^{\rho}{\cal W}_{U},
$$
est un homomorphisme partout non nul de ${\cal O}_{U}$-Modules pour
$\rho =1,\ldots ,r$. (On a not\'{e} ${\cal V}_{U}$ et ${\cal W}_{U}$
les restrictions de ${\cal V}$ et ${\cal W}$ \`{a} $U$.) On dira alors
qu'un tel uplet est un {\it $U$-homomorphisme complet de ${\cal V}$
dans ${\cal W}$}, et on utilisera la notation $u:{\cal V}\Rightarrow
{\cal W}$, si, pour toutes trivialisations locales ${\cal V}\cong
V\otimes_{k}{\cal O}_{S}$ et ${\cal W}\cong W\otimes_{k}{\cal O}_{S}$,
ce uplet est un $U$-homomorphisme complet de $V$ dans $W$.

On d\'{e}finit ainsi un $S$-sch\'{e}ma $\mathop{\rm H\widetilde{om}}
({\cal V},{\cal W})$, muni d'un morphisme repr\'{e}sentable, lisse et
purement de dimension relative $r^{2}$, de champs alg\'{e}briques
$$
\mathop{\rm H\widetilde{om}}({\cal V},{\cal W})\rightarrow S\times_{k}
[{\Bbb A}_{k}^{r-1}/{\Bbb G}_{{\rm m},k}^{r-1}]
$$
Le $S$-sch\'{e}ma quasi-projectif $\mathop{\rm H\widetilde{om}}({\cal
V},{\cal W})$ contient le $S$-sch\'{e}ma $\mathop{\rm Isom}({\cal
V},{\cal W})$ comme un ouvert dense, et le ferm\'{e}
compl\'{e}mentaire est un diviseur \`{a} croisement normaux relatifs
sur $S$, r\'{e}union de $r-1$ diviseurs lisses. 
\vskip 5mm

{\bf 3.2.  Chtoucas it\'{e}r\'{e}s}
\vskip 5mm

\thm D\'{E}FINITION 
\enonce
Un {\rm pr\'{e}-chtouca it\'{e}r\'{e}} de rang $r$ sur un sch\'{e}ma
$U$ est la donn\'{e}e:
\vskip 2mm

\itemitem{-} d'un diagramme
$$
{\cal E}\,\smash{\mathop{ \lhook\joinrel\mathrel{\hbox to
6mm{\rightarrowfill}} }\limits^{\scriptstyle j}}\,{\cal
E}'\,\smash{\mathop{ {\hbox to 6mm{\leftarrowfill}}\joinrel\kern
-0.9mm\mathrel\rhook} \limits^{\scriptstyle j'}}\,{\cal E}''
$$ 
et de deux morphismes $\infty :U\rightarrow X$ et $o:U\rightarrow X$,
o\`{u} ${\cal E}$ est un fibr\'{e} vectoriel de rang $r$ sur $U\times
X$, $j$ est une modification \'{e}l\'{e}mentaire sup\'{e}rieure de
${\cal E}$ le long du graphe de $\infty$ et $j'$ est une modification
\'{e}l\'{e}mentaire inf\'{e}rieure de ${\cal E}'$ le long du graphe de
$o$,
\vskip 2mm

\itemitem{-} de ${\cal O}_{U}$-Modules inversibles ${\cal 
L}_{1},\ldots ,{\cal L}_{r-1}$ munis de sections globales 
$\lambda_{1},\ldots ,\lambda_{r-1}$,
\vskip 2mm

\itemitem{-} pour chaque entier $\rho =1,\ldots ,r$, d'un
homomorphisme de ${\cal O}_{U\times X}$-Modules
$$
u_{\rho}:\wedge^{\rho}({}^{\tau}{\cal E})\otimes{\cal L}_{1}^{\otimes
(q-1)(\rho -1)}\otimes{\cal L}_{2}^{\otimes (q-1)(\rho
-2)}\otimes\cdots\otimes{\cal L}_{\rho -1}^{\otimes
(q-1)}\rightarrow\wedge^{\rho}{\cal E}'',
$$
\vskip 2mm

\noindent de telle sorte que le uplet
$$
u=(u_{1},\ldots ,u_{r};({\cal L}_{1}^{\otimes (q-1)},
\lambda_{1}^{q-1}), \ldots ,({\cal L}_{r-1}^{\otimes (q-1)},
\lambda_{r-1}^{q-1})):{}^{\tau}{\cal E}\Rightarrow {\cal E}''
$$
soit un homomorphisme complet.
\endthm

En faisant varier $U$, on d\'{e}finit de mani\`{e}re \'{e}vidente le
champ ${\cal C}^{r}$ des pr\'{e}-chtoucas it\'{e}r\'{e}s de rang $r$.
Il n'est pas difficile de v\'{e}rifier que c'est un champ
alg\'{e}brique (au sens d'Artin), localement de type fini, muni d'un
morphisme de champs
$$
((\infty ,o),(({\cal L}_{1},\lambda_{1}),\ldots ,({\cal L}_{r-1},
\lambda_{r-1}))):{\cal C}^{r}\rightarrow X\times X\times [{\Bbb
A}^{r-1}/{\Bbb G}_{{\rm m}}^{r-1}].
$$
En particulier, pour tout $R\subset [r-1]$, on dispose de la strate
localement ferm\'{e}e ${\cal C}_{R}^{r}\subset {\cal C}^{r}$ image
r\'{e}ciproque de la strate $[{\Bbb A}_{R}^{r-1}/{\Bbb G}_{{\rm
m}}^{r-1}]\subset [{\Bbb A}^{r-1}/{\Bbb G}_{{\rm m}}^{r-1}]$.

D'apr\`{e}s la proposition du paragraphe 3.1, si $R=\{r_{1},
\ldots ,r_{s-1}\}$ o\`{u} $0=r_{0}<r_{1}<\cdots <r_{s-1}<r_{s}=r$, la
donn\'{e}e d'un $U$-point de ${\cal C}_{R}^{r}$ au-dessus d'un
$U$-point $(({\cal L}_{1},\lambda_{1}),\ldots ,({\cal L}_{r-1},
\lambda_{r-1}))$ de $[{\Bbb A}_{R}^{r-1}/{\Bbb G}_{{\rm m}}^{r-1}]$
\'{e}quivaut aux donn\'{e}es suivantes:
\vskip 2mm

\itemitem{-} un diagramme
$$
{\cal E}\,\smash{\mathop{ \lhook\joinrel\mathrel{\hbox to
6mm{\rightarrowfill}} }\limits^{\scriptstyle j}}\,{\cal
E}'\,\smash{\mathop{ {\hbox to 6mm{\leftarrowfill}}\joinrel\kern
-0.9mm\mathrel\rhook} \limits^{\scriptstyle j'}}\,{\cal E}''
$$ 
et deux morphismes $\infty :U\rightarrow X$ et $o:U\rightarrow X$,
o\`{u} ${\cal E}$ est un fibr\'{e} vectoriel de rang $r$ sur $U\times
X$, $j$ est une modification \'{e}l\'{e}mentaire sup\'{e}rieure de
${\cal E}$ le long du graphe de $\infty$ et $j'$ est une modification
\'{e}l\'{e}mentaire inf\'{e}rieure de ${\cal E}'$ le long du graphe de
$o$,
\vskip 2mm

\itemitem{-} une filtration d\'{e}croissante
$$
{}^{\tau}{\cal E}={\cal F}^{0}\supsetneq {\cal F}^{1}\supsetneq\cdots
\supsetneq {\cal F}^{s-1}\supsetneq {\cal F}^{s}=(0)
$$
de ${}^{\tau}{\cal E}$ par des sous-${\cal O}_{U\times X}$-Modules
localement facteurs directs, de corangs $0=r_{0},r_{1},\ldots
,r_{s-1},r_{s}=r$,
\vskip 2mm

\itemitem{-} une filtration croissante
$$
(0)={\cal E}_{0}''\subsetneq {\cal E}_{1}''\subsetneq\cdots\subsetneq
{\cal E}_{s-1}''\subsetneq {\cal E}_{s}''={\cal E}''
$$
de ${\cal E}''$ par des sous-${\cal O}_{U\times X}$-Modules localement
facteurs directs, de rangs $0=r_{0},r_{1},\ldots ,r_{s-1},r_{s}=r$,
\vskip 2mm

\itemitem{-} une famille d'isomorphismes de ${\cal O}_{U\times
X}$-Modules
$$
v_{\sigma}:({\cal F}^{\sigma -1}/{\cal F}^{\sigma})\otimes\left(
\bigotimes_{\sigma'=1}^{\sigma -1}{}^{\tau}{\cal L}_{r_{\sigma'}}
\right) \buildrel\sim\over\longrightarrow ({\cal E}_{\sigma}''/{\cal
E}_{\sigma -1}'')\otimes\left(\bigotimes_{\sigma'=0}^{\sigma -1}{\cal
L}_{r_{\sigma'}}\right),~\sigma =1,\ldots ,s.
$$
\vskip 2mm

Pour un tel $U$-point notons
$$
{\cal E}_{\sigma}'=j'({\cal E}_{\sigma}'')\subset {\cal E}',~{\cal
E}_{\sigma}=j^{-1}({\cal E}_{\sigma}')\subset {\cal E},~\forall\sigma
=0,1,\ldots ,s-1,\hbox{ et }{\cal E}_{s}'={\cal E}',~{\cal E}_{s}=
{\cal E},
$$
$$
{\cal G}_{1}={\cal E}_{1},~{\cal G}_{1}'={\cal E}_{1}'\hbox{ et }{\cal
G}_{\sigma}={\cal E}_{\sigma}/{\cal E}_{\sigma -1},~{\cal
G}_{\sigma}'= {\cal F}^{\sigma -1}\cap {}^{\tau}{\cal
E}_{\sigma},~\forall\sigma =2,\ldots ,s.
$$
On dispose de l'homomorphisme ${\cal G}_{1}\hookrightarrow {\cal
G}_{1}'$ induit par $j$ et du compos\'{e}
$$
{}^{\tau}{\cal G}_{1}\hookrightarrow {}^{\tau}{\cal E}
\twoheadrightarrow {}^{\tau}{\cal E}/{\cal F}^{1}
\buildrel\sim\over\longrightarrow {\cal E}_{1}''
\,\smash{\mathop{\hbox to 6mm{\rightarrowfill}}
\limits^{\scriptstyle j'}}\,{\cal G}_{1}',
$$
et, pour chaque $\sigma =2,\ldots ,r$, on dispose du compos\'{e}
$$\displaylines{
{\cal G}_{\sigma}'\otimes\left(\bigotimes_{\sigma'=0}^{\sigma
-1}{}^{\tau}{\cal L}_{r_{\sigma'}}\right)\hookrightarrow {\cal
F}^{\sigma -1}\otimes\left(\bigotimes_{\sigma'=0}^{\sigma -1}
{}^{\tau}{\cal L}_{r_{\sigma'}}\right)\twoheadrightarrow ({\cal
F}^{\sigma -1}/{\cal F}^{\sigma})\otimes\left(
\bigotimes_{\sigma'=0}^{\sigma -1}{}^{\tau}{\cal L}_{r_{\sigma'}}
\right)
\hfill\cr\hfill
\buildrel\sim\over\longrightarrow ({\cal E}_{\sigma}''/{\cal
E}_{\sigma -1}'')\otimes\left(\bigotimes_{\sigma'=0}^{\sigma -1}{\cal
L}_{r_{\sigma'}}\right)\,\smash{\mathop{\hbox to 6mm{\rightarrowfill}}
\limits^{\scriptstyle \overline{j'}}}\, ({\cal E}_{\sigma}'/{\cal
E}_{\sigma -1}')\otimes\left(\bigotimes_{\sigma'=0}^{\sigma -1}{\cal
L}_{r_{\sigma'}}\right)}
$$
o\`{u} $\overline{j}'$ est induite par $j'$, de l'homomorphisme
$$
\overline{j}:{\cal G}_{\sigma} \otimes\left(
\bigotimes_{\sigma'=0}^{\sigma -1}{\cal L}_{r_{\sigma'}}\right)
\rightarrow ({\cal E}_{\sigma}'/{\cal E}_{\sigma -1}')\otimes
\left(\bigotimes_{\sigma'=0}^{\sigma -1}{\cal L}_{r_{\sigma'}}\right)
$$
induit par $j$, et du compos\'{e}
$$
{\cal G}_{\sigma}'\hookrightarrow {}^{\tau}{\cal E}_{\sigma}
\twoheadrightarrow {}^{\tau}{\cal E}_{\sigma}/{}^{\tau}{\cal
E}_{\sigma -1}={}^{\tau}{\cal G}_{\sigma}.
$$
On a donc des diagrammes
$$
\widetilde{{\cal G}}_{1}=({\cal G}_{1}\hookrightarrow {\cal G}_{1}'
\hookleftarrow {}^{\tau}{\cal G}_{1})
$$
et
$$
\widetilde{{\cal G}}_{\sigma}=\pmatrix{\displaystyle
{\cal G}_{\sigma} \otimes\left(\bigotimes_{\sigma'=0}^{\sigma
-1}{\cal L}_{r_{\sigma'}}\right)\kern -2mm&\cr
\noalign{\medskip}
\llap{$\scriptstyle $}\left\downarrow \vbox to 4mm{}\right.\rlap{}
\kern -2mm&\cr
\noalign{\medskip}
\hfill\displaystyle ({\cal E}_{\sigma}'/{\cal E}_{\sigma -1}')\otimes
\left(\bigotimes_{\sigma'=0}^{\sigma -1}{\cal L}_{r_{\sigma'}}\right)
\kern -2mm&\hookleftarrow\displaystyle {\cal G}_{\sigma}'
\otimes\left( \bigotimes_{\sigma'=0}^{\sigma -1}{}^{\tau}{\cal
L}_{r_{\sigma'}}\right)\hookrightarrow\displaystyle {}^{\tau}{\cal
G}_{\sigma} \otimes\left(\bigotimes_{\sigma'=0}^{\sigma
-1}{}^{\tau}{\cal L}_{r_{\sigma'}}\right)\kern -1mm\cr}
$$
pour $\sigma =2,\ldots ,s$. 

\thm LEMME 
\enonce
Pour tout $\rho =1,\ldots ,r$, notons $u_{\rho}^{\circ}$ la
restriction de $u_{\rho}$ \`{a} l'ouvert de $U\times X$
compl\'{e}mentaire des graphes des morphismes $\infty$ et $o$, et
identifions les restrictions de ${\cal E}$ et ${\cal E}''$ \`{a} cet
ouvert \`{a} l'aide de $j'^{-1}\circ j$. Avec ces notations les deux
conditions ci-dessous sont \'{e}quivalentes:

\decale{\rm (a)} il n'existe aucun point g\'{e}om\'{e}trique de $U$
tel que la fibre en ce point d'un des $u_{\rho}^{\circ}$ soit
nilpotente.

\decale{\rm (a')} pour tout $\sigma =1,\ldots ,s-1$, on a
$({}^{\tau}\!j)({\cal F}^{\sigma})\cap {}^{\tau}{\cal E}_{\sigma}'
=(0)$ dans ${}^{\tau}{\cal E}'$.
\vskip 2mm

De plus, si ces conditions \'{e}quivalentes sont v\'{e}rifi\'{e}es, le
plus grand ouvert de $U\times X$ o\`{u} toutes les fl\`{e}ches des
diagrammes $\widetilde{{\cal G}}_{\sigma}$ ci-dessus sont des
isomorphismes rencontre chaque fibre de la projection canonique
$U\times X\rightarrow U$ suivant un ouvert dense.
\endthm

La condition (a) a \'{e}t\'{e} introduite par Drinfeld en rang $r=2$.
\vskip 3mm

Suivant Lafforgue consid\'{e}rons les conditions suppl\'{e}mentaires
suivantes:

\decale{\rm (b)} ${\cal E}'/{\cal E}_{\sigma}'$ est un ${\cal
O}_{U\times S}$-Module localement libre pour $\sigma =0,1,\ldots
,s-1$,

\decale{\rm (c)} pour tout $\sigma =1,\ldots ,s$, l'homomorphisme
compos\'{e} ${\cal E}_{\sigma}''\,\smash{\mathop{\lhook\joinrel
\mathrel{\hbox to 6mm{\rightarrowfill}}}\limits^{\scriptstyle
j'}}\,{\cal E}'\twoheadrightarrow {\cal E}'/j({\cal E})$ est
surjectif,

\decale{\rm (d)} pour tout $\sigma =1,\ldots ,s-1$, on a ${\cal 
F}^{\sigma -1}+{}^{\tau}{\cal E}_{\sigma}={}^{\tau}{\cal E}$,
\vskip 2mm

La condition (c) assure en particulier que $\overline{j}: {\cal
G}_{\sigma} \otimes\left( \bigotimes_{\sigma'=0}^{\sigma -1}{\cal
L}_{r_{\sigma'}}\right) \rightarrow ({\cal E}_{\sigma}'/{\cal
E}_{\sigma -1}')\otimes \left(\bigotimes_{\sigma'=0}^{\sigma -1}{\cal
L}_{r_{\sigma'}}\right)$ est un isomorphisme. On peut donc r\'{e}crire
les diagrammes $\widetilde{{\cal G}}_{\sigma}$, $\sigma =2,\ldots ,s$,
sous la forme
$$
\widetilde{{\cal G}}_{\sigma}=\left({\cal G}_{\sigma} \otimes\left(
\bigotimes_{\sigma'=0}^{\sigma -1}{\cal
L}_{r_{\sigma'}}\right)\hookleftarrow {\cal G}_{\sigma}'\otimes\left(
\bigotimes_{\sigma'=0}^{\sigma -1}{}^{\tau}{\cal
L}_{r_{\sigma'}}\right)\hookrightarrow {}^{\tau}{\cal G}_{\sigma}
\otimes\left( \bigotimes_{\sigma'=0}^{\sigma -1}{}^{\tau}{\cal
L}_{r_{\sigma'}}\right)\right).
$$

\thm LEMME
\enonce
{\rm (i)} Sous les conditions {\rm (a), (b), (c)} et {\rm (d)} le
diagramme $\widetilde{{\cal G}}_{1}$ est un chtouca \`{a} droite de
rang $r_{1}$ et de p\^{o}le $\infty$, les diagrammes $\widetilde{{\cal
G}}_{2},\ldots ,\widetilde{{\cal G}}_{s}$ sont des chtoucas \`{a}
gauche, le z\'{e}ro de $\widetilde{{\cal G}}_{\sigma}$ est le p\^{o}le
de $\widetilde{{\cal G}}_{\sigma +1}$ pour $\sigma =1,\ldots ,s-1$ et
le z\'{e}ro de $\widetilde{{\cal G}}_{s}$ est $o$.

\decale{\rm (ii)} Il existe un unique sous-champ ouvert
$\mathop{\overline{\rm Cht}}\nolimits^{\,r}\subset {\cal C}^{r}$ tel
que pour chaque $R=\{r_{1}, \ldots ,r_{s-1}\}$ comme ci-dessus la
trace de cet ouvert sur ${\cal C}_{R}^{r}$ soit d\'{e}finie par les
conditions {\rm (a), (b), (c)} et {\rm (d)}.
\endthm

Le champ $\mathop{\overline{\rm Cht}}\nolimits^{\,r}$ est par
d\'{e}finition le champ des {\it chtoucas it\'{e}r\'{e}s} de rang $r$.
Il contient comme ouvert dense le champ des chtoucas $\mathop{\rm
Cht}\nolimits^{\,r}$. Tout comme un chtouca ordinaire, un chtouca
it\'{e}r\'{e} admet un degr\'{e}, \`{a} savoir le degr\'{e} du
fibr\'{e} vectoriel sous-jacent ${\cal E}$, et on a un d\'{e}coupage de
$\mathop{\overline{\rm Cht}}\nolimits^{\,r}$ en composantes
$$
\mathop{\overline{\rm Cht}}\nolimits^{\,r}=\coprod_{d\in {\Bbb Z}}
\mathop{\overline{\rm Cht}}\nolimits^{\,r,d}
$$
index\'{e}es par ce degr\'{e}. Bien s\^{u}r, pour chaque entier $d$,
$\mathop{\rm Cht}\nolimits^{\,r,d}$ est un ouvert dense de
$\mathop{\overline{\rm Cht}}\nolimits^{\,r,d}$.

Pour chaque $R=\{r_{1}, \ldots ,r_{s-1}\}\subset [r-1]$ on notera
$\mathop{\overline{\rm Cht}}\nolimits_{R}^{\,r}$ l'intersection de
$\mathop{\overline{\rm Cht}}\nolimits^{\,r}$ avec la strate ${\cal
C}_{R}^{r}$. C'est un sous-champ localement ferm\'{e} de
$\mathop{\overline{\rm Cht}}\nolimits^{\,r}$, qui n'est autre que
l'ouvert $\mathop{\rm Cht}^{r}$ pour $R=\emptyset$, et
$\mathop{\overline{\rm Cht}}\nolimits^{\,r}$ est la r\'{e}union
disjointes des $\mathop{\overline{\rm Cht}}\nolimits_{R}^{\,r}$.

Il r\'{e}sulte du dernier lemme que l'on a un morphisme fini,
surjectif et radiciel de champs
$$
\mathop{\overline{\rm Cht}}\nolimits_{R}^{\,r}\rightarrow\mathop{\rm
Cht}\nolimits^{r_{1}}\times_{X} {}^{r_{2}-r_{1}}\!\mathop{\rm Cht}
\cdots\times_{X} {}^{r_{s}-r_{s-1}}\!\mathop{\rm Cht}.
$$
De plus on a un morphisme de Frobenius partiel (au dessus de
l'endomorphisme $\mathop{\rm Id}\nolimits_{X}\times \mathop{\rm
Frob}\nolimits_{X}$ de $X\times X$)
$$\displaylines{
\qquad\mathop{\rm Cht}\nolimits^{r_{1}}
\times_{X}{}^{r_{2}-r_{1}}\!\mathop{\rm Cht}\cdots
\times_{X}{}^{r_{s}-r_{s-1}}\!\mathop{\rm Cht}
\hfill\cr\hfill
\rightarrow\mathop{\rm Cht}\nolimits^{R}:=\mathop{\rm
Cht}\nolimits^{r_{1}}\times_{X}
\mathop{\rm Cht}\nolimits^{r_{2}-r_{1}} \times_{X,\mathop{\rm
Frob}\nolimits_{X}}\cdots\times_{X,\mathop{\rm Frob}
\nolimits_{X}}\mathop{\rm Cht} \nolimits^{r_{s}-r_{s-1}}\qquad}
$$
qui envoie $(\widetilde{{\cal G}}_{1},\widetilde{{\cal G}}_{2},\ldots
,\widetilde{{\cal G}}_{s})$ sur
$$
(\widetilde{{\cal G}}_{1},\mathop{\rm Frob}\nolimits_{o}
(\widetilde{{\cal G}}_{2}),\ldots ,\mathop{\rm Frob}\nolimits_{o}
(\widetilde{{\cal G}}_{s}))
$$
qui est lui aussi fini, surjectif et radiciel. Par composition on a
donc un morphisme fini, surjectif et radiciel de champs
$$
\mathop{\overline{\rm Cht}}\nolimits_{R}^{\,r}\rightarrow\mathop{\rm
Cht}\nolimits^{R}.
$$
Le but de ce dernier morphisme est un champ de Deligne-Mumford
s\'{e}par\'{e} et lisse purement de dimension relative
$$
(2r_{1}-2)+(2(r_{2}-r_{1})-2)+\cdots +(2(r_{s}-r_{s-1})-2)=2r-2s
$$
au dessus de $X\times X^{s-1}\times X$ par le morphisme $(\infty
=\infty_{1},\infty_{2},\ldots ,\infty_{s},o_{s}=o)$. Il se 
d\'{e}compose en
$$
\mathop{\rm Cht}\nolimits^{R}=\coprod_{\scriptstyle d_{\bullet}\in
{\Bbb Z}^{s}} \mathop{\rm Cht}\nolimits^{R,d_{\bullet}}
$$
o\`{u}
$$
\mathop{\rm Cht}\nolimits^{R,d_{\bullet}}:=\mathop{\rm Cht}
\nolimits^{r_{1},d_{1}}\times_{X} \mathop{\rm Cht}
\nolimits^{r_{2}-r_{1},d_{2}} \times_{X,\mathop{\rm Frob}
\nolimits_{X}}\cdots\times_{X,\mathop{\rm Frob}
\nolimits_{X}}\mathop{\rm Cht} \nolimits^{r_{s}-r_{s-1},d_{s}}.
$$
On notera 
$$
\mathop{\overline{\rm Cht}}\nolimits_{R}^{\,r}=\coprod_{\scriptstyle
d_{\bullet}\in {\Bbb Z}^{s}}\mathop{\overline{\rm Cht}}
\nolimits_{R}^{\,r,d_{\bullet}}
$$
la d\'{e}composition correspondante de la source. On v\'{e}rifie que
$$
\mathop{\overline{\rm Cht}}\nolimits_{R}^{\,r,d}:=
\mathop{\overline{\rm Cht}}\nolimits_{R}^{\,r}\cap
\mathop{\overline{\rm Cht}}\nolimits^{\,r,d}=\coprod_{{\scriptstyle
d_{\bullet}\in {\Bbb Z}^{s}\atop\scriptstyle d_{1}+\cdots
d_{s}=d-s+1}}\mathop{\overline{\rm Cht}}
\nolimits_{R}^{\,r,d_{\bullet}}
$$
(le d\'{e}calage $-s+1$ provient des Frobenius partiels).
\vskip 5mm

{\bf 3.3. Chtoucas it\'{e}r\'{e}s et troncature}
\vskip 5mm

Soient $d$ un entier et $p:[0,r]\rightarrow {\Bbb R}$ un param\`{e}tre
de troncature. Pour chaque $\rho =0,1,\ldots ,r$, notons
$\widetilde{p}(\rho)$ l'unique entier appartenant \`{a} l'intervalle
de longueur $1$
$$
]p(\rho)+{\rho\over r}d-1,p(\rho)+{\rho\over r}d].
$$
On d\'{e}finit alors des entiers $d_{1},\ldots ,d_{s}$ par
$$
d_{1}=\widetilde{p}(r_{1})\hbox{ et }d_{\sigma}=
\widetilde{p}(r_{\sigma})-\widetilde{p}(r_{\sigma -1})-1,~
\forall\sigma =2,\ldots ,s,
$$
et des param\`{e}tres de troncature $p_{1}:[0,r_{1}] \rightarrow {\Bbb
R},\ldots ,p_{s}:[0,r_{s}-r_{s-1}]\rightarrow {\Bbb R}$ en imposant
que
$$
p_{1}(\rho_{1})=\widetilde{p}(\rho_{1})-{\rho_{1}\over
r_{1}}d_{1},~\forall \rho_{1}=1,\ldots ,r_{1}-1,
$$
et
$$
p_{\sigma}(\rho_{\sigma})=\widetilde{p}(r_{\sigma -1}+\rho_{\sigma})
-\widetilde{p}(r_{\sigma -1})-1-{\rho_{\sigma}\over
r_{\sigma}-r_{\sigma -1}}d_{\sigma},~ \forall \rho_{\sigma}=1,\ldots 
,r_{\sigma}-r_{\sigma -1}.
$$
\`{A} d\'{e}calage et normalisation pr\`{e}s, les $p_{\sigma}$ sont
essentiellement les restrictions de $p$ aux intervalles
$[r_{\sigma-1},r_{\sigma}]$ comme dans la figure ci-dessous:
$$\kern -69mm
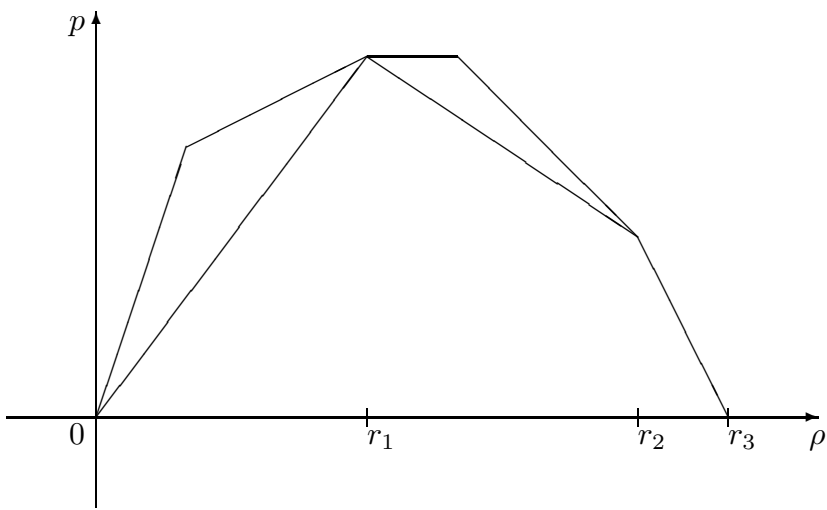
\diagram{\noalign{\vskip 45mm}
\arrow(-10,0)\dir(1,0)\length{90}
\arrow(0,-10)\dir(0,1)\length{55}
\rule(0,0)\dir(1,3)\long{10}
\rule(10,30)\dir(2,1)\long{20}
\rule(30,40)\dir(1,0)\long{10}
\rule(40,40)\dir(1,-1)\long{20}
\rule(60,20)\dir(1,-2)\long{10}
\rule(0,0)\dir(3,4)\long{30}
\rule(30,40)\dir(3,-2)\long{30}
\put(-3,-3){0}
\put(-3,43){p}
\put(30,-3){r_{1}}
\rule(30,-1)\dir(0,1)\long{2}
\put(60,-3){r_{2}}
\rule(60,-1)\dir(0,1)\long{2}
\put(70,-3){r_{3}}
\rule(70,-1)\dir(0,1)\long{2}
\put(79,-3){\rho}
\noalign{\vskip 5mm}}
$$

On note
$$\displaylines{
\qquad\mathop{\rm Cht}\nolimits^{R,d;\,\leq p}=\mathop{\rm Cht}
\nolimits^{r_{1},d_{1};\,\leq p_{1}}\times_{X}\mathop{\rm Cht}
\nolimits^{r_{2}-r_{1},d_{2};\,\leq p_{2}}\times_{X,\mathop{\rm Frob}
\nolimits_{X}}\cdots
\hfill\cr\hfill
\cdots\times_{X,\mathop{\rm Frob}\nolimits_{X}}\mathop{\rm Cht}
\nolimits^{r_{s}-r_{s-1},d_{s};\,\leq p_{2}}\subset\mathop{\rm Cht}
\nolimits^{R}.\qquad}
$$

\thm D\'{E}FINITION 
\enonce
Un param\`{e}tre de troncature $p:[0,r]\rightarrow {\Bbb R}$ est dit 
$\mu$-convexe pour un nombre r\'{e}el $\mu\geq 0$ si l'on a
$$
(p(\rho )-p(\rho -1))-(p(\rho +1)-p(\rho ))\geq\mu,~\forall 
1\leq\rho\leq r-1.
$$

Une propri\'{e}t\'{e} est dite vraie {\og}{pour tout param\`{e}tre de
troncature assez convexe}{\fg} s'il existe un nombre r\'{e}el $\mu\geq
0$ tel que la propri\'{e}t\'{e} soit vraie pour tout param\`{e}tre de
troncature qui est $\mu$-convexe.
\endthm

On remarquera que les param\`{e}tres de troncature $p_{\sigma}$
d\'{e}finis ci-dessus sont automatiquement $(\mu -2)$-convexes d\`{e}s 
que $p$ est $\mu$-convexe.

\thm TH\'{E}OR\`{E}ME 
\enonce
Pour tout entier $d$ et tout param\`{e}tre de troncature
$p:[0,r]\rightarrow {\Bbb R}$ qui est $2$-convexe, il existe un ouvert
$$
\mathop{\overline{\rm Cht}}\nolimits^{\,r,d;\,\leq p}\subset
\mathop{\overline{\rm Cht}}\nolimits^{\,r,d}
$$
de la composante $d$-i\`{e}me du champ des chtoucas it\'{e}r\'{e}s
ayant les propri\'{e}t\'{e}s suivantes:
\vskip 2mm

\itemitem{-} pour tout $R\subset [r-1]$, $\mathop{\overline{\rm
Cht}}\nolimits_{R}^{\,r,d;\,\leq p}:=\mathop{\overline{\rm Cht}}
\nolimits_{R}^{\,r,d}\cap \mathop{\overline{\rm Cht}}
\nolimits^{\,r,d;\,\leq p}$ est l'image r\'{e}ciproque par le
morphisme fini, surjectif et radiciel $\mathop{\overline{\rm Cht}}
\nolimits_{R}^{\,r}\rightarrow\mathop{\rm Cht}\nolimits^{R}$ de
l'ouvert $\mathop{\rm Cht}\nolimits^{R,d;\,\leq p}$;
\vskip 2mm

\itemitem{-} le morphisme de champs $ (\infty ,o):
\mathop{\overline{\rm Cht}}\nolimits^{\,r,d;\,\leq p}=
\mathop{\overline{\rm Cht}}\nolimits^{\,r;\,\leq p}\cap
\mathop{\overline{\rm Cht}}\nolimits^{\,r,d}\rightarrow X\times X$ est
propre {\rm (}et donc s\'{e}par\'{e} et de type fini{\rm )};
\vskip 2mm

\itemitem{-} si l'on suppose de plus que $p$ est assez convexe
relativement au genre de la courbe $X$, le morphisme de champs
$(\infty ,o)$ ci-dessus est lisse purement de dimension relative
$2r-2$ et le ferm\'{e} compl\'{e}mentaire dans $\mathop{\overline{\rm
Cht}}\nolimits^{\,r,d;\,\leq p}$ de la strate ouverte $\mathop{\rm
Cht}\nolimits^{r,d;\,\leq p}=\mathop{\overline{\rm
Cht}}\nolimits_{\emptyset}^{\,r,d;\,\leq p}$ est un diviseur \`{a}
croisements normaux relatif sur $X\times X$, qui est r\'{e}union de
$r-1$ diviseurs lisses et dont la stratification canonique est celle
par les $\mathop{\overline{\rm Cht}}\nolimits_{R}^{\,r,d;\,\leq p}$
pour $R$ parcourant les parties non vides de $[r-1]$.

\endthm

On posera 
$$
\mathop{\overline{\rm Cht}}\nolimits^{\,r;\,\leq p}=\coprod_{d\in
{\Bbb Z}} \mathop{\overline{\rm Cht}}\nolimits^{\,r,d;\,\leq p}\subset
\mathop{\overline{\rm Cht}}\nolimits^{\,r}
$$
et
$$
\mathop{\overline{\rm Cht}}\nolimits_{R}^{\,r;\,\leq p}=\coprod_{d\in
{\Bbb Z}} \mathop{\overline{\rm Cht}}_{R}\nolimits^{\,r,d;\,\leq p}\subset
\mathop{\overline{\rm Cht}}\nolimits^{\,r;\,\leq p}.
$$
\vskip 5mm

{\bf 3.4. Structures de niveau sur les chtoucas it\'{e}r\'{e}s}
\vskip 5mm

{\it Dans cet expos\'{e} nous nous limiterons au cas o\`{u} $N$ est
r\'{e}duit et support\'{e} par un point ferm\'{e} de $X$ rationnel sur
le corps fini de base ${\Bbb F}_{q}$}. On suppose donc $N\cong
\mathop{\rm Spec}({\Bbb F}_{q})$ et on a simplement $\mathop{\rm
GL}\nolimits_{r}^{N}=\mathop{\rm GL}\nolimits_{r}$.

Soit
$$
{\cal E}\,\smash{\mathop{\lhook\joinrel\mathrel{\hbox to
5mm{\rightarrowfill}} }\limits^{\scriptstyle j}}\,{\cal E}'\,
\smash{\mathop{ {\hbox to 5mm{\leftarrowfill}}\joinrel\kern
-0.9mm\mathrel\rhook} \limits^{\scriptstyle j'}}\,{\cal E}''\buildrel
u\over\Longleftarrow {}^{\tau}{\cal E},
$$ 
un chtouca it\'{e}r\'{e} sur un sch\'{e}ma $U$; $u=(u_{1},\ldots
,u_{r};({\cal L}_{1}^{\otimes (q-1)},\lambda_{1}^{q-1}),\ldots ,({\cal
L}_{r-1}^{\otimes (q-1)},\lambda_{r-1}^{q-1}))$ est donc un
homomorphisme complet de ${}^{\tau}{\cal E}$ dans ${\cal E}''$. Si le
p\^{o}le et le z\'{e}ro de ce chtouca it\'{e}r\'{e} ne rencontrent pas
$N$, les restrictions
$$
j_{N}:{\cal E}_{N}\rightarrow {\cal E}_{N}'\hbox{ et }j_{N}':{\cal
E}_{N}''\rightarrow {\cal E}_{N}'
$$
des fl\`{e}ches $j$ et $j'$ sont des isomorphismes de ${\cal
O}_{U}$-Modules localement libres de rang $r$ et la restriction
$$
u_{N}=(u_{1,N},\ldots ,u_{r,N};({\cal L}_{1}^{\otimes
(q-1)},\lambda_{1}^{q-1}),\ldots ,({\cal L}_{r-1}^{\otimes
(q-1)},\lambda_{r-1}^{q-1})):{}^{\tau}{\cal E}_{N}\Longrightarrow
{\cal E}_{N}
$$
est un homomorphisme complet. Le morphisme de restriction \`{a} $N$
$$
\mathop{\rm Cht}\nolimits^{r}\times_{X^{2}}^{}(X-N)^{2}\rightarrow
\mathop{\rm Tr}\nolimits_{N}^{r}
$$
se prolonge donc en un morphisme
$$
\mathop{\overline{\rm Cht}}\nolimits^{\,r}\times_{X^{2}}^{}(X-N)^{2}
\rightarrow \mathop{\overline{\rm Tr}}\nolimits_{N}^{\,r}
$$
o\`{u} $\mathop{\overline{\rm Tr}}\nolimits_{N}^{\,r}$ est le champ qui
associe au sch\'{e}ma $U$ la cat\'{e}gorie des uplets
$$
({\cal F}, u_{1},\ldots ,u_{r};({\cal L}_{1},\lambda_{1}),\ldots
,({\cal L}_{r-1},\lambda_{r-1}))
$$
form\'{e}s d'un fibr\'{e} ${\cal F}$ de rang $r$ sur $N$, de
fibr\'{e}s en droites ${\cal L}_{1},\ldots ,{\cal L}_{r-1}$ sur $U$
munis de sections globales $\lambda_{1},\ldots ,\lambda_{r-1}$ et d'un
homomorphisme complet
$$
u=(u_{1},\ldots ,u_{r};({\cal L}_{1}^{\otimes (q-1)},
\lambda_{1}^{q-1}),\ldots ,({\cal L}_{r-1}^{\otimes (q-1)},
\lambda_{r-1}^{q-1})): {}^{\tau}{\cal F}\Longrightarrow {\cal F}.
$$
L'ouvert $\mathop{\rm Tr}\nolimits_{N}^{r}$ de $\mathop{\overline{\rm
Tr}}\nolimits_{N}^{\,r}$ est l'image r\'{e}ciproque par le morphisme
$$
(({\cal L}_{1},\lambda_{1}),\ldots ,({\cal L}_{r-1},\lambda_{r-1})):
\mathop{\overline{\rm Tr}}\nolimits_{N}^{\,r}\rightarrow [{\Bbb
A}^{r-1}/{\Bbb G}_{{\rm m}}^{r-1}]
$$
de la strate ouverte $\mathop{\rm Spec}({\Bbb F}_{q})=[{\Bbb G}_{{\rm
m}}^{r-1}/{\Bbb G}_{{\rm m}}^{r-1}]\subset [{\Bbb A}^{r-1}/{\Bbb
G}_{{\rm m}}^{r-1}]$.

\thm TH\'{E}OR\`{E}ME 
\enonce
Consid\'{e}rons le morphisme
$$
\mathop{\overline{\rm Cht}}\nolimits^{\,r}\times_{X^{2}}^{}(X-N)^{2}
\rightarrow\mathop{\overline{\rm Tr}}\nolimits_{N}^{\,r}\times
(X-N)^{2}
$$
de composante le morphisme pr\'{e}c\'{e}dant et la projection
canonique. 

Pour tout entier $d$ et tout param\`{e}tre de troncature
$p:[0,r]\rightarrow {\Bbb R}$ assez convexe par rapport \`{a} $X$ et
$N$, la restriction
$$
\mathop{\overline{\rm Cht}}\nolimits^{\,r,d;\,\leq p}\times_{X^{2}}
(X-N)^{2} \rightarrow\mathop{\overline{\rm Tr}}\nolimits_{N}^{r}\times
(X-N)^{2}
$$
de ce morphisme \`{a} l'ouvert $\mathop{\overline{\rm
Cht}}\nolimits^{\,r,d;\,\leq p}\subset \mathop{\overline{\rm
Cht}}\nolimits^{\,r}$ est lisse purement de dimension relative $2r-2$.
\endthm

Le champ $\mathop{\overline{\rm Tr}}\nolimits_{N}^{r}$ est le produit 
fibr\'{e}
$$
\mathop{\overline{\rm Tr}}\nolimits_{N}^{r}=[\mathop{\widetilde{\rm
gl}}\nolimits_{r}/{}^{\tau}\!\!\mathop{\rm Ad}(\mathop{\rm
GL}\nolimits_{r})]\times_{[{\Bbb A}^{r-1}/{\Bbb G}_{{\rm
m}}^{r-1}],\langle q-1\rangle}[{\Bbb A}^{r-1}/{\Bbb G}_{{\rm m}}^{r-1}]
$$
o\`{u} $\mathop{\widetilde{\rm gl}}\nolimits_{r}$ est le sch\'{e}ma
des endomorphismes complets du ${\Bbb F}_{q}$-espace vectoriel ${\Bbb
F}_{q}^{r}$, o\`{u} le morphisme $[\mathop{\widetilde{\rm
gl}}\nolimits_{r}/{}^{\tau}\!\!\mathop{\rm Ad}(\mathop{\rm
GL}\nolimits_{r})]\rightarrow [{\Bbb A}^{r-1}/{\Bbb G}_{{\rm
m}}^{r-1}]$ est induit par le morphisme canonique
$\mathop{\widetilde{\rm gl}}\nolimits_{r}\rightarrow [{\Bbb
A}^{r-1}/{\Bbb G}_{{\rm m}}^{r-1}]$ et o\`{u} $\langle q-1\rangle$ est
l'endomorphisme de $[{\Bbb A}^{r-1}/{\Bbb G}_{{\rm m}}^{r-1}]$ induit
par l'\'{e}l\'{e}vation \`{a} la puissance $q-1$ des coordonn\'{e}es.

Lafforgue d\'{e}finit alors le champ $\mathop{\overline{\rm Cht}}
\nolimits_{N}^{\,r,d;\,\leq p}$ de {\it chtoucas it\'{e}r\'{e}s avec
structure de niveau} $N$ par le carr\'{e} cart\'{e}sien
$$\diagram{
\mathop{\overline{\rm Cht}}\nolimits_{N}^{\,r,d;\,\leq p}&\kern
-12mm\smash{\mathop{\hbox to 19mm {\rightarrowfill}}
\limits^{\scriptstyle }}\kern -4mm&\mathop{\overline{\rm
Tr}}\nolimits_{N}^{\,r,\tau}\cr
\llap{$\scriptstyle $}\left\downarrow \vbox to 4mm{}\right.\rlap{}
&\kern -6mm\square&\llap{}\left\downarrow \vbox to 4mm{}
\right.\rlap{$\scriptstyle $}\cr
\mathop{\overline{\rm Cht}}\nolimits^{\,r,d;\,\leq p}\times_{X^{2}}
(X-N)^{2}&\kern -1mm\smash{\mathop{\hbox to 8mm{\rightarrowfill}}
\limits_{\scriptstyle }}\kern -1mm&\mathop{\overline{\rm
Tr}}\nolimits_{N}^{\,r}\cr}
$$
o\`{u} la fl\`{e}che verticale de droite
$$\displaylines{
\qquad\mathop{\overline{\rm
Tr}}\nolimits_{N}^{\,r,\tau}:=[\mathop{\widetilde{\rm gl}}
\nolimits_{r}^{\tau}/\mathop{\rm GL}\nolimits_{r}]\times_{[{\Bbb
A}^{r-1}/{\Bbb G}_{{\rm m}}^{r-1}],\langle q-1\rangle}[{\Bbb
A}^{r-1}/{\Bbb G}_{{\rm m}}^{r-1}]
\hfill\cr\hfill
\rightarrow [\mathop{\widetilde{\rm gl}}
\nolimits_{r}/{}^{\tau}\!\!\mathop{\rm Ad} (\mathop{\rm GL}
\nolimits_{r})]\times_{[{\Bbb A}^{r-1}/{\Bbb G}_{{\rm
m}}^{r-1}],\langle q-1\rangle}[{\Bbb A}^{r-1}/{\Bbb G}_{{\rm
m}}^{r-1}] =\mathop{\overline{\rm Tr}}\nolimits_{N}^{\,r}\qquad}
$$
est induite par un morphisme projectif \'{e}quivariant $\overline{L}:
\mathop{\widetilde{\rm gl}}\nolimits_{r}^{\tau} \rightarrow
\mathop{\widetilde{\rm gl}}\nolimits_{r}$ qui prolonge l'isog\'{e}nie
de Lang $L:\mathop{\rm GL}\nolimits_{r}\rightarrow \mathop{\rm
GL}\nolimits_{r}$.

Par construction, le sch\'{e}ma $\mathop{\widetilde{\rm
gl}}\nolimits_{r}^{\tau}$ est muni d'un morphisme lisse vers un
{\og}{champ torique}{\fg} (voir la sous-section suivante) ${\cal
T}^{r,\tau}$, qui rel\`{e}ve le morphisme $\mathop{\widetilde{\rm
gl}}\nolimits_{r} \rightarrow [{\Bbb A}^{r-1}/{\Bbb G}_{{\rm
m}}^{r-1}]$. On a donc un morphisme lisse de champs
$$
\mathop{\overline{\rm Tr}}\nolimits_{N}^{\,r,\tau}\rightarrow {\cal
T}^{r,\tau}\times_{[{\Bbb A}^{r-1}/{\Bbb G}_{{\rm m}}^{r-1}],\langle
q-1\rangle}[{\Bbb A}^{r-1}/{\Bbb G}_{{\rm m}}^{r-1}]
$$
et il r\'{e}sulte du th\'{e}or\`{e}me ci-dessus que
$\mathop{\overline{\rm Cht}}\nolimits_{N}^{\,r,d;\,\leq p}$ est lisse
sur le champ
$$
\left({\cal T}^{r,\tau}\times_{[{\Bbb A}^{r-1}/{\Bbb G}_{{\rm
m}}^{r-1}],\langle q-1\rangle}[{\Bbb A}^{r-1}/{\Bbb G}_{{\rm
m}}^{r-1}]\right)\times (X-N)^{2}.
$$
En particulier, comme on sait r\'{e}soudre les singularit\'{e}s d'un
champ torique en raffinant l'\'{e}ventail qui le d\'{e}finit, on sait
aussi r\'{e}soudre les singularit\'{e}s du champ
$\mathop{\overline{\rm Cht}}\nolimits_{N}^{\,r,d;\,\leq p}$.
\vskip 3mm

La section suivante est consacr\'{e} \`{a} la construction de
$\mathop{\widetilde{\rm gl}} \nolimits_{r}^{\tau}$.

\vskip 7mm

\centerline{\bf 4. COMPACTIFICATION DU CLASSIFIANT DE $\mathop{\rm
PGL}\nolimits_{r}$}
\vskip 10mm

{\bf 4.1. Champ des pavages} 
\vskip 5mm

Pour compactifier l'isog\'{e}nie de Lang pour $\mathop{\rm GL}
\nolimits_{r}$ Lafforgue est amen\'{e} \`{a} construire une famille de
champs toriques $({\cal T}^{r,n})_{n\geq 0}$ qu'il appelle les {\it
champs des pavages}.
\vskip 3mm

Une {\it vari\'{e}t\'{e} torique} est un triplet form\'{e} d'un tore
$T$, d'une immersion ouverte $T\hookrightarrow\overline{T}$ d'image
dense de $T$ dans un sch\'{e}ma de type fini normal $\overline{T}$, et
d'une action de $T$ sur $\overline{T}$ qui prolonge l'action de $T$
sur lui-m\^{e}me par translation. Une telle vari\'{e}t\'{e} torique
est associ\'{e}e \`{a} un \'{e}ventail. Rappelons tout d'abord
rapidement cette construction.

Soient $Y$ un ${\Bbb Z}$-module libre de rang fini et $Y_{{\Bbb R}}=
{\Bbb R}\otimes_{{\Bbb Z}}^{}Y$. Un {\it c\^{o}ne convexe
poly\'{e}dral rationnel} est une partie localement ferm\'{e}e $\sigma$
de $Y$ de la forme
$$
\sigma ={\Bbb R}_{>0}\,y_{1}+\cdots +{\Bbb R}_{>0}\,y_{n}\subset
Y_{{\Bbb R}}
$$
pour des vecteurs $y_{1},\ldots ,y_{n}\in Y\subset Y_{{\Bbb R}}$.
L'adh\'{e}rence $\overline{\sigma}$ de $\sigma$ est le c\^{o}ne
convexe ferm\'{e}
$$
\sigma ={\Bbb R}_{\geq 0}\,y_{1}+\cdots +{\Bbb R}_{\geq 0}\,y_{n}
\subset Y_{{\Bbb R}}.
$$
On d\'{e}finit les faces de $\sigma$ comme les c\^{o}nes convexes
poly\'{e}draux rationnels $\tau$ de la forme $\overline{\tau}=
H\cap\overline{\sigma}$ pour un hyperplan $H$ de $Y_{{\Bbb R}}$ ne
rencontrant pas $\sigma$. Un {\it \'{e}ventail} est une famille finie
$E$ de c\^{o}nes convexes poly\'{e}draux rationnels deux \`{a} deux
disjoints, appel\'{e}s les {\it cellules} de $E$ et ayant les
propri\'{e}t\'{e}s suivantes:
\vskip 2mm

\itemitem{-} toute face d'une cellule de $E$ est encore une cellule de
$E$,
\vskip 2mm

\itemitem{-} tout couple $(\sigma_{1},\sigma_{2})$ de cellules
de $E$ admet une face commune $\tau\in E$ telle que
$\overline{\sigma}_{1}\cap\overline{\sigma}_{2}=\overline{\tau}$.
\vskip 2mm

\itemitem{-} le c\^{o}ne convexe poly\'{e}dral rationnel $\{0\}$ est
une cellule de $E$.
\vskip 2mm

Si $\sigma ={\Bbb R}_{>0}\,y_{1}+\cdots +{\Bbb R}_{>0}
\,y_{n}\subset Y_{{\Bbb R}}$ est un c\^{o}ne convexe poly\'{e}dral
rationnel ne contenant aucune droite vectorielle de $Y_{{\Bbb R}}$, la
collection form\'{e}e de $\sigma$ et de toutes ses faces est un
\'{e}ventail. On associe \`{a} cet \'{e}ventail particulier la
vari\'{e}t\'{e} torique affine
$$
T=Y\otimes_{{\Bbb Z}}^{}{\Bbb G}_{{\rm m}}=\mathop{\rm Spec}({\Bbb
F}_{q}[Y^{\vee}])\hookrightarrow \overline{T}(\sigma )=\mathop{\rm
Spec}({\Bbb F}_{q}[Y^{\vee}\cap\sigma^{\vee}])
$$
o\`{u} $Y^{\vee}:=\mathop{\rm Hom}\nolimits_{{\Bbb Z}}^{}(Y,{\Bbb Z})$
et
$$
\sigma^{\vee}=\{y^{\vee}:Y_{{\Bbb R}}\rightarrow {\Bbb R}\mid
y^{\vee}(y_{1}), \ldots ,y^{\vee}(y_{n})\geq 0\}\subset Y_{{\Bbb
R}}^{\vee}
$$
est le c\^{o}ne ferm\'{e} dual de $\sigma$. On v\'{e}rifie que
$Y^{\vee}\cap\sigma^{\vee}$ est un mono\"{\i}de \`{a} engendrement
fini, satur\'{e} et contenant un syst\`{e}me de g\'{e}n\'{e}rateurs du
groupe $Y^{\vee}$. L'exemple type est bien entendu le plongement
torique affine standard ${\Bbb G}_{{\rm m}}^{r-1}\hookrightarrow {\Bbb
A}^{r-1}$ qui est obtenu par la construction ci-dessus pour $Y={\Bbb
Z}^{r-1}$ et $\sigma =({\Bbb R}_{>0})^{r-1}\subset {\Bbb
R}^{r-1}=Y_{{\Bbb R}}$.

Si $\tau$ est une face de $\sigma$, on a $\sigma^{\vee}\subset
\tau^{\vee}$ et l'immersion ouverte $T\hookrightarrow\overline{T}
(\sigma )$ se factorise par l'immersion ouverte
$T\hookrightarrow\overline{T} (\tau )$ en
$$
T\hookrightarrow\overline{T}(\tau )\hookrightarrow\overline{T}(\sigma
).
$$
En particulier, si $\sigma_{1}$ et $\sigma_{2}$ sont deux cellules de
$E$ et si $\overline{\tau}=\overline{\sigma}_{1}\cap
\overline{\sigma}_{2}$, on a une donn\'{e}e de recollement
$$
\overline{T}(\sigma_{1})\hookleftarrow\overline{T}(\tau )
\hookrightarrow\overline{T}(\sigma_{2})
$$
pour les vari\'{e}t\'{e}s toriques affines $T\hookrightarrow
\overline{T}(\sigma_{1})$ et $T\hookrightarrow\overline{T}
(\sigma_{2})$.

Le vari\'{e}t\'{e} torique $T\hookrightarrow\overline{T}$ associ\'{e}e
\`{a} l'\'{e}ventail $E$ est obtenue en recollant les vari\'{e}t\'{e}s
toriques affines $T\hookrightarrow\overline{T}(\sigma )$ pour
$\sigma\in E$ \`{a} l'aide des donn\'{e}es de recollements ci-dessus.

Les orbites pour l'action de $T$ sur $\overline{T}$ sont en nombre
fini et forment une stratification en parties localement ferm\'{e}es
de $\overline{T}$. On peut indexer ces orbites par les
cellules de $E$: l'orbite $\overline{T}_{\sigma}$ est contenue
dans la carte affine $\overline{T}(\sigma )$ et est en fait l'unique
orbite ferm\'{e}e pour l'action de $T$ sur cette vari\'{e}t\'{e}
torique affine. L'orbite $\overline{T}_{\sigma}$ est dans
l'adh\'{e}rence de l'orbite $\overline{T}_{\tau}$ si et seulement si
$\tau$ est une face de $\sigma$. En particulier, $\overline{T}_{\{0\}}
=T$ est l'unique orbite ouverte dans $\overline{T}$. Chaque orbite
$\overline{T}_{\sigma}$ a par construction un point marqu\'{e}
$\overline{t}_{\sigma}$.

Le {\it champ torique} associ\'{e} \`{a} l'\'{e}ventail $E$ est le
champ alg\'{e}brique quotient
$$
\overline{{\cal T}}=[\overline{T}/T].
$$
Ce champ alg\'{e}brique, qui est normal, de type fini et de dimension
$0$, n'a qu'un nombre fini de points, \`{a} savoir les points
d\'{e}finis par les $\overline{t}_{\sigma}\in\overline{T}$.
La gerbe r\'{e}siduelle d'un tel point est le classifiant du
fixateur de $\overline{t}_{\sigma}$ dans $T$.
\vskip 3mm

Soient maintenant $k$ un corps et $n$ un entier $\geq 0$.
Consid\'{e}rons l'espace affine ${\Bbb R}^{r,n}= \{x=(x_{0},\ldots
,x_{n})\in {\Bbb R}^{n+1}\mid x_{0}+x_{1}+\cdots +x_{n}=r\}$ et son
r\'{e}seau ${\Bbb Z}^{r,n}=\{x=(x_{0},\ldots ,x_{n})\in {\Bbb
Z}^{n+1}\mid x_{0}+x_{1}+\cdots +x_{n}=r\}$. On a le simplexe
$$
S_{{\Bbb R}}^{r,n}=\{x=(x_{0},\ldots ,x_{n})\in ({\Bbb R}_{\geq 
0})^{n+1}\mid x_{0}+x_{1}+\cdots +x_{n}=r\}\subset {\Bbb R}^{r,n}
$$
et l'ensemble de ses points entiers
$$
S^{r,n}={\Bbb Z}^{r,n}\cap S_{{\Bbb R}}^{r,n}.
$$
On appelle {\it pav\'{e} entier} toute partie $P$ de $S_{{\Bbb
R}}^{r,n}$ d'int\'{e}rieur $P^{\circ}$ non vide et de la forme
$$
P=\{x\in {\Bbb R}^{r,n}\mid\sum_{j\in J}^{}x_{j}\geq d_{J},~\forall J\}
$$
pour $(d_{J})$ une famille d'entiers index\'{e}s par les parties $J$ 
de $\{0,1,\ldots ,n\}$ qui est convexe au sens o\`{u}
$$
d_{\emptyset}=0,~d_{\{0,1,\ldots ,n\}}=r\hbox{ et }d_{J'}+d_{J''}\leq 
d_{J'\cup J''}+d_{J'\cap J''},~\forall J',J''.
$$
On appelle {\it pavage entier} de $S_{{\Bbb R}}^{r,n}$ toute famille
${\bf P}$ de pav\'{e}s entiers telle que
$$
S_{{\Bbb R}}^{r,n}=\bigcup_{P\in {\bf P}}^{}P
$$
et que $P_{1}^{\circ}\cap P_{2}^{\circ}=\emptyset$ pour tous 
$P_{1}\not= P_{2}$ dans ${\bf P}$.

Une application $S^{r,n}\rightarrow {\Bbb R}$ est dite affine si elle
est la restriction d'une application affine ${\Bbb R}^{r,n}\rightarrow
{\Bbb R}$. Soient $Y_{{\Bbb R}}^{r,n}$ le ${\Bbb R}$-espace vectoriel
quotient de ${\Bbb R}^{S^{r,n}}=\{S^{r,n}\rightarrow {\Bbb R}\}$ par
le sous-espace des applications affines, et $Y^{r,n}$ l'image dans ce
quotient du r\'{e}seau ${\Bbb Z}^{S^{r,n}}\subset {\Bbb R}^{S^{r,n}}$.

Pour tout pavage entier ${\bf P}$ de $S_{{\Bbb R}}^{r,n}$, notons
$\sigma_{{\bf P}}\subset Y_{{\Bbb R}}$ l'ensemble des classes $y$ de
fonctions $S^{r,n} \rightarrow {\Bbb R}$ qui, pour chaque pav\'{e}
$P\in {\bf P}$, admettent un repr\'{e}sentant
$y_{P}:S^{r,n}\rightarrow {\Bbb R}$ \`{a} valeurs positives tel que
$$
P=\{x\in S_{{\Bbb R}}^{r,n}\mid y_{P}(x)=0\}.
$$
Le pavage entier ${\bf P}$ est dit {\it admissible} si $\sigma_{P}$
est non vide.

\thm LEMME (Lafforgue)
\enonce
Pour chaque pavage entier admissible ${\bf P}$ de $S_{{\Bbb
R}}^{r,n}$, $\sigma_{{\bf P}}$ est un c\^{o}ne convexe poly\'{e}dral
rationnel dans $Y_{{\Bbb R}}$. De plus, $\sigma_{{\bf Q}}$ est une
face de $\sigma_{{\bf P}}$ si et seulement le pavage entier admissible
${\bf Q}$ est plus grossier que ${\bf P}$.

La famille des $\sigma_{{\bf P}}$, pour ${\bf P}$ parcourant
l'ensemble des pavages entiers admissibles de $S_{{\Bbb R}}^{r,n}$,
est un \'{e}ventail.
\endthm

On appellera {\it champ des pavages} le $k$-champ torique
$$
{\cal T}^{r,n}=[\overline{T}^{r,n}/T^{r,n}]
$$
associ\'{e} \`{a} cet \'{e}ventail. On v\'{e}rifie que l'on a une 
suite exacte
$$
1\rightarrow {\Bbb G}_{{\rm m},k}\rightarrow {\Bbb G}_{{\rm
m},k}^{n+1}\times {\Bbb G}_{{\rm m},k}\rightarrow {\Bbb G}_{{\rm
m},k}^{S^{r,n}}\rightarrow T^{r,n}\rightarrow 1
$$
o\`{u} la deuxi\`{e}me fl\`{e}che envoie $z$ sur $(z,\ldots ,z;z^{r})$
et la troisi\`{e}me est donn\'{e}e par
$$
(u_{0},u_{1},\ldots ,u_{n};z)\mapsto (u_{0}^{i_{0}}
u_{1}^{i_{1}}\cdots u_{n}^{i_{n}}z^{-1})_{i\in S^{r,n}}.
$$

Par construction, le $k$-sch\'{e}ma $\overline{T}^{r,n}$ est normal de
type fini et l'action de $T^{r,n}$ sur $\overline{T}^{r,n}$ n'a qu'un
nombre fini d'orbites. Ces orbites sont index\'{e}es par les pavages
entiers admissibles de $S^{r,n}$ et chacune a un point marqu\'{e}
$\overline{t}_{{\bf P}}$. Le champ des pavages ${\cal T}^{r,n}$ est
donc alg\'{e}brique normal de type fini et n'a qu'un nombre fini de
points d\'{e}finis par les points $\overline{t}_{{\bf P}}$. Le point
d\'{e}fini par $\overline{t}_{{\bf P}}$ est dans l'adh\'{e}rence du
point d\'{e}fini par $\overline{t}_{{\bf Q}}$ si le pavage entier
admissible ${\bf P}$ raffine le pavage entier admissible ${\bf Q}$. Le
point ouvert correspond au pavage le plus grossier, c'est-\`{a}-dire
au pavage trivial (un seul pav\'{e}).

\rem Exemples
\endrem
{\rm (i)} Pour $n=0$, $S_{{\Bbb R}}^{r,n}=S^{r,n}=\{0\}$ et le champ
des pavages est simplement $\mathop{\rm Spec}(k)$. Pour $n=1$ on peut
identifier $S_{{\Bbb R}}^{r,n}$ \`{a} l'intervalle $[0,r]\subset {\Bbb
R}$ et les pav\'{e}s entiers sont pr\'{e}cis\'{e}ment les intervalles
$[i,j]$ o\`{u} $i<j$ sont des entiers compris entre $0$ et $r$. Un
pavage entier est donc une d\'{e}composition
$$
[0,r]=[0,r_{1}]\cup [r_{1},r_{1}+r_{2}]\cup\cdots\cup [r_{1}+\cdots 
+r_{s-1},r_{1}+\cdots +r_{s-1}+r_{s}]
$$
de l'intervalle $[0,r]$ o\`{u} $r_{1}+\cdots +r_{s}=r$ est une
d\'{e}composition de $r$ en entiers strictement positifs. Le champ des
pavages n'est autre dans ce cas que le champ torique $[{\Bbb
A}_{k}^{r-1}/{\Bbb G}_{{\rm m},k}^{r-1}]$.

\decale{\rm (ii)} Pour $n=2$ on peut identifier $S_{{\Bbb R}}^{r,n}$
\`{a} un triangle \'{e}quilat\'{e}ral dans ${\Bbb R}^{2}$ dont les
trois c\^{o}t\'{e}s sont de longueur $r$. On num\'{e}rote les sommets
$s_{0}=s_{3},s_{1},s_{2}$ de ce triangle et, pour $i=0,1,2$, on
d\'{e}signe par $x_{i}$ la distance au sommet $s_{i}$ d'un point du
c\^{o}t\'{e} $s_{i}s_{i+1}$. Par tout point de $S_{{\Bbb R}}^{r,n}$
passent trois droites $D_{12}$, $D_{20}$ et $D_{01}$ parall\`{e}les
aux c\^{o}t\'{e}s du triangle $s_{1}s_{2}$, $s_{2}s_{0}$ et
$s_{0}s_{1}$, et on peut rep\'{e}rer un tel point par les
coordonn\'{e}es $x_{0}$, $x_{1}$ et $x_{2}$ des intersections
$D_{12}\cap s_{0}s_{1}$, $D_{20}\cap s_{1}s_{2}$ et $D_{01}\cap
s_{2}s_{0}$. Par d\'{e}finition, $S^{r,n}\subset S_{{\Bbb R}}^{r,n}$
est form\'{e} des points \`{a} coordonn\'{e}es $(x_{0},x_{1},x_{2})$
enti\`{e}res. 

Appelons droite enti\`{e}re toute droite qui est parall\`{e}le \`{a}
l'un des trois c\^{o}t\'{e} du triangle $S_{{\Bbb R}}^{r,n}$ et qui
passe par un point de $S^{r,n}$; appelons segment entier tout segment
d'une droite enti\`{e}re dont les deux extr\'{e}mit\'{e}s sont des
points de $S^{r,n}$.

Si l'on trace toutes les droites enti\`{e}res, on obtient un pavage
constitu\'{e} de petits triangles \'{e}quilat\'{e}raux de
c\^{o}t\'{e}s de longueur $1$ qui est entier admissible. Il est plus
fin que tout autre pavage entier: tout pav\'{e} entier est une
r\'{e}union d'une famille de ces petits triangles \'{e}quilat\'{e}raux
et est un hexagone, \'{e}ventuellement d\'{e}g\'{e}n\'{e}r\'{e}, dont
les c\^{o}t\'{e}s sont des segments entiers. La vari\'{e}t\'{e}
torique $T^{r,n}\hookrightarrow \overline{T}^{r,n}$ est donc affine
(Cette derni\`{e}re propri\'{e}t\'{e} n'est plus satisfaite d\`{e}s
que $n\geq 3$).

Des deux pavages entiers ci-dessous, seul celui de gauche (qui joue un
r\^{o}le dans la compactification de l'isog\'{e}nie de Lang) est
admissible.
$$
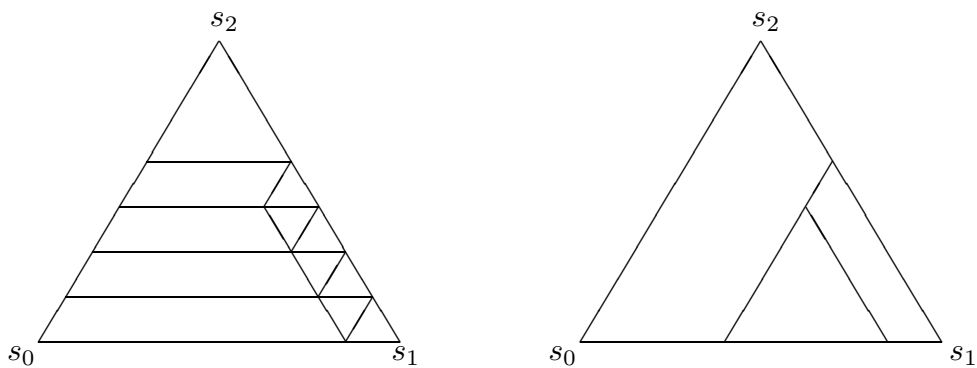
\diagram{
\noalign{\vskip 35mm}
\rule(-50,0)\dir(1,0)\long{40}
\rule(-50,0)\dir(3,5)\long{20}
\rule(-11.1,0)\dir(-3,5)\long{20}
\rule(-47,5)\dir(1,0)\long{34}
\rule(-44,10)\dir(1,0)\long{28}
\rule(-41,15)\dir(1,0)\long{22}
\rule(-38,20)\dir(1,0)\long{16}
\rule(-17.1,0)\dir(-3,5)\long{9}
\rule(-16,0)\dir(3,5)\long{3}
\rule(-19,5)\dir(3,5)\long{3}
\rule(-22,10)\dir(3,5)\long{3}
\rule(-25,15)\dir(3,5)\long{3}
\put (-53.5,-2){s_{0}}
\put (-10.9,-2){s_{1}}
\put (-31,35){s_{2}}
\rule(10,0)\dir(1,0)\long{40}
\rule(10,0)\dir(3,5)\long{20}
\rule(48.9,0)\dir(-3,5)\long{20}
\rule(42.9,0)\dir(-3,5)\long{9}
\rule(26,0)\dir(3,5)\long{12}
\put (6.5,-2){s_{0}}
\put (50.9,-2){s_{1}}
\put (29,35){s_{2}}
}
$$

\decale{\rm (iii)} Si $r=2$, le champ torique ${\cal T}^{r,n}$ est
lisse, mais ce n'est plus le cas d\`{e}s que $n\geq 2$ et $r\geq 3$
\hfill\hfill$\square$
\vskip 3mm

Pour chaque entier $0\leq m<n$ et chaque application injective
$\{0,1,\ldots ,m\}\hookrightarrow \{0,1,\ldots ,n\}$, on a une
identification induite de $S_{{\Bbb R}}^{r,m}$ \`{a} une face de
$S_{{\Bbb R}}^{r,n}$ et tout pavage entier admissible de $S_{{\Bbb
R}}^{r,n}$ induit par restriction un pavage entier admissible de
$S_{{\Bbb R}}^{r,m}$. Dualement, on en d\'{e}duit des projections
$T^{r,n}\rightarrow T^{r,m}$, $\overline{T}^{r,n} \rightarrow
\overline{T}^{r,m}$ et ${\cal T}^{r,n}\rightarrow {\cal T}^{r,m}$.
\vskip 5mm

{\bf 4.2. Compactifications des $\mathop{\rm PGL}\nolimits_{r}^{n+1}/
\mathop{\rm PGL}\nolimits_{r}$}
\vskip 5mm

Consid\'{e}rons un espace vectoriel $V$ de dimension $r$ sur un corps
$k$, que l'on verra aussi comme un $k$-sch\'{e}ma affine isomorphe
\`{a} ${\Bbb A}_{k}^{r}$, et notons $G=\mathop{\rm GL}(V)$ le
$k$-sch\'{e}ma des automorphismes lin\'{e}aires de $V$. Le
$k$-sch\'{e}ma en groupes $G^{n+1}$ agit de mani\`{e}re \'{e}vidente
sur $V^{n+1}$ et donc sur $\bigwedge_{}^{r}V^{n+1}$.

On peut d\'{e}composer la repr\'{e}sentation $\bigwedge_{}^{r}V^{n+1}$
de $G^{n+1}$ en la somme directe de repr\'{e}sentations
irr\'{e}ductibles
$$
\bigwedge^{r}V^{n+1}=\bigoplus_{i\in S^{r,n}}\bigwedge^{i_{0}}V
\otimes\bigwedge^{i_{1}}V\otimes\cdots\otimes \bigwedge^{i_{n}}V
$$
o\`{u} $S^{r,n}$ est l'ensemble des points entiers du simplexe
introduit dans la sous-section pr\'{e}c\'{e}dente. Le tore ${\Bbb G}_{{\rm
m},k}^{S^{r,n}}$ agit sur $\bigwedge^{r}V^{n+1}$ par
$$
(t=(t_{i})_{i\in S^{r,n}},x=\oplus_{i\in S^{r,n}}x_{i})\mapsto
\oplus_{i\in S^{r,n}}t_{i}x_{i},
$$
et cette action commute \`{a} celle de $G^{n+1}$. On en d\'{e}duit une
action du $k$-sch\'{e}ma en groupes produit $G^{n+1}\times {\Bbb G}_{{\rm
m},k}^{S^{r,n}}$ sur l'espace projectif ${\Bbb P}(\bigwedge_{}^{r}
V^{n+1})$ des droites de $\bigwedge_{}^{r}V^{n+1}$.

On envoie le produit $G\times {\Bbb G}_{{\rm m},k}^{n+1}\times {\Bbb
G}_{{\rm m},k}$ dans $G^{n+1}\times {\Bbb G}_{{\rm m},k}^{S^{r,n}}$ par
$$
(g;u_{0},u_{1},\ldots ,u_{n};z)\mapsto (u_{0}^{-1}g,u_{1}^{-1}g,\ldots
,u_{n}^{-1}g;(u_{0}^{i_{0}} u_{1}^{i_{1}}\cdots u_{n}^{i_{n}}
z^{-1})_{i\in S^{r,n}}).
$$
L'homomorphisme ainsi d\'{e}fini n'est pas injectif (il admet pour
noyau ${\Bbb G}_{{\rm m},k}$ plong\'{e} par $z\mapsto (z\mathop{\rm
Id}\nolimits_{V};z,\ldots ,z,z^{r})$), mais on notera quand m\^{e}me
$(G^{n+1}\times {\Bbb G}_{{\rm m},k}^{S^{r,n}})/(G\times {\Bbb
G}_{{\rm m},k}^{n+1}\times {\Bbb G}_{{\rm m},k})$ le quotient de
$G^{n+1}\times {\Bbb G}_{{\rm m},k}^{S^{r,n}}$ par l'image de cet
homomorphisme.

Le plongement diagonal $V\hookrightarrow V^{n+1}$ induit un plongement
$$
\bigwedge^{r}V\hookrightarrow\bigwedge_{}^{r}V^{n+1}
$$
et d\'{e}finit un point $\omega\in {\Bbb P}(\bigwedge_{}^{r}V^{n+1})$
dont le fixateur est pr\'{e}cis\'{e}ment l'image de $G\times {\Bbb
G}_{{\rm m},k}^{n+1}\times {\Bbb G}_{{\rm m},k}$ dans $G^{n+1}\times
{\Bbb G}_{{\rm m},k}^{S^{r,n}}$. On a donc un plongement localement
ferm\'{e} de $k$-sch\'{e}mas
$$
(G^{n+1}\times {\Bbb G}_{{\rm m},k}^{S^{r,n}})/(G\times {\Bbb G}_{{\rm
m},k}^{n+1}\times {\Bbb G}_{{\rm m},k})\hookrightarrow {\Bbb P}
\left(\bigwedge_{}^{r} V^{n+1}\right),~(g,t)\mapsto (g,t)\cdot\omega .
$$

Par construction l'image de ce plongement est contenue dans l'ouvert
$$
\left[\prod_{i\in S^{r,n}}\left(\left(\bigwedge^{i_{0}}V\otimes
\bigwedge^{i_{1}}V\otimes\cdots\otimes\bigwedge^{i_{n}}V\right)
-\{0\}\right)\right]\Big/{\Bbb G}_{{\rm m},k}\subset{\Bbb P}
\left(\bigwedge_{}^{r}V^{n+1}\right).
$$
des droites de $\bigwedge_{}^{r} V^{n+1}$ dont la projection sur
chaque facteur direct $\bigwedge^{i_{0}}V\otimes\bigwedge^{i_{1}}
V\otimes\cdots\otimes\bigwedge^{i_{n}}V$ n'est pas nulle.

Consid\'{e}rons alors le morphisme
$$\displaylines{
(G^{n+1}\times {\Bbb G}_{{\rm m},k}^{S^{r,n}})/(G\times {\Bbb G}_{{\rm
m},k}^{n+1}\times {\Bbb G}_{{\rm m},k})
\hfill\cr\hfill
\hookrightarrow\overline{T}^{r,n}\times_{k}\left[\prod_{i\in S^{r,n}}
\left(\left(\bigwedge^{i_{0}}V\otimes \bigwedge^{i_{1}}V\otimes
\cdots\otimes\bigwedge^{i_{n}}V\right) -\{0\}\right)\right]\Big/
{\Bbb G}_{{\rm m},k}}
$$
dont les composantes sont le compos\'{e} de la projection
$$
(G^{n+1}\times_{k}{\Bbb G}_{{\rm m},k}^{S^{r,n}})/(G\times {\Bbb
G}_{{\rm m},k}^{n+1}\times {\Bbb G}_{{\rm m},k})\twoheadrightarrow
{\Bbb G}_{{\rm m},k}^{S^{r,n}}/({\Bbb G}_{{\rm m},k}^{n+1}\times
{\Bbb G}_{{\rm m},k})=T^{r,n}
$$
et de l'immersion ouverte
$$
T^{r,n}\hookrightarrow\overline{T}^{r,n},
$$
et l'immersion localement ferm\'{e}e ci-dessus. Ce morphisme est une
immersion localement ferm\'{e}e. Notons
$$
\Omega^{n}(V)\subset\left(\left[\prod_{i\in S^{r,n}}\left(\left(
\bigwedge^{i_{0}}V\otimes\bigwedge^{i_{1}}V\otimes\cdots\otimes
\bigwedge^{i_{n}}V\right)-\{0\}\right)\right]\Big/{\Bbb G}_{{\rm m},k}
\right)\times_{k}\overline{T}^{r,n}
$$
l'adh\'{e}rence sch\'{e}matique de son image. et
$$
\overline{\Omega}^{n}(V)\rightarrow\prod_{i\in S^{r,n}}{\Bbb P}\left(
\bigwedge^{i_{0}}V\otimes \bigwedge^{i_{1}}V\otimes\cdots
\otimes\bigwedge^{i_{n}}V\right)
$$
le quotient de $\Omega^{n}(V)$ par l'action libre du tore ${\Bbb
G}_{{\rm m},k}^{S^{r,n}}/{\Bbb G}_{{\rm m},k}$ (o\`{u} ${\Bbb G}_{{\rm
m},k}$ est plong\'{e} diagonalement dans ${\Bbb G}_{{\rm
m},k}^{S^{r,n}}$).

\thm TH\'{E}OR\`{E}ME 
\enonce
{\rm (i)} La premi\`{e}re projection 
$$
\Omega^{n}(V)\rightarrow\left[\prod_{i\in S^{r,n}}\left(\left(
\bigwedge^{i_{0}}V\otimes \bigwedge^{i_{1}}V\otimes\cdots\otimes
\bigwedge^{i_{n}}V\right)-\{0\}\right)\right]\Big/{\Bbb G}_{{\rm m},k}
$$
et son quotient par les actions libres du tore ${\Bbb G}_{{\rm
m},k}^{S^{r,n}}/{\Bbb G}_{{\rm m},k}$
$$
\overline{\Omega}^{n}(V)\rightarrow\prod_{i\in S^{r,n}}{\Bbb P}\left(
\bigwedge^{i_{0}}V\otimes \bigwedge^{i_{1}}V\otimes\cdots
\otimes\bigwedge^{i_{n}}V\right)
$$
sont projectives.

\decale{\rm (ii)} La seconde projection
$$
\Omega^{n}(V)\rightarrow\overline{T}^{r,n},
$$
qui est \'{e}quivariante relativement au morphisme de tores ${\Bbb
G}_{{\rm m},k}^{S^{r,n}}/{\Bbb G}_{{\rm m},k}\twoheadrightarrow
T^{r,n}$, est lisse purement de dimension relative $nr^{2}$; le
morphisme de champs qu'on en d\'{e}duit par passage au quotient
$$
\overline{\Omega}^{n}(V)=\Omega^{n}(V)/({\Bbb
G}_{{\rm m},k}^{S^{r,n}}/{\Bbb G}_{{\rm m},k})\rightarrow 
[\overline{T}^{r,n}/T^{r,n}]={\cal T}^{r,n}
$$
est lisse purement de dimension relative $n(r^{2}-1)$.

Plus g\'{e}n\'{e}ralement, pour tout entier $0\leq m<n$ et toute
application injective $\{0,1,\ldots ,m\}\hookrightarrow \{0,1,\ldots
,n\}$, l'inclusion $S^{r,m}\hookrightarrow S^{r,n}$ et la projection
$G^{n+1}\twoheadrightarrow G^{m+1}$ correspondantes induisent un
\'{e}pimorphisme de tores $T^{r,n}\twoheadrightarrow T^{r,m}$, un
prolongement torique $\overline{T}^{r,n}\rightarrow
\overline{T}^{r,m}$ de cet \'{e}pimorphisme et un morphisme de
sch\'{e}mas $\Omega^{n}(V) \rightarrow\Omega^{m}(V)$ dit {\rm de face}
au-dessus de ce morphisme torique, et les morphismes
$$
\Omega^{n}(V)\rightarrow\Omega^{m}(V)\times_{\overline{T}^{r,m}}
\overline{T}^{r,n}\hbox{ et }\,\overline{\Omega}^{n}(V)\rightarrow
\overline{\Omega}^{m}(V)\times_{{\cal T}^{r,m}}{\cal T}^{r,n}
$$
que l'on en d\'{e}duit, sont lisses purement de dimensions relatives
$(n-m)r^{2}$ et $(n-m)(r^{2}-1)$ respectivement.
\endthm

Le $k$-sch\'{e}ma $\overline{\Omega}^{n}(V)$ est donc une
compactification projective et toro\"{\i}dale de
$\overline{G}^{n+1}/\overline{G}$ o\`{u} on a not\'{e}
$\overline{G}=\mathop{\rm PGL}(V)$. Les $\overline{\Omega}^{n}(V)$ et
les morphismes de face $\overline{\Omega}^{n}(V)\rightarrow
\overline{\Omega}^{m}(V)$ sont sous-jacents \`{a} un $k$-sch\'{e}ma
simplicial $\Omega^{\bullet}(V)$ qui contient comme ouvert dense le
$k$-sch\'{e}ma simplicial classifiant $\mathop{\rm
B}(\overline{G})=\overline{G}^{\bullet +1}/\overline{G}$.

\rem Remarques 
\endrem
{\rm (i)} Pour $n=0$, $\Omega^{n}(V)$ n'est autre que $\mathop{\rm
Spec}(k)$. Pour $n=1$, le carr\'{e} commutatif
$$\diagram{
(\mathop{\rm GL}(V)^{2}\times_{k}{\Bbb G}_{{\rm
m},k}^{r+1})/(\mathop{\rm GL}(V)\times_{k}{\Bbb G}_{{\rm
m},k}^{2}\times_{k}{\Bbb G}_{{\rm m},k})&\kern
-1mm\smash{\mathop{\hbox to 8mm{\rightarrowfill}}
\limits^{\scriptstyle \sim}}\kern -1mm&\mathop{\rm
GL}(V)\times_{k}{\Bbb G}_{{\rm m},k}^{r-1}\cr
\llap{$\scriptstyle $}\left\downarrow
\vbox to 4mm{}\right.\rlap{}&+&\llap{}\left\downarrow
\vbox to 4mm{}\right.\rlap{$\scriptstyle $}\cr
{\Bbb G}_{{\rm m},k}^{r+1}/({\Bbb G}_{{\rm m},k}^{2}\times_{k}{\Bbb
G}_{{\rm m},k})&\kern -20mm\smash{\mathop{\hbox to 36mm
{\rightarrowfill}} \limits^{\scriptstyle \sim}}\kern -10mm&{\Bbb
G}_{{\rm m},k}^{r-1}\cr}
$$
o\`{u} les deux fl\`{e}ches horizontales
$$
(g_{0},g_{1},(\mu_{\rho})_{\rho =0,\ldots ,r})\mapsto
\left({\mu_{1}\over\mu_{0}}g_{0}g_{1}^{-1},{\mu_{0}\mu_{2}\over
\mu_{1}^{2}},\ldots ,{\mu_{r-2}\mu_{r}\over \mu_{r-1}^{2}}\right)
$$
et
$$
(\mu_{\rho})_{\rho =0,1,\ldots ,r}\mapsto\left({\mu_{0}\mu_{2}\over
\mu_{1}^{2}},\ldots ,{\mu_{r-2}\mu_{r}\over \mu_{r-1}^{2}}\right)
$$
sont des isomorphismes, se prolonge en un carr\'{e} commutatif
$$\diagram{
\Omega^{1}(V)&\kern -1mm\smash{\mathop{\hbox to 8mm{\rightarrowfill}}
\limits^{\scriptstyle \sim}}\kern -1mm&\Omega (V,V)\cr
\llap{$\scriptstyle $}\left\downarrow
\vbox to 4mm{}\right.\rlap{}&+&\llap{}\left\downarrow
\vbox to 4mm{}\right.\rlap{$\scriptstyle $}\cr
\overline{T}^{r,1}&\kern -3mm\smash{\mathop{\hbox to
12mm{\rightarrowfill}} \limits^{\scriptstyle \sim}}\kern -3mm&{\Bbb
A}_{k}^{r-1}\cr}
$$
o\`{u}, de nouveau, les deux fl\`{e}ches horizontales sont des
isomorphismes. (La premi\`{e}re fl\`{e}che verticale est celle du
th\'{e}or\`{e}me ci-dessus et la seconde $\Omega (V,V)\rightarrow
{\Bbb A}_{k}^{r-1}$ a \'{e}t\'{e} construite dans la section 3.) En
particulier, on peut identifier le $k$-sch\'{e}ma
$\mathop{\widetilde{\rm End}}(V)$ des endomorphismes complets de $V$
\`{a} un quotient de $\Omega^{1}(V)$ par un tore ${\Bbb G}_{{\rm
m},k}^{r}$ et, en divisant en plus par les homoth\'{e}ties, on en
d\'{e}duit un isomorphisme
$$
\overline{\Omega}^{1}(V)\buildrel\sim\over\longrightarrow
\mathop{\widetilde{\rm End}}(V)/{\Bbb G}_{{\rm m},k}
$$
qui \'{e}change les projections
$$
\overline{\Omega}^{1}(V)\rightarrow {\cal T}^{r,n}\hbox{ et
}\mathop{\widetilde{\rm End}}(V)/ {\Bbb G}_{{\rm m},k}\rightarrow
[{\Bbb A}_{k}^{r-1}/{\Bbb G}_{{\rm m},k}^{r-1}].
$$
Le $k$-sch\'{e}ma projectif $\overline{\Omega}^{1}(V)$ est donc la
compactification de De Concini et Procesi de $\mathop{\rm PGL}(V)$.

\decale{\rm (ii)} Le morphisme $\overline{\Omega}^{n}(V)\rightarrow
(\overline{\Omega}^{1}(V))^{n+1}$ produit des morphismes induits par
les applications injectives $\{0,1\}\hookrightarrow\{0,1,\ldots ,n\}$
donn\'{e}es par $\{0,1\}\mapsto\{0,1\}, \{0,1\}\mapsto\{1,2\},\ldots
,\{0,1\}\mapsto\{n,0\}$, envoie $\overline{G}^{n+1}/\overline{G}$ sur
le ferm\'{e} d'\'{e}quation $g_{0}g_{1}\cdots g_{n}=1$ dans
$(\overline{G}^{2}/\overline{G})^{n+1} \cong \overline{G}^{n+1}$.
Ainsi, $\overline{\Omega}^{n}(V)$ r\'{e}alise une
{\og}{compactification}{\fg} du morphisme de multiplication de $n$
\'{e}l\'{e}ments dans $\overline{G}$.
\hfill\hfill$\square$
\vskip 3mm

Rappelons que le $k$-sch\'{e}ma $\overline{T}^{r,n}$ est stratifi\'{e}
par les orbites de $T^{r,n}$ et que ces orbites sont index\'{e}es par
les pavages entiers admissibles ${\bf P}$ de $S_{{\Bbb R}}^{r,n}$ et
poss\`{e}dent chacune un point marqu\'{e} $\overline{t}_{{\bf P}}\in
\overline{T}^{r,n}$. L'image r\'{e}ciproque par la projection
$\Omega^{n}(V)\rightarrow \overline{T}^{r,n}$ de cette stratification
est \'{e}videmment une stratification de $\Omega^{n}(V)$ par des
parties localement ferm\'{e}es $\Omega_{{\bf P}}^{n}(V)$, index\'{e}es
elles aussi par les pavages entiers admissibles ${\bf P}$ de $S_{{\Bbb
R}}^{r,n}$. L'inclusion dans la strate $\Omega_{{\bf P}}^{n}(V)$ de la
fibre $\Omega^{n}(V)_{\overline{t}_{{\bf P}}}$ de $\Omega^{n}(V)
\rightarrow \overline{T}^{r,n}$ en $\overline{t}_{{\bf P}}$ induit un
isomorphisme de sch\'{e}mas
$$
(\Omega^{n}(V)_{\overline{t}_{{\bf P}}}\times_{k}{\Bbb G}_{{\rm
m},k}^{S^{r,n}})/({\Bbb G}_{{\rm m},k}^{S^{r,n}})_{\overline{t}_{{\bf
P}}}\buildrel\sim\over\longrightarrow \Omega_{{\bf P}}^{n}(V)
$$
o\`{u} $({\Bbb G}_{{\rm m},k}^{S^{r,n}})_{\overline{t}_{{\bf P}}}$ est
le stabilisateur de $\overline{t}_{{\bf P}}$ dans ${\Bbb G}_{{\rm
m},k}^{S^{r,n}}$.

En particulier, si ${\bf P}$ est le pavage trivial,
$\Omega^{n}(V)_{\overline{t}_{{\bf P}}}=G^{n+1}/G$ et $\Omega_{{\bf
P}}^{n}$ est l'ouvert $(G^{n+1}\times {\Bbb G}_{{\rm
m},k}^{S^{r,n}})/(G\times {\Bbb G}_{{\rm m},k}^{n+1}\times {\Bbb
G}_{{\rm m},k})$ de $\Omega^{n}(V)$.

Lafforgue appelle {\it $k$-sch\'{e}ma des graphes recoll\'{e}s dans
$V^{n+1}$} associ\'{e} \`{a} un pavage entier admissible ${\bf P}$ de
$S_{{\Bbb R}}^{r,n}$, et note $\mathop{\rm Gr}\nolimits_{{\bf
P}}^{n}(V)$, le sch\'{e}ma de modules des familles $(W_{P})_{P\in {\bf
P}}$ de sous-espaces de dimension $r$ dans $V^{n+1}$ index\'{e}s par
les pav\'{e}s $P$ de ${\bf P}$ telles que:
\vskip 2mm

\itemitem{-} pour tout pav\'{e} $P\in {\bf P}$ et tout
$J\subset\{0,1,\ldots ,n\}$, on a
$$
\mathop{\rm dim}(W_{P}\cap V^{J})=\mathop{\rm min}_{i\in
P\cap S^{r,n}}\left\{\sum_{j\in J}i_{j}\right\},
$$
\vskip 2mm

\itemitem{-} pour tous pav\'{e}s $P',P''\in {\bf P}$ ayant en commun 
un bord d'\'{e}quation $\sum_{j\in J}x_{j}=d_{J}$ pour un 
certain $J\subset\{0,1,\ldots ,n\}$ et
$$
d_{J}=\mathop{\rm min}_{i'\in P'\cap S^{r,n}}\left\{\sum_{j\in
J}i_{j}'\right\} =\mathop{\rm max}_{i''\in P''\cap
S^{r,n}}\left\{\sum_{j\in J}i_{j}''\right\},
$$
on a
$$
i_{J}^{-1}(W_{P'})=\mathop{\rm pr}\nolimits_{J}(W_{P''})\hbox{ et
}\mathop{\rm pr}\nolimits_{J^{{\rm c}}}(W_{P'})=i_{J^{{\rm c}}}^{-1}(W_{P''})
$$
\itemitem{} o\`{u} $J^{{\rm c}}$ est le compl\'{e}mentaire de $J$ dans
$\{0,1,\ldots ,n\}$ et $i_{J}:V^{J}\hookrightarrow V^{n+1}$,
$i_{J^{{\rm c}}}:V^{J^{{\rm c}}}\hookrightarrow V^{n+1}$ et
$\mathop{\rm pr}\nolimits_{J}: V^{n+1}\twoheadrightarrow V^{J}$,
$\mathop{\rm pr}\nolimits_{J^{{\rm c}}}: V^{n+1}\twoheadrightarrow
V^{J^{{\rm c}}}$ sont les inclusions et projections canoniques.

\thm TH\'{E}OR\`{E}ME 
\enonce
Pour chaque pavage entier admissible ${\bf P}$ de $S_{{\Bbb
R}}^{r,n}$, les coordonn\'{e}es de Pl\"{u}cker induisent un
isomorphisme du $k$-sch\'{e}ma des graphes recoll\'{e}s $\mathop{\rm
Gr}\nolimits_{{\bf P}}^{n}(V)$ sur la fibre
$\Omega^{n}(V)_{\overline{t}_{{\bf P}}}$.
\endthm
\vskip 5mm

{\bf 4.3. Compactification de l'isog\'{e}nie de Lang}
\vskip 5mm

Le corps de base est de nouveau ${\Bbb F}_{q}$. On notera encore
simplement $\tau$ l'endomorphisme de Frobenius d'un sch\'{e}ma.

Rappelons que l'on a d\'{e}fini un \'{e}ventail $(\sigma_{{\bf
P}})_{{\bf P}}$ dans l'espace vectoriel r\'{e}el $Y_{{\Bbb R}}^{r,2}$
quotient de ${\Bbb R}^{S^{r,2}}=\{S^{r,2}\rightarrow {\Bbb R}\}$ par
le sous-espace des applications affines et que l'on a not\'{e}
$T^{r,2}\hookrightarrow\overline{T}^{r,2}$ la vari\'{e}t\'{e} torique
correspondante.

Consid\'{e}rons maintenant le sous-espace de ${\Bbb R}^{S^{r,2}}$
form\'{e} des applications
$$
y:S^{r,2}=\{(i_{0},i_{1},i_{2})\in {\Bbb N}^{3}\mid 
i_{0}+i_{1}+i_{2}=r\}\rightarrow {\Bbb R}
$$
telles que
$$
y(0,i_{1},i_{2})=qy(i_{1},0,i_{2}),~\forall (0,i_{1},i_{2})\in 
S^{r,2},
$$
et son image $Y_{{\Bbb R}}^{r,\tau}$ dans $Y_{{\Bbb R}}^{r,2}$. On 
remarquera que la donn\'{e}e d'une telle application $y$ \'{e}quivaut 
\`{a} la donn\'{e}e de sa restriction \`{a}
$$
S^{r,\tau}=\{(i_{0},i_{1},i_{2})\in {\Bbb N}^{3}\mid 
i_{0}+i_{1}+i_{2}=r\hbox{ et }i_{0}\not=0\}\subset S^{r,2}
$$
Appelons {\it pavage entier $q$-admissible} tout pavage ${\bf P}$ de
$S_{{\Bbb R}}^{r,2}$ tel que le c\^{o}ne convexe poly\'{e}dral
rationnel $\sigma_{{\bf P}}\subset Y_{{\Bbb R}}^{r,2}$ rencontre le
sous-espace $Y_{{\Bbb R}}^{r,\tau}$. La famille des intersections
$\sigma_{{\bf P}}^{\tau}:=\sigma_{{\bf P}}\cap Y_{{\Bbb R}}^{r,\tau}$
pour ${\bf P}$ parcourant l'ensemble des pavages entiers $q$-admissibles
est un \'{e}ventail dans $Y_{{\Bbb R}}^{r,\tau}$ et d\'{e}finit donc
une vari\'{e}t\'{e} torique normale
$$
T^{r,\tau}\hookrightarrow \overline{T}^{r,\tau}.
$$

De la m\^{e}me fa\c{c}on que l'on a une suite exacte
$$
1\rightarrow ({\Bbb G}_{{\rm m}}^{3}\times {\Bbb G}_{{\rm m}})/{\Bbb
G}_{{\rm m}}\rightarrow {\Bbb G}_{{\rm m}}^{S^{r,2}}\rightarrow
T^{r,2}\rightarrow 1
$$
on a une suite exacte
$$
1\rightarrow {\Bbb G}_{{\rm m}}\rightarrow {\Bbb G}_{{\rm
m}}^{S^{r,\tau}}\rightarrow T^{r,\tau}\rightarrow 1
$$
o\`{u} la deuxi\`{e}me fl\`{e}che envoie $u$ sur
$(u^{i_{0}+qi_{1}})_{i\in S^{r,\tau}}$. On v\'{e}rifie que le
plongement
$$
{\Bbb G}_{{\rm m}}^{S^{r,\tau}}\hookrightarrow {\Bbb G}_{{\rm
m}}^{S^{r,2}}
$$
qui envoie $(t_{i})_{i\in S^{r,\tau}}$ sur $(t_{i})_{i\in S^{r,2}}$,
o\`{u} $t_{(0,i_{1},i_{2})}=t_{(i_{1},0,i_{2})}^{q}$ si $i_{1}\not=0$
et $t_{(0,0,r)}=1$, induit un plongement $T^{r,\tau}\hookrightarrow
T^{r,2}$ qui se prolonge en une immersion ferm\'{e}e
$$
\overline{T}^{r,\tau}\hookrightarrow\overline{T}^{r,2}.
$$
On identifiera dans la suite $\overline{T}^{r,\tau}$ \`{a} l'image de
cette immersion ferm\'{e}e. 

On a vu que les trois applications strictement croissantes
$\{0,1\}\hookrightarrow\{0,1,2\}$, $(0,1)\mapsto (1,2),(0,2),(0,1)$,
induisent trois morphismes toriques $p_{0},p_{1},p_{2}:
\overline{T}^{r,2}\rightarrow\overline{T}^{r,1}$ et trois morphismes
de face $\pi_{0},\pi_{1},\pi_{2}:\Omega^{r,2}\rightarrow\Omega^{r,1}$
au-dessus de ceux-ci. Par construction, la relation
$$
p_{0}=\tau\circ p_{1}
$$
est v\'{e}rifi\'{e}e sur le ferm\'{e} $\overline{T}^{r,\tau}$ de
$\overline{T}^{r,2}$.

Soit alors
$$
\Omega^{r,\tau}\subset\Omega^{r,2}\times_{\overline{T}^{r,2}}
\overline{T}^{r,\tau}
$$
le sous-sch\'{e}ma ferm\'{e} d\'{e}fini par l'\'{e}quation
$\pi_{0}=\tau\circ\pi_{1}$ dans $\Omega^{r,1}$. L'action \'{e}vidente
du sous-tore ${\Bbb G}_{{\rm m}}^{S^{r,\tau}}\subset {\Bbb
G}^{S^{r,2}}$ sur $\Omega^{r,2} \times_{\overline{T}^{r,2}}
\overline{T}^{r,\tau}$ stabilise ce ferm\'{e}. On notera
$\pi_{1}^{\tau},\pi_{2}^{\tau}$ les restrictions des morphismes
$\pi_{1},\pi_{2}$ \`{a} $\Omega^{r,\tau}$. Elles rel\`{e}vent les
restrictions $p_{1}^{\tau},p_{2}^{\tau}: \overline{T}^{r,\tau}
\rightarrow\overline{T}^{r,1}$ de $p_{1},p_{2}$ et sont ${\Bbb
G}_{{\rm m}}^{S^{r,\tau}}$-\'{e}quivariantes.

\thm TH\'{E}OR\`{E}ME 
\enonce
La seconde projection $\Omega^{r,\tau}\rightarrow\overline{T}^{r,
\tau}$ est ${\Bbb G}_{{\rm m}}^{S^{r,\tau}}$-\'{e}quivariante et lisse
purement de dimension relative $r^{2}$.
\endthm

La stratification par les $T^{r,\tau}$-orbites de
$\overline{T}^{r,\tau}$ induit par image r\'{e}ciproque une
stratification de $\Omega^{r,\tau}$ dont les strates $\Omega_{{\bf
P}}^{r,\tau}$ sont index\'{e}es par les pavages entiers
$q$-admissibles. Un tel pavage ${\bf P}$ induit des pavages entiers
convexes admissibles ${\bf Q}_{0}={\bf Q}_{1}$ et ${\bf Q}_{2}$ sur
les c\^{o}t\'{e}s $\{x_{0}=0\}$, $\{x_{1}=0\}$ et $\{x_{2}=0\}$ de
$S_{{\Bbb R}}^{r,2}$ identifi\'{e}s \`{a} $S_{{\Bbb R}}^{r,1}$ et on a
des morphismes de face $p_{i}:\Omega_{{\bf P}}^{r,\tau}
\rightarrow\Omega_{{\bf Q}_{i}}^{r,1}$ pour $i=0,1,2$. Ces strates et
ces morphismes de face admettent des descriptions en termes des
sch\'{e}mas de graphes recoll\'{e}s $\mathop{\rm Gr}\nolimits_{{\bf
P}}^{r,2}$ et $\mathop{\rm Gr}\nolimits_{{\bf Q}}^{r,1}$.

En particulier, pour le pavage trivial, on v\'{e}rifie que la strate
ouverte $\Omega_{{\bf P}}^{r,\tau}$ est isomorphe au ferm\'{e}
$$
(\{(g_{0},g_{1},g_{2})\mid g_{1}g_{2}^{-1}=\tau (g_{0}g_{2}^{-1})\}
\times {\Bbb G}_{{\rm m}}^{S^{r, \tau}})/(\mathop{\rm GL}
\nolimits_{r}\times {\Bbb G}_{{\rm m}})\subset (\mathop{\rm GL}
\nolimits_{r}^{3}\times {\Bbb G}_{{\rm m}}^{S^{r, \tau}})/
(\mathop{\rm GL}\nolimits_{r}\times {\Bbb G}_{{\rm m}})
$$
o\`{u} le quotient est pris pour le plongement $(g,u)\mapsto
((u^{-1}g,u^{-q}g,g),(u^{i_{0}+qi_{1}})_{i\in S^{r,\tau}})$. On peut
identifier ce ferm\'{e} au sch\'{e}ma
$$
\mathop{\rm GL}\nolimits_{r}\times ({\Bbb G}_{{\rm m}}^{S^{r,
\tau}}/{\Bbb G}_{{\rm m}})
$$
par le morphisme $((g_{0},g_{1},g_{2}),t)\mapsto (g_{2}g_{0}^{-1},t)$.
Si l'on identifie $\Omega^{r,1}$ \`{a} $\mathop{\rm GL}
\nolimits_{r}\times ({\Bbb G}_{{\rm m}}^{S^{r,1}}/{\Bbb G}_{{\rm m}})$
par le morphisme $((g_{0},g_{1}),t)\mapsto (g_{1}g_{0}^{-1},t)$, les
projections $p_{1},p_{2}:\Omega_{{\bf P}}^{r,\tau}\rightarrow
\Omega_{{\bf Q}}^{r,1}$, o\`{u} ${\bf Q}={\bf Q}_{0}={\bf Q}_{1}={\bf
Q}_{2}$ est le pavage trivial de $S_{{\Bbb R}}^{r,1}$, sont induites
par les morphismes $\mathop{\rm GL}\nolimits_{r}\rightarrow\mathop{\rm
GL}\nolimits_{r}$ qui envoient $g$ sur $g$ et sur $\tau (g)^{-1}g$
respectivement et les morphismes ${\Bbb G}_{{\rm m}}^{S^{r,
\tau}}\rightarrow {\Bbb G}_{{\rm m}}^{S^{r, 1}}$ qui envoient $t$ sur
$(t_{(i_{0},0,i_{1})})_{i\in S^{r,1}}$ (avec $t_{(0,0,r)}=1$) et sur
$(t_{(i_{0},i_{1},0)})_{i\in S^{r,1}}$ (avec
$t_{(0,r,0)}=t_{r,0,0}^{q}$) respectivement.

Soit $S^{r,\tau ,\circ}=S^{r,\tau}-\{(1,0,r-1)\}$. L'inclusion ${\Bbb
G}_{{\rm m}}^{S^{r,\tau ,\circ}}\subset {\Bbb G}_{{\rm
m}}^{S^{r,\tau}}$ est une section de l'\'{e}pimorphisme ${\Bbb
G}_{{\rm m}}^{S^{r,\tau}}\twoheadrightarrow T^{r,\tau}$.

\thm TH\'{E}OR\`{E}ME 
\enonce
Le quotient $\mathop{\widetilde{\rm gl}}_{r}^{\tau}$ de
$\Omega^{r,\tau}$ par l'action libre du tore ${\Bbb G}_{{\rm
m}}^{S^{r,\tau ,\circ}}$ contient $\mathop{\rm GL}\nolimits_{r}$ comme
un ouvert dense. Il est lisse sur le champ torique ${\cal
T}^{r,\tau}=[\overline{T}^{r,\tau}/T^{r,\tau}]$ et il est muni de deux
morphismes projectifs sur le sch\'{e}ma des homomorphismes complets
$\mathop{\widetilde{\rm gl}}_{r}$ qui s'ins\`{e}rent dans les deux
carr\'{e}s commutatifs
$$\diagram{
\mathop{\rm GL}\nolimits_{r}&\,\smash{\mathop{
\lhook\joinrel\mathrel{\hbox to 8mm{\rightarrowfill}}
}\limits^{\scriptstyle }}\,&\mathop{\widetilde{\rm gl}}_{r}^{\tau}\cr
\llap{$\scriptstyle \mathop{\rm Id}$}\left\downarrow
\vbox to 4mm{}\right.\rlap{}&+&\llap{}\left\downarrow
\vbox to 4mm{}\right.\rlap{$\scriptstyle$}\cr
\mathop{\rm GL}\nolimits_{r}&\,\smash{\mathop{
\lhook\joinrel\mathrel{\hbox to 8mm{\rightarrowfill}}
}\limits^{\scriptstyle }}\,&\mathop{\widetilde{\rm gl}}_{r}\cr}
\kern 5mm\hbox{ et }\kern 5mm
\diagram{
\mathop{\rm GL}\nolimits_{r}&\,\smash{\mathop{
\lhook\joinrel\mathrel{\hbox to 8mm{\rightarrowfill}}
}\limits^{\scriptstyle }}\,&\mathop{\widetilde{\rm gl}}_{r}^{\tau}\cr
\llap{$\scriptstyle L$}\left\downarrow
\vbox to 4mm{}\right.\rlap{}&+&\llap{}\left\downarrow
\vbox to 4mm{}\right.\rlap{$\scriptstyle$}\cr
\mathop{\rm GL}\nolimits_{r}&\,\smash{\mathop{
\lhook\joinrel\mathrel{\hbox to 8mm{\rightarrowfill}}
}\limits^{\scriptstyle }}\,&\mathop{\widetilde{\rm gl}}_{r}\cr}
$$
o\`{u} $L:g\rightarrow \tau (g)^{-1}g$ est l'isog\'{e}nie de Lang.
\endthm

\rem Remarque 
\endrem
En passant au quotient par l'action libre du tore ${\Bbb G}_{{\rm
m}}^{r,\tau}$ tout entier, on obtient un sch\'{e}ma {\it projectif}
$\overline{\Omega}^{r,\tau}$, muni d'un morphisme vers la
compactification de De Concini et Procesi de $\mathop{\rm
PGL}\nolimits_{r}$ qui prolonge l'isog\'{e}nie de Lang $g\mapsto \tau
(g)^{-1}g$ de ce sch\'{e}ma en groupes.
\hfill\hfill$\square$

\vskip 7mm

\centerline{\bf 5. ALLURE DES NOMBRES DE LEFSCHETZ}
\vskip 10mm

{\bf 5.1. Correspondances de Hecke et endomorphismes de Frobenius
partiels}
\vskip 5mm

On fixe dans toute cette section un niveau $N=\mathop{\rm Spec}({\cal
O}_{N})\hookrightarrow X$.

Tout \'{e}l\'{e}ment $g\in \mathop{\rm GL}\nolimits_{r} ({\Bbb A})$
d\'{e}finit une correspondance dite {\it de Hecke}
$$
c=(c_{1},c_{2}):\mathop{\rm Cht}\nolimits_{N}^{r}(g)\rightarrow
\mathop{\rm Cht}\nolimits_{N}^{r}\times_{(X-N)^{2}}\mathop{\rm
Cht}\nolimits_{N}^{r}
$$
o\`{u} $\mathop{\rm Cht}\nolimits_{N}^{r}(g)$ est un champ de
Deligne-Mumford et $c_{1},c_{2}$ sont des morphismes
re\-pr\'{e}\-sen\-ta\-bles \'{e}tales. Cette correspondance ne
d\'{e}pend que de la double classe $K_{N}gK_{N}$ o\`{u} $K_{N}$ est le
sous-groupe d'indice fini $\mathop{\rm Ker} (\mathop{\rm GL}
\nolimits_{r}({\cal O}) \twoheadrightarrow\mathop{\rm GL}
\nolimits_{r}({\cal O}_{N}))$ de $K=\mathop{\rm GL}\nolimits_{r}
({\cal O})$. Si $S_{g}$ est la r\'{e}union du support $|N|\subset |X|$
de $N$ et de l'ensemble fini des places $x$ de $X$ telles que
$g_{x}\notin F_{x}^{\times}K_{N,x} \subset \mathop{\rm
GL}\nolimits_{r}(F_{x})$, $c_{1}$ et $c_{2}$ sont finies au-dessus de
l'ouvert $(X-S_{g})^{2}$ de $(X-N)^{2}$.

Si $a\in {\Bbb A}^{\times}$, la double classe $F^{\times}a\mathop{\rm
Ker}({\cal O}^{\times}\twoheadrightarrow {\cal O}_{N}^{\times})$ dans
$F^{\times}\backslash {\Bbb A}^{\times}/\mathop{\rm Ker}({\cal
O}^{\times}\twoheadrightarrow {\cal O}_{N}^{\times})$ est la classe
d'isomorphie d'un ${\cal O}_{X}$-Module inversible ${\cal L}$ muni
d'une trivialisation de sa restriction \`{a} $N$. Si $g\in
aK_{N}\subset \mathop{\rm GL}\nolimits_{r}({\Bbb A})$, on a
$\mathop{\rm Cht}\nolimits_{N}^{r}(g)= \mathop{\rm
Cht}\nolimits_{N}^{r}$, $c_{1}$ est l'identit\'{e} et $c_{2}$ est
l'automorphisme de $\mathop{\rm Cht}\nolimits_{N}^{r}$ qui envoie un
chtouca $\widetilde{{\cal E}}$ muni d'une structure de niveau $N$ sur
le produit tensoriel ${\cal L}\otimes \widetilde{{\cal E}}$ muni de
la structure de niveau $N$ \'{e}vidente. On notera simplement $a$ cet
automorphisme.

Soit $\Lambda_{X}$ le sch\'{e}ma intersection de tous les ouverts de
$X\times X$ compl\'{e}mentaires des images r\'{e}ciproques de la
diagonales par les endomorphismes $\mathop{\rm Frob}\nolimits_{X}^{n}
\times\mathop{\rm Id}\nolimits_{X}$ et $\mathop{\rm Id}\nolimits_{X}
\times\mathop{\rm Frob}\nolimits_{X}^{n}$ de $X\times X$ pour tous les
entiers $n\geq 0$. Au-dessus de $\Lambda_{X}$, il n'y a plus de
diff\'{e}rence entre chtoucas \`{a} gauche et chtoucas \`{a} droite:
un chtouca \`{a} droite ${\cal E}\,\smash{\mathop{
\lhook\joinrel\mathrel{\hbox to 6mm{\rightarrowfill}}
}\limits^{\scriptstyle j}}\,{\cal E}'\,\smash{\mathop{ {\hbox to
6mm{\leftarrowfill}}\joinrel\kern -0.9mm\mathrel\rhook}
\limits^{\scriptstyle t}}\,{}^{\tau}{\cal E}$ dont le couple
p\^{o}le-z\'{e}ro $(\infty ,o)$ se factorise par $\Lambda_{X}\subset
X\times X$ d\'{e}finit un chtouca \`{a} gauche ${\cal
E}\,\smash{\mathop{ {\hbox to 6mm{\leftarrowfill}}\joinrel\kern
-0.9mm\mathrel\rhook} \limits^{\scriptstyle t'}}\,{\cal
E}''\,\smash{\mathop{ \lhook\joinrel\mathrel{\hbox to
6mm{\rightarrowfill}} }\limits^{\scriptstyle j'}}\,{}^{\tau}{\cal E}$
de p\^{o}le $\infty'=\infty$ et de z\'{e}ro $o'=o$, o\`{u} ${\cal
E}''={\cal E}\times_{j,{\cal E}',t}{}^{\tau}{\cal E}$ et $t',j'$ sont
les deux projections, et vice versa. Au-dessus de $\Lambda_{X}$, les
morphismes de Frobenius partiels introduits en 2.1 peuvent donc
\^{e}tre vus comme des endomorphismes, dits encore de Frobenius
partiels,
$$
\mathop{\rm Frob}\nolimits_{\infty}\hbox{ et }\mathop{\rm Frob}
\nolimits_{o}:\Lambda_{X}\times_{X^{2}}\mathop{\rm Cht}
\nolimits_{N}^{r}\rightarrow \Lambda_{X}\times_{X^{2}}\mathop{\rm Cht}
\nolimits_{N}^{r}
$$
au dessus de $\mathop{\rm Frob}\nolimits_{X}\times\mathop{\rm Id}
\nolimits_{X}$ et $\mathop{\rm Id}\nolimits_{X}\times \mathop{\rm
Frob}\nolimits_{X}$ respectivement.

On fixe dans la suite un \'{e}l\'{e}ment $a\in {\Bbb A}^{\times}$ de
degr\'{e} $\mathop{\rm deg}(a)>0$ et on consid\`{e}re le quotient
$$
\mathop{\rm Cht}\nolimits_{N}^{r}/a^{{\Bbb Z}}\cong
\coprod_{d=1}^{r\mathop{\rm deg}(a)}\mathop{\rm Cht}
\nolimits_{N}^{r,d}.
$$
Les correspondances de Hecke et les endomorphismes de Frobenius
partiels que nous venons d'introduire passent au quotient par
$a^{{\Bbb Z}}$.

Pour tout param\`{e}tre de troncature $p:[0,r]\rightarrow {\Bbb R}$ et
tout $\alpha\in [0,1]$, $a^{{\Bbb Z}}$ stabilise l'ouvert $\mathop{\rm
Cht}\nolimits_{N}^{r; \,\leq_{\alpha}p}$ de $\mathop{\rm Cht}
\nolimits_{N}^{r}$ et l'ouvert quotient
$$
\mathop{\rm Cht}\nolimits_{N}^{r;\,\leq_{\alpha}p}/a^{{\Bbb Z}}\subset
\mathop{\rm Cht}\nolimits_{N}^{r}/a^{{\Bbb Z}}
$$
est de type fini. {\it Par contre, cet ouvert n'est stabilis\'{e} ni
par les endomorphismes de Frobenius partiels $\mathop{\rm
Frob}\nolimits_{\infty}$, $\mathop{\rm Frob} \nolimits_{o}$, ni par
les correspondances de Hecke} ({\it sauf celles qui sont associ\'{e}es
\`{a} des $g\in {\Bbb A}^{\times}K$}). C'est la difficult\'{e} majeure
qu'a d\^{u} surmonter Lafforgue.
\vskip 5mm

{\bf 5.2. Nombres de Lefschetz} 
\vskip 5mm

Si $\infty$ et $o$ sont deux points ferm\'{e}s de $X$ le
sous-sch\'{e}ma fini $\infty\times o\subset X\times X$ a exactement
$\delta (\infty ,o)$ points ferm\'{e}s, o\`{u} $\delta (\infty ,o)$
est le plus grand commun diviseur de $\mathop{\rm deg}(\infty )$ et
$\mathop{\rm deg}(o)$. Le corps r\'{e}siduel $\kappa (\xi )$ de chaque
$\xi\in\infty\times o$ est une extension compos\'{e}e de $\kappa
(\infty )$ et $\kappa (o)$, et son degr\'{e} $\mathop{\rm deg}(\xi )$
sur ${\Bbb F}_{q}$ est le plus petit commun multiple
$$
\mu (\infty ,o)={\mathop{\rm deg}(\infty )\mathop{\rm
deg}(o)\over\delta (\infty ,o)}
$$
de $\mathop{\rm deg}(\infty )$ et $\mathop{\rm deg}(o)$. Pour chaque
$\xi\in\infty\times o$, les suites de points ferm\'{e}s de
$X\times X$,
$$
n\mapsto (\mathop{\rm Frob}\nolimits_{X}^{n}\times\mathop{\rm
Id}\nolimits_{X})(\xi)\hbox{ et }n\mapsto (\mathop{\rm
Id}\nolimits_{X}\times\mathop{\rm Frob}\nolimits_{X}^{n})(\xi)
$$
sont en fait \`{a} valeurs dans $\infty\times o$ et sont
p\'{e}riodiques de p\'{e}riodes $\mathop{\rm deg}(\infty )$ et
$\mathop{\rm deg}(o)$ respectivement.
\vskip 3mm

Fixons $g\in\mathop{\rm GL}\nolimits_{r}({\Bbb A})$, deux points
ferm\'{e}s distincts $\infty$ et $o$ de $X-S_{g}$, un entier
$t$, deux entiers $n_{\infty},n_{o}\geq 1$ tels que $\mathop{\rm
deg}(o)n_{o}-\mathop{\rm deg}(\infty )n_{\infty}=t$, un point
ferm\'{e} $\xi$ de $\infty\times o$, un param\`{e}tre de troncature
$p:[0,r]\rightarrow {\Bbb R}$ et un nombre r\'{e}el $\alpha\in [0,1]$.

On note $\mathop{\rm Cht}\nolimits_{N,\xi}^{r}/a^{{\Bbb Z}}$ la fibre
en $\xi$ de la projection canonique $\mathop{\rm Cht}\nolimits_{N}^{r}
/a^{{\Bbb Z}}\rightarrow (X-N)^{2}$. C'est un champ alg\'{e}brique de
Deligne-Mumford lisse purement de dimension relative $2r-2$ sur le
corps fini $\kappa (\xi )$. On note
$$
c_{\xi}=(c_{1,\xi},c_{2,\xi}):\mathop{\rm Cht}\nolimits_{N,\xi}^{r}(g)
/a^{{\Bbb Z}}\rightarrow\mathop{\rm Cht}\nolimits_{N,\xi}^{r}/
a^{{\Bbb Z}} \times_{\kappa (\xi )}\mathop{\rm Cht}
\nolimits_{N,\xi}^{r}/a^{{\Bbb Z}}
$$
la fibre en $\xi$ de la correspondance de Hecke d\'{e}finie par $g$.
Les deux composantes $c_{1,\xi}$ et $c_{2,\xi}$ sont
repr\'{e}sentables finies \'{e}tales.

Puisque $\infty$ et $o$ sont distincts, le point ferm\'{e} $\xi$ est
dans $\Lambda_{X}\subset X^{2}$. Les puissances $\mathop{\rm
Frob}\nolimits_{\infty}^{\mathop{\rm deg}(\infty )}$ et $\mathop{\rm
Frob}\nolimits_{o}^{\mathop{\rm deg}(o)}$ des endomorphismes de
Frobenius partiels sont donc d\'{e}finies sur la fibre $\mathop{\rm
Cht} \nolimits_{N,\xi}^{r}/a^{{\Bbb Z}}$ et elles la stabilisent.

\thm D\'{E}FINITION 
\enonce
Le nombre de Lefschetz
$$
\mathop{\rm Lef}\nolimits_{\xi}(g\times\mathop{\rm Frob}
\nolimits_{\infty}^{\mathop{\rm deg}(\infty )n_{\infty}}
\times\mathop{\rm Frob} \nolimits_{o}^{\mathop{\rm deg}(o)n_{o}},
\mathop{\rm Cht}\nolimits_{N}^{r; \,\leq_{\alpha} p}/a^{{\Bbb Z}})
$$
est la somme
$$
\sum_{y}{1\over |\mathop{\rm Aut}(y)|}
$$
o\`{u} $y$ parcourt l'ensemble des classes d'isomorphie de points dans
$\mathop{\rm Cht}\nolimits_{N,\xi}^{r}(g)/a^{{\Bbb Z}}$ munis d'un 
isomorphisme
$$\eqalign{
c_{1,\xi}(y)&\cong (\mathop{\rm Frob}\nolimits_{\infty}^{\mathop{\rm
deg}(\infty )n_{\infty}} \circ\mathop{\rm Frob}
\nolimits_{o}^{\mathop{\rm deg}(o)n_{o}})(c_{2,\xi}(y))\cr
&\cong (\mathop{\rm Frob}\nolimits_{o}^{\mathop{\rm deg}(o)n_{o}}\circ
\mathop{\rm Frob}\nolimits_{\infty}^{\mathop{\rm deg}(\infty )
n_{\infty}})(c_{2,\xi}(y))\cr}
$$
et tels que
$$
c_{1,\xi}(y)\in\mathop{\rm Cht}\nolimits_{N,\xi}^{r;\,\leq_{\alpha} p}
/a^{{\Bbb Z}}\subset\mathop{\rm Cht}\nolimits_{N,\xi}^{r}/a^{{\Bbb
Z}},
$$
et o\`{u} $\mathop{\rm Aut}(y)$ est le groupe fini des automorphismes 
du point fixe $y$.
\endthm

En utilisant
\vskip 2mm

\itemitem{-} la description ad\'{e}lique de Drinfeld des points de
$\mathop{\rm Cht}\nolimits_{N,\xi}^{r}$ ([Dr 5], [La 1]),
\vskip 2mm

\itemitem{-} le cas particulier du {\og}{Lemme Fondamental}{\fg}
prouv\'{e} par Drinfeld ([Lau 1]),
\vskip 2mm

\itemitem{-} la d\'{e}composition spectrale de Langlands ([Lan 2],
[Mo-Wa]) et la formule des traces d'Arthur-Selberg ([Ar], [La 1]),
\vskip 2mm

\itemitem{-} la formule de Cauchy pour calculer les int\'{e}grales 
qui apparaissent dans la formule des traces,
\vskip 2mm

\noindent Lafforgue d\'{e}montre alors:

\thm TH\'{E}OR\`{E}ME
\enonce
Si $\mathop{\rm deg}(\infty )$ et $\mathop{\rm deg}(o)$ sont assez
grands en fonction de la double classe $K_{N}gK_{N}$ et de l'entier
$t$ et si $p$ est assez convexe en fonction du niveau $N$ et de la
double classe $K_{N}gK_{N}$, la moyenne de chacune des deux suites
p\'{e}riodiques
$$
k\mapsto \mathop{\rm Lef}\nolimits_{(\mathop{\rm Frob}
\nolimits_{X}^{k}\times \mathop{\rm Id}\nolimits_{X})(\xi)}
(g\times\mathop{\rm Frob}\nolimits_{\infty}^{\mathop{\rm deg}(\infty
)n_{\infty}} \times\mathop{\rm Frob} \nolimits_{o}^{\mathop{\rm
deg}(o)n_{o}},\mathop{\rm Cht}\nolimits_{N}^{r; \,\leq_{\alpha}
p}/a^{{\Bbb Z}})
$$
et
$$
k\mapsto\mathop{\rm Lef}\nolimits_{(\mathop{\rm Id}
\nolimits_{X}\times \mathop{\rm Frob}\nolimits_{X}^{k})(\xi)}
(g\times\mathop{\rm Frob} \nolimits_{\infty}^{\mathop{\rm deg}(\infty
)n_{\infty}} \times\mathop{\rm Frob} \nolimits_{o}^{\mathop{\rm
deg}(o)n_{o}},\mathop{\rm Cht}\nolimits_{N}^{r; \,\leq_{\alpha}
p}/a^{{\Bbb Z}}),
$$
de p\'{e}riode le p.g.c.d. de $\mathop{\rm deg}(o)$, $\mathop{\rm
deg}(\infty )$ et $r!$, est \'{e}gale \`{a} l'expression spectrale
suivante
$$\eqalign{
&\sum_{{\scriptstyle\pi\in {\cal A}_{r}(K_{N})\atop\scriptstyle
\omega_{\pi}(a)=1}}\mathop{\rm Tr}\nolimits_{\pi}(1_{K_{N}gK_{N}})
q^{(\mathop{\rm deg}(\infty )n_{\infty}+\mathop{\rm
deg}(o)n_{o}){r-1\over 2}}S_{\infty}^{(-n_{\infty})}\!(\pi )
S_{o}^{(n_{o})}\!(\pi )\cr
&+\kern -2mm\sum_{{\scriptstyle 1\leq r'<r\atop\scriptstyle 1\leq
r''<r}}\sum_{{\scriptstyle\pi'\in {\cal A}_{r'}(K_{N})\atop
\scriptstyle\pi''\in {\cal A}_{r''}(K_{N})}}\mathop{\rm Tr}
\nolimits_{\pi',\pi''}^{\leq_{\alpha}p}(1_{K_{N}gK_{N}},\mathop{\rm
deg}(\infty )n_{\infty}) S_{\infty}^{(-n_{\infty})}\!
(\pi')S_{o}^{(n_{o})}\!(\pi'')\cr}
$$
o\`{u}
\vskip 2mm

\itemitem{-} ${\cal A}_{r}(K_{N})\subset {\cal A}_{r}$ est un syst\`{e}me
de repr\'{e}sentants des classes d'isomorphie de repr\'{e}sentations
automorphes cuspidales irr\'{e}ductibles de $\mathop{\rm
GL}\nolimits_{r}({\Bbb A})$ qui admettent au moins un vecteur non nul
fixe sous $K_{N}$,
\vskip 2mm

\itemitem{-} $1_{K_{N}gK_{N}}$ est la fonction caract\'{e}ristique de
$K_{N}gK_{N}\subset \mathop{\rm GL}\nolimits_{r}({\Bbb A})$ et ${\rm
Tr}_{\pi}(1_{K_{N}gK_{N}})$ est la trace de l'op\'{e}rateur $\pi
(1_{K_{N}gK_{N}})$,
\vskip 2mm

\itemitem{-} pour $x=\infty ,o$ et $\nu$ un entier positif,
$S_{x}^{(\nu )}(\pi )=z_{1}(\pi_{x})^{\nu}+ \cdots
+z_{r}(\pi_{x})^{\nu}$ est la somme des $\nu$-i\`{e}mes puissances des
valeurs propres de Hecke de la composante locale non ramifi\'{e}e
$\pi_{x}$ de $\pi$,
\vskip 2mm

\itemitem{-} $\nu\mapsto\mathop{\rm Tr}\nolimits_{\pi',
\pi''}^{\leq_{\alpha}p}(1_{K_{N}gK_{N}},\nu )$ est une fonction
complexe de l'entier $\nu$, qui ne d\'{e}pend pas des places
$o,\infty\in |X|-S_{g}$, qui est de la forme
$$
\sum_{z\in {\Bbb C}^{\times}}P_{z}(\nu )z^{\nu}
$$
pour une famille \`{a} support fini de polyn\^{o}mes $P_{z}(u)\in
{\Bbb C}[u]$, $z\in {\Bbb C}^{\times}$, et qui est identiquement nulle
sauf pour un nombre fini de couples $(\pi',\pi'')$.

\hfill\hfill$\square$
\endthm

\vskip 7mm

\centerline{\bf 6. REPR\'{E}SENTATIONS GALOISIENNES}
\vskip 10mm

{\bf 6.1. Repr\'{e}sentations de Galois $r$-n\'{e}gligeables} 
\vskip 5mm

Rappelons que $F$ est le corps des fonctions rationnelles sur la
courbe $X$ et que l'on a fix\'{e} une cl\^{o}ture s\'{e}parable
$\overline{F}$ de $F$. On d\'{e}signera par $\overline{{\Bbb F}}_{q}$
la cl\^{o}ture alg\'{e}brique de ${\Bbb F}_{q}$ dans $\overline{F}$.
On identifiera le groupe de Galois $\Gamma_{{\Bbb F}_{q}}:=\mathop{\rm
Gal}(\overline{{\Bbb F}}_{q}/{\Bbb F}_{q})$ \`{a} $\widehat{{\Bbb Z}}$
en choisissant pour g\'{e}n\'{e}rateur topologique de $\Gamma_{{\Bbb
F}_{q}}$ l'\'{e}l\'{e}ment de Frobenius g\'{e}om\'{e}trique (l'inverse
de l'\'{e}l\'{e}vation \`{a} la puissance $q$-i\`{e}me dans
$\overline{{\Bbb F}}_{q}$).

Notons $E$ le corps des fractions de $F\otimes F$, c'est-\`{a}-dire le
corps des fonctions rationnelles sur la surface $X\times X$, et fixons
une cl\^{o}ture alg\'{e}brique $\overline{E}$ de $E$ et un plongement
$\overline{F}\otimes_{\overline{{\Bbb F}}_{q}}\overline{F}
\hookrightarrow \overline{E}$ au dessus du plongement $F\otimes
F\hookrightarrow E$. On a donc d\'{e}fini un point g\'{e}om\'{e}trique
$\overline{\delta}:\mathop{\rm Spec}(\overline{E}) \rightarrow X\times
X$ localis\'{e} en le point g\'{e}n\'{e}rique $\delta =\mathop{\rm
Spec}(E)$ de $X\times X$ et tel que ses images par les deux
projections canoniques de $X\times X$ se factorisent par le point
g\'{e}om\'{e}trique $\overline{\eta}:\mathop{\rm Spec}(\overline{F})
\rightarrow X$ localis\'{e} au point g\'{e}n\'{e}rique $\eta$ de $X$.

Pour tout ouvert non vide $U$ de $X$, le groupe fondamental de
Grothendieck $\pi_{1}(U,\overline{\eta})$ est le quotient de
$\Gamma_{F}=\mathop{\rm Gal}(\overline{F}/F)$ qui classifie les
extensions finies de $F$ dans $\overline{F}$ qui sont non
ramifi\'{e}es en toutes les places $x\in |U|$. Il admet comme quotient
le groupe de Galois $\Gamma_{{\Bbb F}_{q}}=\widetilde{{\Bbb Z}}$. De
m\^{e}me, pour tout ouvert non vide $V$ de $X\times X$, le groupe
fondamental de Grothendieck $\pi_{1}(V,\overline{\delta})$ est un
quotient du groupe de Galois $\mathop{\rm Gal}(\overline{E}/E)$ et
admet $\widetilde{{\Bbb Z}}$ comme quotient.

Soit $V$ un ouvert non vide de $X\times X$ de la forme $U_{1}\times
U_{2}$ pour des ouverts $U_{1}$ et $U_{2}$ de $X$. Cet ouvert est
stable par les endomorphismes $\mathop{\rm
Frob}\nolimits_{X}\times\mathop{\rm Id}\nolimits_{X}$ et $\mathop{\rm
Id}\nolimits_{X}\times\mathop{\rm Frob} \nolimits_{X}$ de $X\times X$
et les foncteurs {\og}{images r\'{e}ciproques}{\fg} par ces
endomorphismes sont des \'{e}quivalences, quasi-inverses l'une de
l'autre, de la cat\'{e}gorie des rev\^{e}tements finis \'{e}tales de
$V$ sur elle-m\^{e}me (le foncteur image r\'{e}ciproque par
$\mathop{\rm Frob}\nolimits_{V}=\mathop{\rm Frob}
\nolimits_{U_{1}}\times \mathop{\rm Frob}\nolimits_{U_{2}}$ est
canoniquement isomorphe \`{a} l'identit\'{e}). On peut donc
consid\'{e}rer la cat\'{e}gorie des rev\^{e}tements finis \'{e}tales
$V'$ de $V$ munis d'un isomorphisme $(\mathop{\rm
Frob}\nolimits_{U_{1}}\times\mathop{\rm Id}\nolimits_{U_{2}})^{\ast}
V'\buildrel\sim\over\longrightarrow V'$ ou, ce qui revient au
m\^{e}me, d'un isomorphisme $(\mathop{\rm Id}
\nolimits_{U_{1}}\times\mathop{\rm Frob}\nolimits_{U_{2}})^{\ast}
V'\buildrel\sim\over\longrightarrow V'$. Le foncteur fibre en
$\overline{\delta}$ identifie cette cat\'{e}gorie \`{a} la
cat\'{e}gorie des ensembles finis munis de l'action d'un groupe
pro-fini $\widetilde{\pi}_{1}(U_{1}\times U_{2},\overline{\delta})$
qui s'ins\`{e}re dans un diagramme commutatif
$$
\matrix{1 & \longrightarrow & \pi_{1}(U_{1}\times U_{2},
\overline{\delta}) & \longrightarrow & \widetilde{\pi}_{1}(U_{1}\times
U_{2},\overline{\delta}) & \longrightarrow & \widehat{{\Bbb Z}} &
\longrightarrow & 0\cr \noalign{\smallskip}
& & \llap{$\scriptstyle $}\left\downarrow \vbox to 3mm{}\right.\rlap{}
& & \llap{$\scriptstyle $}\left\downarrow \vbox to 3mm{}\right.\rlap{}
& & \big|\!\big| &&\cr
\noalign{\smallskip}
0 & \smash{\mathop{\hbox to 15mm{\rightarrowfill}}
\limits_{\scriptstyle }} \kern -9mm & \widehat{{\Bbb Z}} & \kern
-9mm\smash{\mathop{\hbox to 23mm{\rightarrowfill}}
\limits_{\scriptstyle }} \kern -8mm & \widehat{{\Bbb Z}}^{2} & \kern
-8mm\smash{\mathop{\hbox to 14mm{\rightarrowfill}}
\limits_{\scriptstyle }} & \widehat{{\Bbb Z}} & \rightarrow & 0\cr}
$$
\`{a} lignes exactes et \`{a} fl\`{e}ches verticales surjectives. Tout
rev\^{e}tement fini \'{e}tale $V'$ de $V=U_{1}\times U_{2}$ de la
forme $V'=U_{1}'\times U_{2}'$ o\`{u} chaque $U_{i}'$ est un
rev\^{e}tement fini \'{e}tale de $U_{i}$, est canoniquement muni d'un
isomorphisme $(\mathop{\rm Frob}\nolimits_{U_{1}}\times \mathop{\rm
Id} \nolimits_{U_{2}})^{\ast}V'\buildrel\sim\over \longrightarrow V'$.
On a donc un homomorphisme de $\widetilde{\pi}_{1}(U_{1}\times
U_{2},\overline{\delta})$ dans $\pi_{1}(U_{1},\overline{\eta})\times
\pi_{1}(U_{2},\overline{\eta}) $ qui prolonge l'homomorphisme
$\pi_{1}(U_{1}\times U_{2}, \overline{\delta})\rightarrow
\pi_{1}(U_{1},\overline{\eta})\times \pi_{1}(U_{2}, \overline{\eta})$
d\'{e}fini par les deux projections canoniques de $U_{1}\times U_{2}$,
et qui induit l'identit\'{e} entre les quotients $\widehat{{\Bbb
Z}}^{2}$ de sa source et de son but.

\thm PROPOSITION (Drinfeld)
\enonce
L'homomorphisme d\'{e}fini ci-dessus
$$
\widetilde{\pi}_{1}(U_{1}\times U_{2},\overline{\delta})\rightarrow
\pi_{1}(U_{1},\overline{\eta})\times\pi_{1}(U_{2},\overline{\eta}),
$$
est un isomorphisme.
\endthm

On a fix\'{e} dans la section 1 un nombre premier $\ell$ distinct de
la caract\'{e}ristique de ${\Bbb F}_{q}$ et une cl\^{o}ture
alg\'{e}brique $\overline{{\Bbb Q}}_{\ell}$ de ${\Bbb Q}_{\ell}$.

Rappelons que, pour tout ouvert non vide $U$ de $X$ (resp. $V$ de
$X\times X$), le foncteur fibre en $\overline{\eta}$ (resp.
$\overline{\delta}$) est une \'{e}quivalence de la cat\'{e}gorie des
syst\`{e}mes locaux $\ell$-adiques (c'est-\`{a}-dire des
$\overline{{\Bbb Q}}_{\ell}$-faisceaux lisses) sur $U$ (resp. $V$) sur
la cat\'{e}gorie des repr\'{e}sentations $\ell$-adiques de
$\pi_{1}(U,\overline{\eta})$ (resp. $\pi_{1}(V,\overline{\delta})$).
Pour tout syst\`{e}me local $\ell$-adique $L$ sur $U$ (resp. $M$ sur
$V$) et tout point ferm\'{e} $x$ de $U$ (resp. $y$ de $V$), la classe
de conjugaison dans $\pi_{1}(U,\overline{\eta})$ (resp.
$\pi_{1}(V,\overline{\delta})$) des \'{e}l\'{e}ments de Frobenius
g\'{e}om\'{e}triques en $x$ (resp. $y$) d\'{e}finit donc une classe de
conjugaison d'endomorphismes de $L_{\overline{\eta}}$ (resp.
$M_{\overline{\delta}}$). On note $\mathop{\rm Tr}(\mathop{\rm
Frob}\nolimits_{x},L)$ (resp. $\mathop{\rm Tr}(\mathop{\rm
Frob}\nolimits_{y},M)$) sa trace et $\mathop{\rm det}(1-T\mathop{\rm
Frob}\nolimits_{x},L)$ (resp. $\mathop{\rm det}(1-T\mathop{\rm
Frob}\nolimits_{y},M)$) son polyn\^{o}me caract\'{e}ristique.

Soient $U$ un ouvert non vide de $X$ et $M$ un syst\`{e}me local
$\ell$-adique {\it irr\'{e}ductible} sur $U\times U$. On consid\`{e}re
les propri\'{e}t\'{e}s suivantes:
\vskip 1mm

\itemitem{(i)} $M$ est de la forme $L_{1}\boxtimes L_{2}$ pour des
syst\`{e}mes locaux $\ell$-adiques (irr\'{e}ductibles) $L_{1}$ sur
$L_{2}$ sur $U$ et, en particulier, il existe des fonctions $\Delta_{1},
\Delta_{2}: |U\times U|\rightarrow 1+T\overline{{\Bbb Q}}_{\ell}[T]$
telles que
$$
\mathop{\rm det}(1-T\mathop{\rm Frob}\nolimits_{\xi},M)
=\Delta_{1}(\infty )\ast\Delta_{2}(o)
$$
pour tout $\xi\in |U\times U|$ de projections $\infty ,o\in 
|U|$, o\`{u} l'op\'{e}ration $\ast$ est d\'{e}finie par
$$
\prod_{i_{1}=1}^{r_{1}}(1-\alpha_{1,i_{1}}T)\ast
\prod_{i_{2}=1}^{r_{2}}(1-\alpha_{2,i_{2}}T)=
\prod_{i_{1},i_{2}}(1-\alpha_{1,i_{1}}\alpha_{2,i_{2}}T).
$$
\vskip 1mm

\itemitem{(ii)} Il existe un isomorphisme $(\mathop{\rm Frob}
\nolimits_{U}\times\mathop{\rm Id}\nolimits_{U})^{\ast}M
\buildrel\sim\over\longrightarrow M$ et, en particulier, la fonction
$\xi\mapsto \mathop{\rm det}(1-T\mathop{\rm Frob}\nolimits_{\xi},M)$
est constante sur $|\infty\times o|$ pour chaque couple $(\infty ,o)$ de
points ferm\'{e}s de $U$.
\vskip 1mm

La propri\'{e}t\'{e} (i) implique \'{e}videmment la propri\'{e}t\'{e}
(ii) et il devrait r\'{e}sulter de la proposition de Drinfeld que les
propri\'{e}t\'{e}s (i) et (ii) sont en fait \'{e}quivalentes.

Lafforgue d\'{e}montre un r\'{e}sultat partiel dans cette direction.

\thm PROPOSITION 
\enonce
On suppose de plus que $M$ est ponctuellement pur au sens de {\rm
[De]} et qu'il existe des fonctions $\Delta_{1},\Delta_{2}:
|U\times U|\rightarrow 1+T\overline{{\Bbb Q}}_{\ell}[T]$ telles que
$$
\mathop{\rm det}(1-T\mathop{\rm Frob}\nolimits_{\xi},M)
=\Delta_{1}(\infty )\ast\Delta_{2}(o)
$$
pour tout $\xi\in |U\times U|$ dont les projections $\infty ,o\in |U|$
sont {\rm distinctes}.

Alors, la propri\'{e}t\'{e} {\rm (ii)} implique la propri\'{e}t\'{e}
{\rm (i)}.
\endthm

\thm D\'{E}FINITION 
\enonce
Si $U$ est un ouvert de $X$, un syst\`{e}me local $\ell$-adique sur
$U\times U$ est dit {\rm $r$-n\'{e}gligeable} si tous ses
sous-quotients irr\'{e}ductibles sont de la forme $L_{1}\boxtimes
L_{2}$ o\`{u} $L_{1}$ et $L_{2}$ sont des syst\`{e}mes locaux
$\ell$-adiques {\rm (}irr\'{e}ductibles{\rm )} {\rm de rangs} $\leq
r-1$ sur $U$.
\endthm

Si $\lambda$ est une unit\'{e} $\ell$-adique, on notera par
$\overline{{\Bbb Q}}_{\ell}^{(\lambda )}$ le $\overline{{\Bbb
Q}}_{\ell}$-faisceau sur $\mathop{\rm Spec}({\Bbb F}_{q})$ d\'{e}fini
par le caract\`{e}re $\overline{{\Bbb Z}}\rightarrow\overline{{\Bbb
Q}}_{\ell}^{\times},~1\mapsto\lambda$, du groupe de Galois
$\Gamma_{{\Bbb F}_{q}}$. Si $F$ est un faisceau $\ell$-adique sur un
(${\Bbb F}_{q}$-)sch\'{e}ma $S$, on notera suivant Deligne
$F^{(\lambda )}$ le produit tensoriel de $F$ par l'image
r\'{e}ciproque sur $S$ de $\overline{{\Bbb Q}}_{\ell}^{(\lambda )}$.
En particulier, $F^{(q^{-1})}$ n'est autre que la torsion \`{a} la
Tate $F(1)$ de $F$. Rappelons que tout syst\`{e}me local $\ell$-adique
$L$ sur un sch\'{e}ma $S$ normal de type fini peut s'\'{e}crire sous
la forme $L=M^{(\lambda )}$ o\`{u} $\lambda$ est une unit\'{e}
$\ell$-adique et $M$ est un syst\`{e}me local $\ell$-adique dont le
d\'{e}terminant (c'est-\`{a}-dire la puissance ext\'{e}rieure
maximale) est {\it d'ordre fini}.
\vskip 5mm

{\bf 6.2. La r\'{e}currence}
\vskip 5mm

Soient $N\subset X$ un niveau, $p:[0,r]\rightarrow {\Bbb R}$ un
param\`{e}tre de troncature assez convexe par rapport au genre de la
courbe $X$ et \`{a} $N$, et $d$ un entier. On a vu dans la section 3
comment Lafforgue construit un morphisme de champs alg\'{e}briques
$$
\mathop{\overline{\rm Cht}}\nolimits_{N}^{\,r,d;\,\leq
p}\rightarrow\mathop{\overline{\rm Tr}}\nolimits_{N}^{\,r,\tau}\times
(X-N)^{2}
$$
lisse purement de dimension relative $2r-2$, dont le compos\'{e} avec
la projection sur $(X-N)^{2}$ est propre, et dont la restriction
au-dessus de l'ouvert $\mathop{\rm Spec}({\Bbb F}_{q})=\mathop{\rm
Tr}\nolimits_{N}^{r,\tau}\subset \mathop{\overline{\rm
Tr}}\nolimits_{N}^{\,r,\tau}$ est le morphisme p\^{o}le-z\'{e}ro
$\mathop{\rm Cht}\nolimits_{N}^{r,d;\,\leq p}\rightarrow (X-N)^{2}$.

Par construction (cf. 3.4), le champ $\mathop{\widetilde{\rm
Tr}}\nolimits_{N}^{\,r,\tau}$ est lisse sur ${\cal T}^{r,\tau}
\times_{[{\Bbb A}^{r-1}/{\Bbb G}_{{\rm m}}^{r-1}],\langle
q-1\rangle}[{\Bbb A}^{r-1}/{\Bbb G}_{{\rm m}}^{r-1}]$ o\`{u} ${\cal
T}^{r,\tau}$ est un champ torique. On sait donc construire des
r\'{e}solutions des singularit\'{e}s $\mathop{\widetilde{\rm
Tr}}\nolimits_{N}^{\,r,\tau} \rightarrow\mathop{\overline{\rm
Tr}}\nolimits_{N}^{\,r,\tau}$. Fixons-en une et notons
$$
\mathop{\widetilde{\rm Cht}}\nolimits_{N}^{\,r,d;\,\leq p}=
\mathop{\overline{\rm Cht}}\nolimits_{N}^{\,r,d;\,\leq p}
\times_{\mathop{\overline{\rm Tr}}\nolimits_{N}^{\,r,\tau}}
\mathop{\widetilde{\rm Tr}}\nolimits_{N}^{\,r,\tau}.
$$

Le morphisme de champs
$$
\mathop{\widetilde{\rm Cht}}\nolimits_{N}^{\,r,d;\,\leq p}\rightarrow
(X-N)^{2}
$$
prolonge le morphisme p\^{o}le-z\'{e}ro $\mathop{\rm Cht}
\nolimits_{N}^{r,d;\,\leq p}\rightarrow (X-N)^{2}$ et est encore
propre. Il est de plus lisse purement de dimension relative $2r-2$, et
le ferm\'{e} compl\'{e}mentaire de l'ouvert $\mathop{\rm Cht}
\nolimits_{N}^{r,d; \,\leq p}\subset\mathop{\widetilde{\rm Cht}}
\nolimits_{N}^{\,r,d; \,\leq p}$ est un diviseur \`{a} croisements
normaux relatif sur $(X-N)^{2}$ qui est r\'{e}union de diviseurs
lisses sur $(X-N)^{2}$.

En prenant la r\'{e}union disjointe sur les entiers $d$ des champs
$\mathop{\widetilde{\rm Cht}} \nolimits_{N}^{\,r,d; \,\leq p}$ et en
divisant par l'action de $a^{{\Bbb Z}}$, on obtient un morphisme de
champs
$$
\mathop{\widetilde{\rm Cht}}\nolimits_{N}^{\,r;\,\leq 
p}/a^{{\Bbb Z}}\rightarrow (X-N)^{2}
$$
qui prolonge le morphisme p\^{o}le-z\'{e}ro $\mathop{\rm
Cht}\nolimits_{N}^{r;\,\leq p}/a^{{\Bbb Z}}\rightarrow (X-N)^{2}$,
qui est propre et lisse purement de dimension relative $2r-2$, et pour
lequel le ferm\'{e} compl\'{e}mentaire de l'ouvert $\mathop{\rm
Cht}\nolimits_{N}^{r; \,\leq p}/a^{{\Bbb Z}}\subset
\mathop{\widetilde{\rm Cht}}\nolimits_{N}^{\,r; \,\leq p}/a^{{\Bbb
Z}}$ est un diviseur \`{a} croisements normaux relatif sur
$(X-N)^{2}$ qui est r\'{e}union de diviseurs lisses sur $(X-N)^{2}$.

Le th\'{e}or\`{e}me suivant est une description partielle de la 
cohomologie $\ell$-adique des vari\'{e}t\'{e}s de chtoucas 
compactifi\'{e}es.

\thm TH\'{E}OR\`{E}ME 
\enonce
Pour tout niveau $N$ et tout param\`{e}tre de troncature $p$ assez
convexe par rapport au genre de la courbe $X$ et \`{a} $N$, les images
directes sup\'{e}rieures du faisceau constant de valeur
$\overline{{\Bbb Q}}_{\ell}$ par le morphisme canonique
$\mathop{\widetilde{\rm Cht}}\nolimits_{N}^{\,r; \,\leq p}/a^{{\Bbb
Z}}\rightarrow (X-N)^{2}$ sont toutes des syst\`{e}mes locaux
$\ell$-adiques $(r+1)$-n\'{e}gligeables.
\endthm

Lafforgue prouve le th\'{e}or\`{e}me principal de la section 1 et le
th\'{e}or\`{e}me ci-dessus par une r\'{e}currence simultan\'{e}e sur
$r$. Plus pr\'{e}cis\'{e}ment, pour tout entier $r'\geq 1$, il
consid\`{e}re les propri\'{e}t\'{e}s suivantes:
\vskip 2mm

\itemitem{}${\rm P}_{1}(r')$ Quel que soit $\pi'\in {\cal A}_{r'}$, il
existe $\sigma (\pi')\in {\cal G}_{r'}$ tel que l'on ait
l'\'{e}galit\'{e} de facteurs $L$ locaux
$$
\iota (L(\sigma (\pi')_{x},T))=L(\pi_{x}',T)
$$
pour presque toute place $x\notin S_{\sigma (\pi')}\cup S_{\pi'}$.
\vskip 2mm

\itemitem{}${\rm P}_{2}(r')$ Quel que soit $\sigma'\in {\cal G}_{r'}$,
il existe $\pi (\sigma')\in {\cal A}_{r'}$ tel que l'on ait
l'\'{e}galit\'{e} de facteurs $L$ locaux
$$
L(\pi (\sigma')_{x},T)=\iota (L(\sigma_{x}',T))
$$
pour presque toute place $x\notin S_{\pi (\sigma')}\cup S_{\sigma'}$.
\vskip 2mm

\itemitem{}${\rm P}_{3}(r')$ Quels que soient $\pi'\in {\cal A}_{r'}$
et $x\in |X|-S_{\pi'}$, les valeurs propres de Hecke $z_{1}(\pi'),
\ldots ,z_{r'}(\pi')$ sont toutes de valeur absolue $1$.
\vskip 2mm

\itemitem{}${\rm P}_{4}(r')$ Le th\'{e}or\`{e}me ci-dessus est
d\'{e}montr\'{e} en rang $r'$.
\vskip 3mm

\noindent Et il prouve:

\thm TH\'{E}OR\`{E}ME 
\enonce
Fixons un entier $r\geq 1$. Alors les propri\'{e}t\'{e}s ${\rm
P}_{1}(r')$, ${\rm P}_{3}(r')$ et ${\rm P}_{4}(r')$ pour $r'=1,\ldots
,r-1$ impliquent ces m\^{e}mes propri\'{e}t\'{e}s pour $r'=r$.
\endthm

Comme on l'a d\'{e}j\`{a} signal\'{e} en 1.4, les propri\'{e}t\'{e}s
${\rm P}_{1}(1),\ldots ,{\rm P}_{1}(r-1)$ impliquent les
propri\'{e}t\'{e}s ${\rm P}_{2}(1),\ldots ,{\rm P}_{2}(r-1),{\rm
P}_{2}(r)$, et en particulier la correspondance de Langlands en tout
rang $r'<r$. Le th\'{e}or\`{e}me ci-dessus implique donc bien le
th\'{e}or\`{e}me principal de 1.3. Ces m\^{e}mes propri\'{e}t\'{e}s
impliquent aussi des propri\'{e}t\'{e}s renforc\'{e}es o\`{u} l'on
impose que $S_{\sigma (\pi')}=S_{\pi'}$ et les \'{e}galit\'{e}s de
facteurs $L$ locaux pour toutes les places $x\notin S_{\pi'}$.

Dans la suite de cette section nous allons commencer la
d\'{e}monstration du th\'{e}or\`{e}me pr\'{e}c\'{e}dent en nous
concentrant sur les aspects galoisiens.
\vskip 5mm

{\bf 6.3. Les strates du bord sont $r$-n\'{e}gligeables}
\vskip 5mm

Dans toute cette sous-section on fixe un param\`{e}tre de troncature
$p:[0,r]\rightarrow {\Bbb R}$ assez convexe en fonction du genre de la
courbe $X$ et de $N$.

Pour tout entier $\nu$ on note
$$
H_{N}^{r;\,\leq p,\,\nu},~\widetilde{H}_{N}^{\,r;\,\leq 
p,\,\nu}\hbox{ et }\widetilde{H}_{N,\partial}^{\,r;\,\leq p,\,\nu}
$$
les $\nu$-i\`{e}mes images directes sup\'{e}rieures \`{a} supports
compacts du faisceau constant $\overline{{\Bbb Q}}_{\ell}$ par les
morphismes de champs
$$
\mathop{\rm Cht}\nolimits_{N}^{r;\,\leq p}/a^{{\Bbb Z}}\rightarrow
(X-N)^{2},~ \mathop{\widetilde{\rm Cht}}\nolimits_{N}^{\,r; \,\leq
p}/a^{{\Bbb Z}}\rightarrow (X-N)^{2}
$$
et
$$
\left(\mathop{\widetilde{\rm Cht}}\nolimits_{N}^{r;\,\leq p}
-\mathop{\rm Cht}\nolimits_{N}^{r; \,\leq p}\right)/a^{{\Bbb
Z}}\rightarrow (X-N)^{2}.
$$
Ce sont des syst\`{e}mes locaux $\ell$-adiques sur $(X-N)^{2}$, qui
sont nuls pour $\nu\notin\{0,1,\ldots ,4r-4\}$ (et m\^{e}me pour
$\nu\notin\{0,1,\ldots ,4r-6\}$ dans le cas de
$\widetilde{H}_{N,\partial}^{\,r;\,\leq p,\,\nu}$), et on a une suite
exacte longue
$$
\cdots\rightarrow H_{N}^{r;\,\leq p,\nu}\rightarrow
\widetilde{H}_{N}^{\,r;\,\leq p,\,\nu}\rightarrow
\widetilde{H}_{N,\partial}^{\,r;\,\leq p,\,\nu}\rightarrow
H_{N}^{r;\,\leq p,\,\nu +1}\rightarrow\cdots .
$$

Le diviseur \`{a} croisements normaux relatif
$\left(\mathop{\widetilde{\rm Cht}}\nolimits_{N}^{r;\,\leq p}
-\mathop{\rm Cht}\nolimits_{N}^{r; \,\leq p}\right)/a^{{\Bbb Z}}$ sur
$(X-N)^{2}$ est r\'{e}union d'une famille finie de diviseurs lisses.
Il est donc muni d'une stratification canonique. Une {\it strate
ferm\'{e}e du bord} est une intersection d'une sous-famille non vide
de cette famille. La relation d'inclusion d\'{e}finit une relation
d'ordre sur les strates ferm\'{e}es du bord. Une {\it strate ouverte
du bord} est le compl\'{e}mentaire dans une strate ferm\'{e}e du bord
de la r\'{e}union des strates ferm\'{e}es du bord strictement plus
petite. Par exemple, pour $N=\emptyset$, les strates ouvertes du bord
sont les $\mathop{\overline{\rm Cht}}\nolimits_{R}^{\,r;\,\leq
p}/a^{{\Bbb Z}}$ pour $\emptyset\not= R=\{r_{1},\ldots
,r_{s-1}\}\subset [r-1]$ d\'{e}finies en 3.3, et les strates
ferm\'{e}es du bord sont les adh\'{e}rences des strates ouvertes du
bord.

Toute strate du bord (ferm\'{e}e ou ouverte) est lisse sur
$(X-N)^{2}$, et les images directes sup\'{e}rieures \`{a} supports
compacts de $\overline{{\Bbb Q}}_{\ell}$ par la restriction \`{a} une
telle strate du morphisme p\^{o}le-z\'{e}ro sont toutes des
syst\`{e}mes locaux $\ell$-adiques sur $(X-N)^{2}$, qu'il est commode
d'appeler {\it images directes sup\'{e}rieures de la strate}.

\thm PROPOSITION 
\enonce
Supposons que les propri\'{e}t\'{e}s ${\rm P}_{1}(1)\ldots ,{\rm
P}_{1}(r-1)$ soient v\'{e}rifi\'{e}es pour ce niveau $N$. Alors, les
syst\`{e}mes locaux $\ell$-adiques
\vskip 2mm

\itemitem{-} images directes sup\'{e}rieures des strates ferm\'{e}es 
du bord,
\vskip 2mm

\itemitem{-} images directes sup\'{e}rieures des strates ouvertes 
du bord
\vskip 2mm

\itemitem{-} $\widetilde{H}_{N,\partial}^{\,r;\,\leq p,\,\nu}$,
$0\leq\nu\leq 4r-6$, 
\vskip 2mm

\noindent sont tous $r$-n\'{e}gligeables. 
\endthm

Esquissons la preuve de cette proposition pour $N=\emptyset$. On
proc\`{e}de par r\'{e}currence sur $r$. Compte tenu de la suite exacte
longue ci-dessus, les hypoth\`{e}ses de la proposition et de
r\'{e}currence assurent que les images directes sup\'{e}rieures par
les morphismes
$$
\mathop{\rm Cht}\nolimits^{r',d';\,\leq p'}\rightarrow X\times X
$$
sont toutes $(r'+1)$-n\'{e}gligeables quels que soient les entiers
$1\leq r'<r$ et $d'$ et le param\`{e}tre de troncature $p'$ assez
convexe par rapport au genre de $X$.

Chaque strate $\mathop{\overline{\rm Cht}}\nolimits_{R}^{r,d;\,\leq
p}$, $\emptyset\not= R=\{r_{1},\ldots ,r_{s-1}\}\subset [r-1]$, est
finie et radicielle sur un champ $\mathop{\rm Cht}
\nolimits^{R,d;\,\leq p}$ et lui est donc \'{e}quivalente du point de
vue de la cohomologie $\ell$-adique. Par hypoth\`{e}se de
r\'{e}currence, tous les sous-quotients irr\'{e}ductibles des images
directes sup\'{e}rieures par le morphisme
$$
(\infty =\infty_{1},o_{1}=\infty_{2},\ldots ,o_{s-1}=\infty_{s},
o_{s}=o):\mathop{\rm Cht}\nolimits^{R,d;\,\leq p}\rightarrow X\times
X^{s-1}\times X
$$
sont de la forme
$$
L_{1}'\boxtimes\left(\mathop{\boxtimes}_{j=1}^{s-1}(L_{j}''\otimes
L_{j+1}')\right)\boxtimes L_{s}''
$$
pour des syst\`{e}mes locaux $\ell$-adiques irr\'{e}ductibles
$L_{1}',L_{1}''$ de rangs $\leq r_{1}<r$, $L_{2}',L_{2}''$ de rangs
$\leq r_{2}-r_{1}<r$, ... et $L_{s}',L_{s}''$ de rangs $\leq
r-r_{s-1}<r$ sur la courbe $X$. Par suite, les images directes
sup\'{e}rieures par le morphisme $(\infty ,o): \mathop{\overline{\rm
Cht}}\nolimits_{R}^{r,d;\,\leq p}\rightarrow X\times X$ se
d\'{e}composent en syst\`{e}mes locaux $\ell$-adiques
$$
(L_{1}'\boxtimes L_{s}'')\otimes H^{\nu}\left(\overline{{\Bbb F}}_{q}
\otimes_{{\Bbb F}_{q}}X^{s-1},\mathop{\boxtimes}_{j=1}^{s-1}
(L_{j}''\otimes L_{j+1}')\right),
$$
et ont tous leurs sous-quotients irr\'{e}ductibles de la forme
$(L_{1}'\boxtimes L_{s}'')^{(\lambda)}$ o\`{u} $\lambda$ est une
valeur propre de $\mathop{\rm Frob}\nolimits_{q}$ agissant sur
$H^{\nu}\left(\overline{{\Bbb F}}_{q} \otimes_{{\Bbb
F}_{q}}X^{s-1},\boxtimes_{j=1}^{s-1} (L_{j}''\otimes
L_{j+1}')\right)$; la proposition s'en suit.
\vskip 2mm

Le cas g\'{e}n\'{e}ral est un peu plus difficile mais similaire en
utilisant le fait essentiel que les $\mathop{\overline{\rm Cht}}
\nolimits_{N}^{\,r;\,\leq p}/a^{{\Bbb Z}}$ et $\mathop{\widetilde{\rm
Cht}}\nolimits_{N}^{\,r;\,\leq p}/a^{{\Bbb Z}}$ sont construits \`{a}
partir des $\mathop{\overline{\rm Cht}}\nolimits^{\,r;\,\leq
p}/a^{{\Bbb Z}}$ en formant un produit fibr\'{e} au dessus du champ
$\mathop{\overline{\rm Tr}}\nolimits_{N}^{\,r}$. On a besoin d'une
formule de comptage suppl\'{e}mentaire o\`{u} interviennent les
fonctions d'Euler-Poincar\'{e} introduites par Kottwitz (cf. [Lau 1])
et un argument inspir\'{e} de [Ha-Ta].
\vskip 5mm

{\bf 6.4. S\'{e}paration de la partie cuspidale de la cohomologie}
\vskip 5mm

Soit $S$ un sch\'{e}ma de type fini sur ${\Bbb F}_{q}$. Un syst\`{e}me
local $\ell$-adique {\it virtuel} $M$ sur $U\times U$ est une
combinaison lin\'{e}aire formelle $M=\sum_{L}m_{M}(L)[L]$ o\`{u} $L$
parcourt un ensemble de repr\'{e}sentants des classes d'isomorphie de
syst\`{e}mes locaux $\ell$-adiques irr\'{e}ductibles sur $S$ (fix\'{e}
une fois pour toute) et o\`{u} $L\mapsto m_{M}(L)$ est une fonction
sur cet ensemble \`{a} valeurs {\it rationnelles} et \`{a} support
fini. Le coefficient $m_{M}(L)\in {\Bbb Q}$ est appel\'{e} la {\it
multiplicit\'{e}} de $L$ dans $M$. On dit que $L$ {\it intervient
dans} $M$ ou encore que $L$ {\it est un constituant de} $M$ si
$m_{M}(L)\not=0$. Un syst\`{e}me local $\ell$-adique $M$ d\'{e}finit
un syst\`{e}me local $\ell$-adique virtuel $[M]$: $m_{M}(L)$ est la
multiplicit\'{e} de $L$ dans n'importe quelle filtration de
Jordan-H\"{o}lder de $M$.

On dira qu'un syst\`{e}me local $\ell$-adique virtuel $M=\sum_{L}
m_{M}(L)[L]$ sur $S$ est {\rm effectif} si $m_{M}(L)\geq 0$ quel que
soit $L$. On dira qu'il est {\it pur de poids} un entier $w$ si tous
les $L$ qui interviennent dans $M$ sont ponctuellement purs de poids
$w$. On dira qu'il est {\it mixte de poids contenus dans} un ensemble
$W$ d'entiers si chaque $L$ qui intervient dans $M$ est ponctuellement
pur de poids appartenant \`{a} $W$.

Soit maintenant $U$ un ouvert non vide de $X$. Un syst\`{e}me local
$\ell$-adique {\it virtuel} $M$ sur $U\times U$ est dit {\it
$r$-n\'{e}gligeable} si tous les $L$ {\it qui interviennent dans} $M$
sont $r$-n\'{e}gligeables. Si $M=\sum_{L}m_{M}(L)[L]$ est un
syst\`{e}me local $\ell$-adique virtuel sur $U\times U$, on appelera
{\it morceau non $r$-n\'{e}gligeable} de $M$ et on notera $M_{{\rm
cusp}}= \sum_{L}m_{M_{{\rm cusp}}}(L)[L]$ le syst\`{e}me local
$\ell$-adique virtuel sur $U\times U$ d\'{e}fini par $m_{M_{{\rm
cusp}}}(L)=0$ si $L$ est $r$-n\'{e}gligeable et $m_{M_{{\rm
cusp}}}(L)=m_{M}(L)$ sinon.
\vskip 3mm

En confrontant la formule des points fixes de Grothendieck-Lefschetz
([Gr])
$$
\mathop{\rm Lef}\nolimits_{\xi}\left(\mathop{\rm Frob}
\nolimits_{\mathop{\rm Cht}\nolimits_{N}^{r;\,\leq p}/a^{{\Bbb
Z}}}^{\mathop{\rm deg}(\xi )n}\right)=\sum_{\nu}(-1)^{\nu}\mathop{\rm
Tr}(\mathop{\rm Frob}\nolimits_{\xi}^{n},H_{N}^{r;\leq p,\,\nu})
$$
et le cas particulier $g=1$ et $\mathop{\rm deg}(\infty )n_{\infty}=
\mathop{\rm deg}(o)n_{o}=\mathop{\rm deg}(\xi )n$ de la formule pour
la moyenne de nombres de Lefschetz \'{e}nonc\'{e}e en 5.2, on voit que
sans changer la valeur de l'expression
$$
\sum_{{\scriptstyle 1\leq r'<r\atop\scriptstyle 1\leq r''<r}}
\sum_{{\scriptstyle\pi'\in {\cal A}_{r'}(K_{N})\atop\scriptstyle
\pi''\in {\cal A}_{r''}(K_{N})}}\mathop{\rm Tr}
\nolimits_{\pi',\pi''}^{\leq_{\alpha}p}(1_{K_{N}},\mathop{\rm deg}(\xi
)n)S_{\infty}^{\left(-{\mathop{\rm deg}(\xi )\over \mathop{\rm
deg}(\infty )}n\right)}\! (\pi')S_{o}^{\left({\mathop{\rm deg}(\xi
)\over \mathop{\rm deg}(o)}n\right)}\!(\pi'')
$$
on peut y remplacer les fonctions $\nu\mapsto\mathop{\rm Tr}
\nolimits_{\pi',\pi''}^{\leq p}(1_{K_{N}},\nu )$ par des fonctions de
la forme $\sum_{z\in {\Bbb C}^{\times}}c_{z}z^{\nu}$ pour une famille
\`{a} support fini de scalaires $c_{z}\in {\Bbb Q}$, $z\in {\Bbb
C}^{\times}$. 

L'hypoth\`{e}se de r\'{e}currence permet de remplacer les $\pi'$ et
$\pi''$ qui interviennent dans cette m\^{e}me expression par des
syst\`{e}mes locaux $\ell$-adiques irr\'{e}ductibles et purs de poids
$0$ sur $X-N$. Compte tenu du th\'{e}or\`{e}me de puret\'{e} de
Deligne ([De]) on en d\'{e}duit:

\thm PROPOSITION 
\enonce
Si $p:[0,r]\rightarrow {\Bbb R}$ est un param\`{e}tre de troncature
assez convexe en fonction du genre de la courbe $X$ et de $N$, il
existe un syst\`{e}me local $\ell$-adique virtuel mixte de poids
entiers $\widetilde{H}_{N,{\rm cusp}}^{\,r;\,\leq p}$ sur $(X-N)^{2}$
tel que la diff\'{e}rence
$$
\widetilde{H}_{N,{\rm cusp}}^{\,r;\,\leq p}-{1\over r!}\sum_{k=1}^{r!}
\sum_{\nu}(-1)^{\nu}[(\mathop{\rm Frob}\nolimits_{X}^{k}
\times\mathop{\rm Id}\nolimits_{X})^{\ast}H_{N}^{r;\, \leq p,\,\nu}]
$$
soit $r$-n\'{e}gligeable et que, pour tout couple $(\infty ,o)$ de points
ferm\'{e}s {\rm distincts} de $X-N$, tout point ferm\'{e}
$\xi\in\infty\times o$ et tout entier $n\geq 1$, on ait
$$
\mathop{\rm Tr}(\mathop{\rm Frob}\nolimits_{\xi}^{n},
\widetilde{H}_{N,{\rm cusp}}^{\,r;\,\leq p})=\sum_{{\scriptstyle
\pi\in {\cal A}_{r}(K_{N}) \atop\scriptstyle\omega_{\pi}(a)=1}}
q^{(r-1)\mathop{\rm deg}(\xi )n}\mathop{\rm Tr}\nolimits_{\pi}
(1_{K_{N}}) S_{\infty}^{\left(-{\mathop{\rm deg}(\xi )\over
\mathop{\rm deg}(\infty)}n\right)}\!(\pi )S_{o}^{\left({\mathop{\rm
deg}(\xi )\over\mathop{\rm deg}(o)}n\right)}\!(\pi ).
$$
\endthm

L'encadrement de Jacquet et Shalika ([J-S 1])
$$
q^{-{1\over 2}}<|z_{i}(\pi_{x})|^{{1\over \mathop{\rm deg}(x)}}
<q^{{1\over 2}},~\forall i=1,\ldots ,r,~\forall x\in |X-N|,
$$
pour les valeurs propres de Hecke de tout $\pi\in {\cal A}_{r}(K_{N})$
assure que $\widetilde{H}_{N,{\rm cusp}}^{\,r;\,\leq p}$ est mixte de
poids contenus dans $\{2r-3,2r-2,2r-1\}$ et qu'il a des constituants
de poids $2r-3$ si et seulement si il en a de poids $2r-1$.

Par des arguments de fonctions $L$ pour les syst\`{e}me locaux
$\ell$-adiques sur la surface $(X-N)^{2}$, Lafforgue d\'{e}duit de la
proposition ci-dessus, du th\'{e}or\`{e}me de puret\'{e} de Deligne et
des r\'{e}sultats de 6.3:

\thm COROLLAIRE
\enonce
{\rm (i)} Le syst\`{e}me local $\ell$-adique virtuel
$\widetilde{H}_{N,{\rm cusp}}^{\,r;\,\leq p}$ sur $(X-N)^{2}$ est
effectif, pur de poids $2r-2$ et n'admet aucun constituant
$r$-n\'{e}gligeable.

\decale{\rm (ii)} Les syst\`{e}mes locaux $\ell$-adiques
$\widetilde{H}_{N}^{\,r;\,\leq p,\,\nu}$, $\nu\not=2r-2$, sont tous
$r$-n\'{e}gligeables et $\widetilde{H}_{N,{\rm cusp}}^{\,r;\,\leq p}$
est exactement le morceau non $r$-n\'{e}gligeable du syst\`{e}me local
$\ell$-adique virtuel
$$
{1\over r!}\sum_{k=1}^{r!} [(\mathop{\rm Frob}\nolimits_{X}^{k}\times
\mathop{\rm Id}\nolimits_{X})^{\ast}\widetilde{H}_{N}^{\,r;\,\leq
p,\,2r-2}].
$$
\endthm

Pour tous param\`{e}tres de troncatures $p\leq q$ qui sont assez
convexes par rapport au genre de $X$ et au niveau $N$, on a un
carr\'{e} commutatif
$$\diagram{
H_{N}^{r;\,\leq q,\,\ast}&\kern -1mm\smash{\mathop{\hbox to
8mm{\rightarrowfill}}\limits^{\scriptstyle }}\kern
-1mm&\widetilde{H}_{N}^{r;\,\leq q,\,\ast}\cr
\llap{$\scriptstyle $}\left\uparrow
\vbox to 4mm{}\right.\rlap{}&+&\llap{}\left\downarrow
\vbox to 4mm{}\right.\rlap{$\scriptstyle $}\cr
H_{N}^{r;\,\leq p,\,\ast}&\kern -1mm\smash{\mathop{\hbox to
8mm{\rightarrowfill}} \limits_{\scriptstyle }}\kern
-1mm&\widetilde{H}_{N}^{r;\,\leq p,\,\ast}\cr}
$$
o\`{u} les deux fl\`{e}ches horizontales et la fl\`{e}che verticale
montante de gauche sont les prolongements par z\'{e}ro pour les
inclusions de $\mathop{\rm Cht}\nolimits_{N}^{r;\,\leq p}\subset
\mathop{\widetilde{\rm Cht}}\nolimits_{N}^{\,r;\,\leq p}$,
$\mathop{\rm Cht}\nolimits_{N}^{r;\,\leq q}\subset
\mathop{\widetilde{\rm Cht}}\nolimits_{N}^{\,r;\,\leq q}$ et
$\mathop{\rm Cht}\nolimits_{N}^{r;\,\leq p}\subset \mathop{\rm
Cht}\nolimits_{N}^{\,r;\,\leq q}$ et la fl\`{e}che verticale
descendante de droite est induite par la correspondance dans
$\mathop{\widetilde{\rm Cht}}\nolimits_{N}^{\,r;\,\leq p}\times
\mathop{\widetilde{\rm Cht}}\nolimits_{N}^{\,r;\,\leq q}$ qui est
l'adh\'{e}rence du graphe de l'inclusion $\mathop{\rm
Cht}\nolimits_{N}^{r;\,\leq p} \subset\mathop{\rm
Cht}\nolimits_{N}^{r;\,\leq q}$. La proposition de 6.3, le corollaire
ci-dessus et l'existence de ce carr\'{e} commutatif impliquent
facilement le lemme suivant qui est crucial pour la suite.

\thm LEMME 
\enonce
Pour tous param\`{e}tres de troncatures $p\leq q$ qui sont assez 
convexes par rapport au genre de $X$ et au niveau $N$, les noyau et 
conoyau de l'homomorphisme de prolongement par z\'{e}ro
$$
H_{N}^{r;\,\leq p,\,2r-2}\rightarrow H_{N}^{r;\,\leq q,\,2r-2}
$$
sont $r$-n\'{e}gligeables.
\endthm

En utilisant ce lemme, Lafforgue prouve que le ind-syst\`{e}me local
$\ell$-adique {\og}{de rang infini}{\fg} sur $(X-N)^{2}$
$$
H_{N}^{r;\,2r-2}=\lim_{{\scriptstyle\longrightarrow\atop
\scriptstyle p}}H_{N}^{r;\,\leq p,\,2r-2},
$$
o\`{u} $p$ parcourt l'ensemble des param\`{e}tres de troncature qui
sont assez convexes par rapport au genre de la courbe $X$ et \`{a}
$N$, a la propri\'{e}t\'{e} suivante:

\thm LEMME
\enonce
Il existe une unique filtration finie
$$
F^{\bullet}=((0)=F^{0}\subseteq F^{1}\subsetneq F^{2}\subsetneq\cdots
\subsetneq F^{2u+1}\subsetneq F^{2u}\subsetneq\cdots\subsetneq
F^{T}=H_{N}^{r;\,2r-2}
$$
telle que:
\vskip 2mm

\itemitem{-} pour tout entier $u\geq 0$ tel que $2u+1\leq T$,
$F^{2u+1}/F^{2u}$ est la somme de tous les sous-syst\`{e}mes locaux
$\ell$-adiques {\rm (}de rang fini{\rm )} $r$-n\'{e}gligeables de
$H_{N}^{r;\,2r-2}/F^{2u}$,
\vskip 2mm

\itemitem{-} pour tout entier $u\geq 0$ tel que $2u+2\leq T$,
$F^{2u+2}/F^{2u+1}$ est la somme de tous les sous-syst\`{e}mes locaux
$\ell$-adiques {\rm (}de rang fini{\rm )} $H_{N}^{r;\,2r-2}/F^{2u+1}$
dont aucun sous-quotient n'est $r$-n\'{e}gligeable,
\vskip 2mm

\itemitem{-} si $p$ est un param\`{e}tre de troncature assez convexe
par rapport au genre de $X$ et au niveau $N$ et si on note $F^{\leq
p,\,\bullet}$ la filtration sur $H_{N}^{r;\,\leq p,\,2r-2}$ induite
par $F^{\bullet}$, alors, pour tout entier $u\geq 0$, le plongement
$$
F^{\leq p,\,2u+2}/F^{\leq p,\,2u+1}\hookrightarrow F^{2u+2}/F^{2u+1}
$$
est un isomorphisme.

\endthm

Consid\'{e}rons alors le syst\`{e}me local $\ell$-adique (de rang fini) sur
$(X-N)^{2}$
$$
H_{N,{\rm cusp}}^{r}=\bigoplus_{u\geq 0}F^{2u+2}/F^{2u+1}.
$$
Rassemblant les r\'{e}sultats de cette sous-section, on obtient
finalement:

\thm PROPOSITION 
\enonce
{\rm (i)} Pour tout param\`{e}tre de troncature $p$ assez convexe par
rapport au genre de $X$ et au niveau $N$, on a l'\'{e}galit\'{e} de
syst\`{e}mes locaux $\ell$-adiques virtuels
$$
\widetilde{H}_{N,{\rm cusp}}^{\,r;\,\leq p}={1\over r!}\sum_{k=1}^{r!}
[(\mathop{\rm Frob}\nolimits_{X-N}^{k}\times\mathop{\rm
Id}\nolimits_{X-N})^{\ast}H_{N,{\rm cusp}}^{r}].
$$

\decale{\rm (ii)} Le syst\`{e}me local $\ell$-adique $H_{N,{\rm
cusp}}^{r}$ est ponctuellement pur de poids $2r-2$.

\decale{\rm (iii)} Pour tout couple $(\infty ,o)$ de points ferm\'{e}s
distincts de $X-N$, tout point ferm\'{e} $\xi$ de $\infty\times o$ et
tout entier $n$, on a
$$\displaylines{
\qquad {1\over r!}\sum_{k=1}^{r!}\mathop{\rm Tr}(\mathop{\rm Frob}
\nolimits_{\xi}^{n},(\mathop{\rm Frob}\nolimits_{X-N}^{k}\times
\mathop{\rm Id}\nolimits_{X-N})^{\ast}H_{N,{\rm cusp}}^{r})
\hfill\cr\hfill
=q^{(r-1)\mathop{\rm deg}(\xi )n}\kern -2mm
\sum_{{\scriptstyle\pi\in {\cal A}_{r}(K_{N})\atop\scriptstyle
\omega_{\pi}(a)=1}}\mathop{\rm Tr}\nolimits_{\pi}(1_{K_{N}})
S_{\infty}^{\left(-{\mathop{\rm deg}(\xi )\over\mathop{\rm deg}
(\infty )}n\right)}\!(\pi ) S_{o}^{\left({\mathop{\rm deg}(\xi )\over
\mathop{\rm deg}(o)}n\right)}\!(\pi ).\qquad}
$$
\endthm

\vskip 7mm

\centerline{\bf 7. UNE VARIANTE D'UN TH\'{E}OR\`{E}ME DE PINK}
\vskip 10mm

{\bf 7.1. Le cas des sch\'{e}mas} 
\vskip 5mm

Soient $S$ un sch\'{e}ma de type fini et lisse (sur ${\Bbb F}_{q}$) et
$p:X\rightarrow S$ un morphisme de sch\'{e}mas propre et lisse
purement de dimension relative $d$. On se donne un diviseur $Y\subset
X$ \`{a} croisements normaux relatif sur $S$ qui est r\'{e}union
$Y=Y_{1}\cup\cdots\cup Y_{m}$ de $m$ diviseurs lisses sur $S$. Pour
toute partie $I$ de $\{1,\ldots ,m\}$ on note $Y_{I}=\bigcap_{i\in
I}Y_{i}$ et $X_{I}=Y_{I}-\bigcup_{J\supsetneq I}Y_{J}$. On a donc
$X_{\emptyset}=X-Y\subset Y_{\emptyset}=X$ et $X$ est la r\'{e}union
disjointe de ses parties localement ferm\'{e}es $X_{I}$. De plus, pour
chaque $I\subset\{1,\ldots ,m\}$, la restriction $p_{I}$ de $p$
\`{a} $Y_{I}$ est propre et lisse purement de dimension relative 
$d-|I|$, et $X_{I}$ est un ouvert dense de $Y_{I}$.

On consid\`{e}re un endomorphisme $f:S\rightarrow S$ et un morphisme 
fini de sch\'{e}mas
$$
\Gamma_{\emptyset}\rightarrow Z_{\emptyset}:=X_{\emptyset}
\times_{f,S}X_{\emptyset}\subset X_{\emptyset} \times X_{\emptyset}
$$
o\`{u} {\og}{$\times_{f,S}$}{\fg} est une notation abr\'{e}g\'{e}e
pour {\og}{$\times_{p_{\emptyset}\circ f,S,p_{\emptyset}}$}{\fg}. On
suppose que la premi\`{e}re projection $\mathop{\rm
pr}\nolimits_{\emptyset ,1}:\Gamma_{\emptyset} \rightarrow
X_{\emptyset}$ est \'{e}tale. Le $S$-sch\'{e}ma est donc lisse
purement de dimension relative $d$ sur $S$. En particulier il est
normal.

Alors, si $s$ est un point ferm\'{e} de $S$ et si $n\geq 1$ est un
entier tel que
$$
\mathop{\rm Frob}\nolimits_{S}^{n}(f(s))=f(\mathop{\rm Frob}
\nolimits_{S}^{n}(s))=s,
$$
la fibre $\Gamma_{\emptyset ,s}\rightarrow Z_{\emptyset ,s}\subset
X_{\emptyset ,f(s)}\times X_{\emptyset ,s}$ en $s$ de
$\Gamma_{\emptyset}$ coupe transversalement la fibre $X_{\emptyset
,f(s)}\rightarrow X_{\emptyset ,s}\times X_{\emptyset ,f(s)}$ en
$(s,f(s))$ du graphe $\delta_{\emptyset}^{n}=(\mathop{\rm Frob}
\nolimits_{X_{\emptyset}}^{n},\mathop{\rm Id}
\nolimits_{X_{\emptyset}}): X_{\emptyset}\rightarrow
X_{\emptyset}\times X_{\emptyset}$. On note
$$
\mathop{\rm Lef}\nolimits_{s}(\Gamma_{\emptyset}\times\mathop{\rm
Frob}\nolimits_{X_{\emptyset}}^{n},X_{\emptyset})
$$
le nombre des points d'intersection.

On note $\Gamma\rightarrow Z:=X\times_{f,S}X\subset X\times X$ le
morphisme fini de sch\'{e}mas obtenu par normalisation de $Z$ dans
$\Gamma_{\emptyset}$. On se propose de donner une interpr\'{e}tation
cohomologique de $\mathop{\rm Lef}\nolimits_{s}
(\Gamma_{\emptyset}\times\mathop{\rm Frob}\nolimits_{X}^{n},
X_{\emptyset})$ sous une hypoth\`{e}se de {\og}{stabilit\'{e} de
l'ouvert $X_{\emptyset}\subset X$ au voisinage des points fixes de
$\Gamma$}{\fg}.
\vskip 5mm

{\bf 7.2. Stabilit\'{e} au voisinage des points fixes}
\vskip 5mm

On note $\mathop{\rm pr}\nolimits_{1},\mathop{\rm pr}\nolimits_{2}:
\Gamma\rightarrow X$ les deux projections canoniques de $\Gamma$ et,
pour chaque entier $n\geq 0$, $\delta^{n}$ le graphe $(\mathop{\rm
Frob}\nolimits_{X}^{n},\mathop{\rm Id} \nolimits_{X}): X\rightarrow
X\times X$.

\thm D\'{E}FINITION 
\enonce
On dira que $\Gamma$ stabilise $X_{\emptyset}$ au voisinage de ses
points fixes s'il existe un ouvert $U\subset X\times X$ contenant
toutes les intersections de l'image de $\Gamma$ dans $Z\subset X\times
X$ avec les images des graphes de Frobenius $\delta^{n}(X)\subset
X\times X$, $n\geq 0$, et tel que
$$
\mathop{\rm pr}\nolimits_{1}(\mathop{\rm
pr}\nolimits_{2}^{-1}(X_{\emptyset})\cap U_{\Gamma})\subset
X_{\emptyset}
$$
o\`{u} $U_{\Gamma}\subset\Gamma$ est la trace de $U$ sur $\Gamma$.
\endthm

\thm LEMME
\enonce
Pour que $\Gamma$ stabilise $X_{\emptyset}$ au voisinage de ses
points fixes, il faut et il suffit que $\Gamma$ satisfasse au crit\`{e}re 
valuatif suivant:
\vskip 2mm

$(\ast)$ Pour tout anneau de valuation discr\`{e}te $A$ de corps des
fractions $K$ et de corps r\'{e}siduel $k$ et pour tout point $\gamma
:\mathop{\rm Spec}(A)\rightarrow\Gamma$ tel que $\mathop{\rm
pr}\nolimits_{1}(\gamma (\mathop{\rm Spec}(K)))\in Y(K)$ et
$\mathop{\rm pr}\nolimits_{2}(\gamma (\mathop{\rm Spec}(K)))\in
X_{\emptyset}(K)$, l'image de $\gamma (\mathop{\rm Spec}(k))\in\Gamma
(k)$ dans $Z\subset X\times X$ n'est dans l'image d'aucun des graphes
de Frobenius $\delta^{n}$, $n\in {\Bbb N}$.
\endthm

Comme l'a fait Pink dans [Pi], Lafforgue introduit l'\'{e}clatement
$\pi :\widetilde{Z}\rightarrow Z$ de $Z=X\times_{f,S}X$ le long du
ferm\'{e} r\'{e}union des $Y_{i}\times_{f,S}Y_{i}$, $1\leq i\leq m$.
Cet \'{e}clatement est aussi le produit fibr\'{e} au-dessus de $Z$ des
\'{e}clatements $\widetilde{Z}_{i}\rightarrow Z$, $1\leq i\leq m$, de
$Z$ le long des $Y_{i}\times_{f,S}Y_{i}$. Il est donc muni de
diviseurs $E_{i}$, $1\leq i\leq m$, qui sont les images
r\'{e}ciproques dans $\widetilde{Z}$ des diviseurs exceptionnels dans
les $\widetilde{Z}_{i}$. Pour toute partie non vide
$I\subset\{1,\ldots ,m\}$, on v\'{e}rifie que
$$
E_{I}:=\bigcap_{i\in I}E_{i}\subset\bigcap_{i\in I}\pi^{-1}
(Y_{i}\times_{f,S}Y_{i})=\pi^{-1}(Y_{I}\times_{f,S}Y_{I}),
$$
que $E_{I}$ est le produit fibr\'{e} sur $Z$ des diviseurs
exceptionnels dans les $\widetilde{Z}_{i}$ pour $i\in I$ et des
$\widetilde{Z}_{i}$ pour $i\notin I$, que $E_{I}$ est lisse sur $S$
purement de dimension relative $2d-|I|$ et que la restriction de
$\pi_{I}:E_{I}\rightarrow Y_{I}\times_{f,S}Y_{I}$ au-dessus de
l'ouvert $X_{I}\times_{f,S}X_{I}$ est un fibr\'{e} projectif de rang
$|I|$. Les $E_{i}$ sont donc des diviseurs lisses sur $S$ et leur
r\'{e}union $E=E_{1}\cup\cdots \cup E_{m}$ est un diviseur \`{a}
croisements normaux relatif sur $S$.

Par construction $\pi :\widetilde{Z}\rightarrow Z$ est un isomorphisme
au dessus de $Z_{\emptyset}$. On notera $\widetilde{\Gamma}\rightarrow
\widetilde{Z}$ la normalisation de $\Gamma_{\emptyset}$ dans
$\widetilde{Z}$. Bien entendu $\pi$ {\og}{envoie}{\fg}
$\widetilde{\Gamma}$ dans $\Gamma$.

Pour tout point ferm\'{e} $s\in S$ la fibre $\pi_{s}:\widetilde{Z}_{s}
\rightarrow Z_{s}$ en $s$ de $\pi$ est l'\'{e}clatement de $Z_{s}$ le
long de la r\'{e}union des ferm\'{e}s $Y_{i,s}\times_{f,f(s)}
Y_{i,f(s)}$ de $Z_{s}=X_{s}\times_{f,f(s)}X_{f(s)}$. Pour
tout point ferm\'{e} $s\in S$ et tout entier $n\geq 0$ tel que
$\mathop{\rm Frob}\nolimits_{S}^{n}(f(s))=f(\mathop{\rm
Frob}\nolimits_{S}^{n}(s))=s$, la fibre $\delta_{(s,f(s))}^{n}:
X_{f(s)}\rightarrow Z_{s}\subset X_{s}\times X_{f(s)}$ en $(s,f(s))$
du graphe de Frobenius se rel\`{e}ve canoniquement en un morphisme
$$
\widetilde{\delta}_{(s,f(s))}^{n}:X_{f(s)}\rightarrow\widetilde{Z}_{s}
$$
puisque, pour chaque $i=1,\ldots ,m$, l'image r\'{e}ciproque par
$\delta_{(s,f(s))}^{n}$ du ferm\'{e} $Y_{i,s}\times_{f,f(s)}
Y_{i,f(s)}\subset Z_{s}$ est le diviseur $Y_{i,f(s)}\subset X_{f(s)}$.

La proposition suivante prolonge un r\'{e}sultat de Pink ([Pi]).

\thm PROPOSITION 
\enonce
Supposons que $\Gamma$ stabilise $X_{\emptyset}$ au voisinage de ses
points fixes. Alors, quitte \`{a} remplacer $f$ par
$$
f^{(n_{0})}:=\mathop{\rm Frob} \nolimits_{S}^{n_{0}}\circ
f=f\circ\mathop{\rm Frob} \nolimits_{S}^{n_{0}}
$$
et $\Gamma_{\emptyset}$ par
$$
\Gamma_{\emptyset}^{(n_{0})}:=(X_{\emptyset}\times X_{\emptyset})
\times_{\mathop{\rm Frob}\nolimits_{X_{\emptyset}}^{n_{0}}\times
\mathop{\rm Id}\nolimits_{X_{\emptyset}},X_{\emptyset}\times
X_{\emptyset}}\Gamma_{\emptyset}.
$$
pour un entier $n_{0}\geq 0$ assez grand, la normalisation
$\widetilde{\Gamma}$ v\'{e}rifie la propri\'{e}t\'{e} suivante: pour
tout point ferm\'{e} $s\in S$ et tout entier $n>0$ tel que
$\mathop{\rm Frob}\nolimits_{S}^{n}(f(s))=f(\mathop{\rm Frob}
\nolimits_{S}^{n}(s))=s$, l'image de $\widetilde{\delta}_{f,f(s)}^{n}$
dans $\widetilde{Z}_{s}$ ne rencontre pas l'image de
$\widetilde{\Gamma}_{s}$ en dehors de $Z_{\emptyset ,s}\subset
\widetilde{Z}_{s}$.
\endthm
\vskip 5mm

{\bf 7.3. Formule des points fixes}
\vskip 5mm

Rappelons que l'on a fix\'{e} un nombre premier $\ell$ distinct de la
caract\'{e}ristique de ${\Bbb F}_{q}$. Consid\'{e}rons la cohomologie
$\ell$-adique de $X$ et $Z$ au-dessus de $S$
$$
Rp_{\ast}{\Bbb Q}_{\ell}\hbox{ et }Rq_{\ast}{\Bbb Q}_{\ell}\in D_{{\rm
c}}^{{\rm b}}(S,{\Bbb Q}_{\ell})
$$
o\`{u} $q:X\times_{f,S}X\rightarrow S$ est la projection canonique, et
plus g\'{e}n\'{e}ralement les cohomologies $\ell$-adique des $Y_{I}$
au-dessus de $S$
$$
Rp_{I,\ast}{\Bbb Q}_{\ell}\hbox{ et }Rq_{I,\ast}{\Bbb Q}_{\ell}\in
D_{{\rm c}}^{{\rm b}}(S,{\Bbb Q}_{\ell})
$$
o\`{u} $q_{I}:Y_{I}\times_{f,S}Y_{I}\rightarrow S$ est la projection
canonique. La formation de ces cohomologies commute aux changements
de base $S'\rightarrow S$. En particulier, la fibre en un point
g\'{e}om\'{e}trique $\overline{s}$ de $S$ de $Rp_{I,\ast}{\Bbb
Q}_{\ell}$ est la cohomologie $\ell$-adique $R\Gamma
(Y_{I,\overline{s}},{\Bbb Q}_{\ell})$ de la fibre de $p_{I}$ en ce
point.

Nos hypoth\`{e}ses assurent que les faisceaux de cohomologie
$R^{j}p_{I,\ast}{\Bbb Q}_{\ell}$ (resp. $R^{j}q_{I,\ast}{\Bbb
Q}_{\ell}$) sont tous lisses sur $S$ et qu'ils sont nuls pour les $j$
qui ne sont pas dans l'intervalle $[0,2(d-|I|)]$ (resp.
$[0,4(d-|I|]$).

Soient $I\subset\{0,1,\ldots ,m\}$, $n$ un entier $\geq 1$ et
$u_{I}\in H^{0}(S,R^{2(d-|I|)}q_{I,\ast}{\Bbb Q}_{\ell})(d-|I|)$ une
classe de cohomologie. Alors
\vskip 2mm

\itemitem{-} l'endomorphisme de Frobenius $\mathop{\rm
Frob}\nolimits_{Y_{I}}^{n}$ qui rel\`{e}ve $\mathop{\rm
Frob}\nolimits_{S}^{n}$, induit un homomorphisme
$$
\mathop{\rm Frob}\nolimits_{Y_{I}}^{n}:(\mathop{\rm Frob}
\nolimits_{S}^{n})^{\ast}Rp_{I,\ast}{\Bbb Q}_{\ell}\rightarrow
Rp_{I,\ast}{\Bbb Q}_{\ell},
$$
\vskip 2mm

\itemitem{-} d'apr\`{e}s la formule de K\"{u}nneth et la dualit\'{e} de
Poincar\'{e}, $u_{I}$ peut \^{e}tre vu comme un homomorphisme
$$
u_{I}:f^{\ast}Rp_{I,\ast}{\Bbb Q}_{\ell}\rightarrow Rp_{I,\ast}{\Bbb
Q}_{\ell},
$$
\vskip 2mm

\itemitem{-} les homomorphismes compos\'{e}s
$$
\mathop{\rm Frob}\nolimits_{Y_{I}}^{n}\circ (\mathop{\rm Frob}
\nolimits_{S}^{n})^{\ast}(u_{I}): (\mathop{\rm Frob}
\nolimits_{S}^{n})^{\ast}f^{\ast}Rp_{I,\ast}{\Bbb Q}_{\ell}
\rightarrow Rp_{I,\ast}{\Bbb Q}_{\ell}
$$
et
$$
u_{I}\circ f^{\ast}(\mathop{\rm Frob}\nolimits_{Y_{I}}^{n}):
f^{\ast}(\mathop{\rm Frob}\nolimits_{S}^{n})^{\ast} Rp_{I,\ast}{\Bbb
Q}_{\ell}\rightarrow Rp_{I,\ast}{\Bbb Q}_{\ell}
$$
co\"{\i}ncident.
\vskip 2mm

\noindent On notera 
$$
u_{I}\times\mathop{\rm Frob}\nolimits_{Y_{I}}^{n}: (\mathop{\rm
Frob}\nolimits_{S}^{n})^{\ast}f^{\ast}Rp_{I,\ast}{\Bbb Q}_{\ell}=
f^{\ast}(\mathop{\rm Frob}\nolimits_{S}^{n})^{\ast} Rp_{I,\ast}{\Bbb
Q}_{\ell}\rightarrow Rp_{I,\ast}{\Bbb Q}_{\ell}
$$
ces homomorphismes compos\'{e}s. Pour tout point g\'{e}om\'{e}trique
$\overline{s}$ de $S$ qui v\'{e}rifie $f(\mathop{\rm
Frob}\nolimits_{S}^{n} (\overline{s}))=\mathop{\rm
Frob}\nolimits_{S}^{n} (f(\overline{s}))=\overline{s}$,
$u_{I}\times\mathop{\rm Frob}\nolimits_{Y_{I}}^{n}$ induit un
endomorphisme de la cohomologie $R\Gamma (Y_{I,\overline{s}},{\Bbb
Q}_{\ell})$ de la fibre en $\overline{s}$ de $p_{I}:Y_{I}\rightarrow
S$ dont on notera
$$
\mathop{\rm tr}(u_{I}\times\mathop{\rm
Frob}\nolimits_{Y_{I}}^{n},R\Gamma (Y_{I,\overline{s}},{\Bbb
Q}_{\ell}))\in {\Bbb Q}_{\ell}
$$
la trace.
\vskip 3mm

Le {\og}{cycle de codimension $d$}{\fg} $\Gamma$ sur $X\times_{f,S}X$
admet une classe de cohomologie
$$
\mathop{\rm cl}(\Gamma)\in H^{0}(S,R^{2d}q_{\ast}{\Bbb Q}_{\ell})(d).
$$
De m\^{e}me, si on note $\widetilde{q}= q\circ\pi :\widetilde{Z}
\rightarrow S$, $\widetilde{\Gamma}$ admet une classe de cohomologie
$$
\mathop{\rm cl}(\widetilde{\Gamma})\in H^{0}(S,R^{2d}
\widetilde{q}_{\ast}{\Bbb Q}_{\ell})(d).
$$
Pour chaque $I\subset\{1,\ldots ,m\}$ non vide, on a un homomorphisme
de faisceaux
$$
R^{2d}\widetilde{q}_{\ast}{\Bbb Q}_{\ell}(d)\rightarrow
R^{2d}\widetilde{q}_{I,\ast}{\Bbb Q}_{\ell}(d)\rightarrow
R^{2(d-|I|)}q_{I,\ast}{\Bbb Q}_{\ell}(d-|I|)
$$
o\`{u} $\widetilde{q}_{I}:E_{I}\rightarrow S$ est la projection 
canonique, o\`{u} la premi\`{e}re fl\`{e}che est l'homomorphisme de 
restriction au ferm\'{e} $E_{I}\subset \widetilde{Z}$ et o\`{u} la 
seconde fl\`{e}che est duale de l'homomorphisme
$$
\pi_{I}^{\ast}:R^{2(d-|I|)}q_{I,\ast}{\Bbb Q}_{\ell}(d-|I|)\rightarrow
R^{2(d-|I|)}\widetilde{q}_{I,\ast}{\Bbb Q}_{\ell}(d-|I|).
$$
On note $\mathop{\rm cl}(\Gamma )_{I}$ l'image de $\mathop{\rm
cl}(\widetilde{\Gamma})$ par cet homomorphisme compos\'{e}. 

Plus g\'{e}n\'{e}ralement on peut refaire cette construction apr\`{e}s
avoir remplac\'{e} $f$ par $f^{(n_{0})}$ et $\Gamma_{\emptyset}$ par
$\Gamma_{\emptyset}^{(n_{0})}$ pour un entier $n_{0}\geq 0$ (cf. la
proposition de 7.2). Pour chaque $I\subset\{1,\ldots ,m\}$ non vide,
on obtient une section globale $\mathop{\rm cl}(\Gamma^{(n_{0})})_{I}$
de $R^{2(d-|I|)}q_{I,\ast}^{(n_{0})}{\Bbb Q}_{\ell}(d-|I|)$, o\`{u}
$q_{I}^{(n_{0})}:Y_{I}\times_{f^{(n_{0})},S}Y_{I} \rightarrow S$ est
la projection canonique. L'image directe de cette section globale par
le morphisme radiciel $\mathop{\rm Frob}\nolimits_{Y_{I}}^{n_{0}}
\times\mathop{\rm Id}\nolimits_{Y_{I}}:Y_{I}\times_{f^{(n_{0})},S}
Y_{I}\rightarrow Y_{I}\times_{f,S}Y_{I}$ est une section globale
not\'{e}e $\mathop{\rm cl}(\Gamma )_{I}^{(n_{0})}$ de
$R^{2(d-|I|)}q_{I,\ast}{\Bbb Q}_{\ell}(d-|I|)$.
\vskip 3mm

Avec ces notations, Lafforgue prouve la g\'{e}n\'{e}ralisation
suivante de la formule des points fixes de Grothendieck-Lefschetz
([Gr]).

\thm TH\'{E}OR\`{E}ME 
\enonce
Supposons que $\Gamma$ stabilise l'ouvert $X_{\emptyset}$ au voisinage
de ses points fixes et soit $n_{0}\geq 0$ un entier v\'{e}rifiant la
conclusion de la proposition de $7.3$.

Alors, pour tout point ferm\'{e} $s$ de $S$ et tout entier $n>n_{0}$
tels que $\mathop{\rm Frob}\nolimits_{S}^{n}(f(s)) =f(\mathop{\rm
Frob}\nolimits_{S}^{n}(s))=s$, on a
$$\eqalign{
\mathop{\rm Lef}\nolimits_{s}(\Gamma_{\emptyset}\times \mathop{\rm
Frob}\nolimits_{X_{\emptyset}}^{n},X_{\emptyset})&=\mathop{\rm
tr}(\mathop{\rm cl}(\Gamma)\times\mathop{\rm Frob}\nolimits_{X}^{n},
R\Gamma (X_{\overline{s}},{\Bbb Q}_{\ell}))\cr
&\kern 4mm +\sum_{{\scriptstyle I\subset\{1,\ldots ,m\}\atop
\scriptstyle I\not=\emptyset}}(-1)^{|I|}\mathop{\rm tr}(\mathop{\rm
cl}(\Gamma )_{I}^{(n_{0})}\times\mathop{\rm Frob}
\nolimits_{Y_{I}}^{n},R\Gamma (Y_{I,\overline{s}},{\Bbb
Q}_{\ell})).\cr}
$$
\endthm
\vskip 5mm

{\bf 7.4. Extension \`{a} certains champs alg\'{e}briques} 
\vskip 5mm

Comme le montre Lafforgue, les r\'{e}sultats de cette section valent
encore apr\`{e}s avoir remplac\'{e} le $S$-sch\'{e}ma $p:X\rightarrow
S$, le diviseur $Y=Y_{1}\cup\cdots\cup Y_{m}$ et le morphisme de
sch\'{e}mas $\Gamma_{\emptyset}\rightarrow Z_{\emptyset}$ par un
$S$-champ alg\'{e}brique $p:{\cal X}\rightarrow S$, un diviseur ${\cal
Y}$ dans ${\cal X}$ et un morphisme de champs alg\'{e}briques
$\Gamma_{\emptyset}\rightarrow {\cal Z}_{\emptyset}$ qui v\'{e}rifient
les propri\'{e}t\'{e}s suivantes:
\vskip 2mm

\itemitem{-} le morphisme $p:{\cal X}\rightarrow S$ est de type fini,
\vskip 2mm

\itemitem{-} au voisinage de chacun de ses points le champ ${\cal X}$
admet un recouvrement fini et plat par un espace alg\'{e}brique muni
de l'action d'un groupe fini qui agit transitivement sur les fibres de
ce recouvrement,
\vskip 2mm

\itemitem{-} le morphisme $p$ satisfait au crit\`{e}re valuatif de
propret\'{e},
\vskip 2mm

\itemitem{-} le morphisme $p$ est lisse purement de dimension relative
$d$,
\vskip 2mm

\itemitem{-} ${\cal Y}$ est une r\'{e}union de diviseurs ${\cal
Y}_{1},\ldots ,{\cal Y}_{m}$ qui sont lisses sur $S$ et est \`{a}
croisement normaux relatif sur $S$,
\vskip 2mm

\itemitem{-} ${\cal X}_{\emptyset}$ est un $S$-champ alg\'{e}brique de
Deligne-Mumford, le morphisme $\Gamma_{\emptyset}\rightarrow {\cal
Z}_{\emptyset}= {\cal X}_{\emptyset}\times_{f,S}{\cal X}_{\emptyset}$
est repr\'{e}sentable fini et la premi\`{e}re projection de
$\Gamma_{\emptyset}$ sur ${\cal X}_{\emptyset}$ est \'{e}tale.

\vskip 7mm

\centerline{\bf 8. FIN DE LA R\'{E}CURRENCE}
\vskip 10mm

{\bf 8.1. Ce qu'il reste \`{a} d\'{e}montrer}
\vskip 5mm

Soient $N\subset X$ un niveau. Pour tout $g\in\mathop{\rm GL}
\nolimits_{r}({\Bbb A})$ la correspondance de Hecke
$$
c=(c_{1},c_{2}):\mathop{\rm Cht}\nolimits_{N}^{r}(g)/a^{{\Bbb
Z}}\rightarrow (\mathop{\rm Cht}\nolimits_{N}^{r}/a^{{\Bbb
Z}})\times_{(X-N)^{2}}(\mathop{\rm Cht} \nolimits_{N}^{r}/a^{{\Bbb
Z}})
$$
envoie $\mathop{\rm Cht} \nolimits_{N}^{r;\,\leq p}/a^{{\Bbb Z}}$ dans
$\mathop{\rm Cht}\nolimits_{N}^{r;\,\leq q}/a^{{\Bbb Z}}$ d\`{e}s que
$q-p$ est assez convexe en fonction de $K_{N}gK_{N}$. Par
cons\'{e}quent, le ind-syst\`{e}me local $\ell$-adique
$H_{N}^{r;\,2r-2}$ d\'{e}fini en 6.4 est muni d'une action de la
$\overline{{\Bbb Q}}_{\ell}$-alg\`{e}bre de Hecke ${\cal H}(K_{N})=
{\cal C}_{{\rm c}} (\mathop{\rm GL}\nolimits_{r}({\Bbb A})//K_{N})$
des fonctions \`{a} valeurs dans $\overline{{\Bbb Q}}_{\ell}$,
$K_{N}$-invariantes \`{a} gauche et \`{a} droite et \`{a} supports
compacts sur $\mathop{\rm GL} \nolimits_{r}({\Bbb A})$.

De m\^{e}me, les morphismes de Frobenius partiels induisent un
isomorphisme
$$
(\mathop{\rm Frob}\nolimits_{X-N}\times\mathop{\rm Id}
\nolimits_{X-N})^{\ast}H_{N}^{r;\,2r-2} \buildrel\sim\over
\longrightarrow H_{N}^{r;\,2r-2}
$$
qui commute \`{a} l'action de ${\cal H}(K_{N})$.

Comme la construction faite en 6.4 du syst\`{e}me local $\ell$-adique
$H_{N,{\rm cusp}}^{r}$ \`{a} partir de $H_{N}^{r;\,2r-2}$ est
canonique, $H_{N,{\rm cusp}}^{r}$ est lui aussi muni d'une action de
l'alg\`{e}bre de Hecke ${\cal H}(K_{N})$ et d'un isomorphisme
$(\mathop{\rm Frob}\nolimits_{X-N}\times\mathop{\rm Id}
\nolimits_{X-N})^{\ast}H_{N,{\rm cusp}}^{r}\buildrel\sim\over
\longrightarrow H_{N,{\rm cusp}}^{r}$. Dans l'\'{e}nonc\'{e} de la
derni\`{e}re proposition de la section 6.4, on peut donc remplacer la
moyenne ${1\over r!}\sum_{k=1}^{r!} [(\mathop{\rm Frob}
\nolimits_{X-N}^{k}\times\mathop{\rm Id}\nolimits_{X-N})^{\ast}
H_{N,{\rm cusp}}^{r}]$ par $[H_{N,{\rm cusp}}^{r}]$.

Dans cette derni\`{e}re section nous allons montrer:

\thm PROPOSITION 
\enonce
Pour tout $g\in\mathop{\rm GL}\nolimits_{r}({\Bbb A})$, tout couple
$(\infty ,o)$ de points ferm\'{e}s distincts de $X-S_{g}$ {\rm
(}o\`{u} $S_{g}\subset |N|$ est l'ensemble fini de points ferm\'{e}s
de $X$ d\'{e}fini en {\rm 5.1)}, tout point ferm\'{e} $\xi$ de
$\infty\times o$ et tout entier $n$, on a
$$\displaylines{
\qquad\mathop{\rm Tr}(1_{K_{N}gK_{N}}\times\mathop{\rm 
Frob}\nolimits_{\xi}^{n},H_{N,{\rm cusp}}^{r})
\hfill\cr\hfill
=q^{(r-1)\mathop{\rm deg}(\xi )n}\kern -2mm\sum_{{\scriptstyle\pi\in
{\cal A}_{r}(K_{N})\atop\scriptstyle \omega_{\pi}(a)=1}}\mathop{\rm
Tr}\nolimits_{\pi}(1_{K_{N}gK_{N}})S_{\infty}^{\left(-{\mathop{\rm
deg}(\xi )\over \mathop{\rm deg}(\infty )}n\right)}\!(\pi )
S_{o}^{\left({\mathop{\rm deg}(\xi )\over\mathop{\rm
deg}(o)}n\right)}\!(\pi ).\qquad}
$$
\endthm

Cette proposition g\'{e}n\'{e}ralise la partie (iii) de la
derni\`{e}re proposition de la sous-section 6.4 et des arguments
standard permettent d'en d\'{e}duire les propri\'{e}t\'{e}s ${\rm
P}_{1}(r)$, ${\rm P}_{2}(r)$ et ${\rm P}_{3}(r)$.
\vskip 5mm

{\bf 8.2. Trace des op\'{e}rateurs de Hecke sur $\widetilde{H}_{N,{\rm
cusp}}^{\,r;\,\leq p}$}
\vskip 5mm

\thm LEMME 
\enonce
Soient $U$ un ouvert non vide de $X$ et $M$ un syst\`{e}me local
$\ell$-adique sur $U\times U$ muni d'un endomorphisme
$h:M\rightarrow M$. Alors, un fois fix\'{e} un ensemble de
repr\'{e}sentants des classes d'isomorphie de syst\`{e}mes locaux
$\ell$-adiques irr\'{e}ductibles sur $U\times U$, il existe une unique
application \`{a} support fini $L\mapsto\lambda_{L}$ de cet ensemble
dans $\overline{{\Bbb Q}}_{\ell}$ telle que
$$
\mathop{\rm Tr}(h\circ\gamma ,M_{\overline{\delta}})=
\sum_{L}\lambda_{L} \mathop{\rm Tr}(\gamma ,L_{\overline{\delta}})
$$
pour tout $\gamma\in\pi_{1}(U\times U,\overline{\delta})$.
\endthm

Si $[M]=\sum_{L}m_{M}(L)[L]$ est le syst\`{e}me local $\ell$-adique
virtuel sur $U\times U$ associ\'{e} \`{a} $M$ et si $M_{{\rm
cusp}}=\sum_{L}m_{M_{{\rm cusp}}}(L)[L]$ est son morceau non
$r$-n\'{e}gligeable (cf. 6.4), on notera
$$
\mathop{\rm Tr}\nolimits_{M_{{\rm cusp}}}(h\circ\gamma)=
\sum_{{\scriptstyle L\atop\scriptstyle m_{M}(L)\not=0}}{m_{M_{{\rm
cusp}}}(L)\over m_{M}(L)}\lambda_{L}\mathop{\rm Tr}(\gamma
,L_{\overline{\delta}})
$$
pour tout $\gamma\in\pi_{1}(U\times U,\overline{\delta})$.
\vskip 3mm

Bien entendu, on va appliquer ceci au morceau non $r$-n\'{e}gligeable
$\widetilde{H}_{N,{\rm cusp}}^{\,r;\,\leq p}$ de la cohomologie
$\ell$-adique $\widetilde{H}_{N}^{\,r;\,\leq p,\,\ast}$ et \`{a}
l'endomorphisme de $\widetilde{H}_{N}^{\,r;\,\leq p,\,2r-2}$ induit
par une corres\-pondance de Hecke.

Plus pr\'{e}cis\'{e}ment, soit $p:[0,r]\rightarrow {\Bbb R}$ un
param\`{e}tre de troncature assez convexe par rapport au genre de la
courbe $X$ et \`{a} $N$ et soit $g\in \mathop{\rm GL}\nolimits_{r}
({\Bbb A})$. D'une part, on a d\'{e}j\`{a} vu que la correspondance de
Hecke $c:\mathop{\rm Cht} \nolimits_{N}^{r}(g)\rightarrow \mathop{\rm
Cht}\nolimits_{N}^{r} \times_{(X-N)^{2}}\mathop{\rm
Cht}\nolimits_{N}^{r}$ envoie $\mathop{\rm Cht}\nolimits_{N}^{r;\,\leq
p}$ dans $\mathop{\rm Cht}\nolimits_{N}^{r;\,\leq q}$ pour un
param\`{e}tre de troncature $q\geq p$ convenable et qu'elle induit
un morphisme entre les cohomologies \`{a} supports compacts
$h:H_{N}^{r;\,\leq p,\,\ast} \rightarrow H_{N}^{r;\,\leq q,\,\ast}$.
D'autre part, cette correspondance induit par restriction un cycle
$$
c^{\leq p}:\mathop{\rm Cht}\nolimits_{N}^{r;\,\leq p}(g)\rightarrow
\mathop{\rm Cht}\nolimits_{N}^{r;\,\leq p}\times_{(X-N)^{2}}
\mathop{\rm Cht}\nolimits_{N}^{r;\,\leq p},
$$
puis par normalisation une correspondance
$$
\widetilde{c}^{\leq p}:\mathop{\widetilde{\rm Cht}}\nolimits_{N}^{r;
\,\leq p}(g)\rightarrow \mathop{\widetilde{\rm Cht}}
\nolimits_{N}^{r;\,\leq p}\times_{(X-N)^{2}} \mathop{\widetilde{\rm
Cht}}\nolimits_{N}^{r;\,\leq p},
$$
et donc un endomorphisme not\'{e} $\widetilde{h}$ de la cohomologie
$\ell$-adique $\widetilde{H}_{N}^{r;\,\leq p,\,\ast}$, d'o\`{u} une
trace
$$
\mathop{\rm Tr}\nolimits_{\widetilde{H}_{N,{\rm cusp}}^{\,r;\,\leq p}}
(\widetilde{h}\circ\gamma)
$$
pour tout $\gamma\in\pi_{1}((X-N)^{2},\overline{\delta})$.

Quitte \`{a} agrandir $q$, on peut supposer que $q$ est assez convexe
par rapport au genre de la courbe $X$ et \`{a} $N$. On alors un
diagramme commutatif
$$\diagram{
H_{N}^{r;\,\leq q,\,\ast}&\kern -1mm\smash{\mathop{\hbox to
37mm{\rightarrowfill}}\limits^{\scriptstyle }}\kern -1mm
&\widetilde{H}_{N}^{r;\,\leq q,\,\ast}\cr
\llap{$\scriptstyle h$}\left\uparrow \vbox to
4mm{}\right.\rlap{}&+&\llap{}\left\downarrow \vbox to
4mm{}\right.\rlap{$\scriptstyle $}\cr
H_{N}^{r;\,\leq p,\,\ast}&\kern -1mm\smash{\mathop{\hbox to
8mm{\rightarrowfill}} \limits_{\scriptstyle }}\kern 2mm
\widetilde{H}_{N}^{r;\,\leq p,\,\ast}\kern 2mm\smash{\mathop{\hbox to
8mm{\rightarrowfill}} \limits_{\scriptstyle \widetilde{h}}}\kern -1mm
&\widetilde{H}_{N}^{r;\,\leq p,\,\ast}\cr}
$$
o\`{u} les fl\`{e}ches $h$ et $\widetilde{h}$ sont celles d\'{e}finies
ci-dessus et les autres fl\`{e}ches sont les m\^{e}mes que dans le
carr\'{e} commutatif de 6.4. Par un argument facile utilisant les
r\'{e}sultats de 6.4, Lafforgue en d\'{e}duit:

\thm PROPOSITION 
\enonce
Pour tout param\`{e}tre de troncature $p$ assez convexe en fonction 
de $N$, tout point ferm\'{e} $\xi$ de $(X-N)^{2}$ et tout entier $n$, on 
a
$$
\mathop{\rm Tr}(1_{K_{N}gK_{N}}\times\mathop{\rm Frob}
\nolimits_{\xi}^{n},H_{N,{\rm cusp}}^{r})=\mathop{\rm Tr}
\nolimits_{\widetilde{H}_{N,{\rm cusp}}^{\,r;\,\leq p}}(\widetilde{h}
\circ\mathop{\rm Frob}\nolimits_{\xi}^{n}).
$$
\endthm
\vskip 5mm

{\bf 8.3. Application de la variante du th\'{e}or\`{e}me de Pink}
\vskip 5mm

Fixons un param\`{e}tre de troncature $p:[0,r]\rightarrow {\Bbb R}$
assez convexe par rapport au genre de la courbe $X$ et \`{a} $N$, et 
un \'{e}l\'{e}ment $g\in\mathop{\rm GL}\nolimits_{r}({\Bbb A})$. On a 
d\'{e}fini un ensemble fini $S_{g}\supset |N|$ de points ferm\'{e}s 
de $X$.

On cherche \`{a} appliquer les r\'{e}sultats g\'{e}n\'{e}raux du
chapitre 7 \`{a}
$$
S=(X-S_{g})^{2},~f=\mathop{\rm Id}\nolimits_{S}:S\rightarrow S,
$$
$$
{\cal X}_{\emptyset}=(X-S_{g})^{2}\times_{(X-N)^{2}}(\mathop{\rm
Cht}\nolimits_{N}^{r;\,\leq p}/a^{{\Bbb Z}})\subset (X-S_{g})^{2}
\times_{(X-N)^{2}}(\mathop{\widetilde{\rm Cht}}
\nolimits_{N}^{\,r;\,\leq p}/a^{{\Bbb Z}})={\cal X},
$$
$$
(c_{1},c_{2}):(X-S_{g})^{2}\times_{(X-N)^{2}}(\mathop{\rm Cht}
\nolimits_{N}^{r}(g)/a^{{\Bbb Z}})=\Gamma_{\emptyset}\rightarrow
{\cal X}_{\emptyset}\times_{S}{\cal X}_{\emptyset},
$$
et
$$
\Gamma\rightarrow {\cal X}\times_{S}{\cal X}
$$
la normalisation de ${\cal X}_{\emptyset}\times_{S}{\cal
X}_{\emptyset}$ dans $\Gamma_{\emptyset}$. Le ferm\'{e}
compl\'{e}mentaire ${\cal Y}={\cal X}-{\cal X}_{\emptyset}$ est un
diviseur \`{a} croisements normaux relatif sur $(X-S_{g})^{2}$, qui
est r\'{e}union ${\cal Y}={\cal Y}_{1}\cup\cdots\cup {\cal Y}_{m}$ de
diviseurs lisses sur $(X-S_{g})^{2}$.

Pour chaque entier $n_{0}\geq 0$, on dispose des classes de
cohomologie suivantes:
\vskip 2mm

\itemitem{-} $\mathop{\rm cl}(\Gamma )$ dans la cohomologie de ${\cal
X}\times_{S}{\cal X}$ au-dessus de $(X-S_{g})^{2}$,
\vskip 2mm

\itemitem{-} $\mathop{\rm cl}(\Gamma )_{I}^{(n_{0})}$ dans la
cohomologie de ${\cal Y}_{I}\times_{S}{\cal Y}_{I}$ au-dessus de
$(X-S_{g})^{2}$ pour tout $I\subset\{1,\ldots ,m\}$ non vide.
\vskip 2mm

Pour tout ouvert $U$ de $X$ notons $\Lambda_{U}$ le sch\'{e}ma
intersection de tous les ouverts de $U\times U$ compl\'{e}mentaires
des images r\'{e}ciproques de la diagonales par les endomorphismes
$\mathop{\rm Frob}\nolimits_{U}^{n}\times\mathop{\rm Id}\nolimits_{U}$
et $\mathop{\rm Id}\nolimits_{U} \times\mathop{\rm Frob}
\nolimits_{U}^{n}$ de $U\times U$ pour tous les entiers $n\geq 0$,
c'est-\`{a}-dire l'intersection de $\Lambda_{X}\subset X\times X$ avec
$U\times U$.

Nous avons besoin du th\'{e}or\`{e}me suivant de Lafforgue, qui
prolonge en rang arbitraire un r\'{e}sultat \'{e}tabli par Drinfeld en
rang $2$ ([Dr 7]).

\thm TH\'{E}OR\`{E}ME 
\enonce
Supposons que $p$ est assez convexe et $U\subset X-S_{g}$ est assez
petit relativement \`{a} $N$ et \`{a} $K_{N}gK_{N}$. Il existe alors
un ouvert de Zariski $V\subset U\times U$, contenant $\Lambda_{U}$,
tel que, apr\`{e}s avoir remplac\'{e} $S=(X-S_{g})^{2}$ par $V$ et
${\cal X}_{\emptyset}\subset {\cal X}$ et $\Gamma_{\emptyset} \subset
\Gamma$ par leurs restrictions \`{a} cet ouvert de $S$, la
correspondance $\Gamma$ stabilise l'ouvert ${\cal X}_{\emptyset}$ au
voisinage de ses points fixes.
\endthm

En d'autres termes, l'hypoth\`{e}se du th\'{e}or\`{e}me de la
sous-section 7.3 est v\'{e}rifi\'{e}e au-dessus de
$V\supset\Lambda_{U}$. Par cons\'{e}quent, si $n_{0}\geq 0$ est
l'entier de loc. cit., pour tout point g\'{e}om\'{e}trique
$\overline{\xi}$ localis\'{e} en un point ferm\'{e} $\xi$ de
$\Lambda_{U}$ et tout entier $n$ tel que $\mathop{\rm deg}(\xi
)n>n_{0}$, on a la formule de points fixes
$$\eqalign{
\mathop{\rm Lef}\nolimits_{\xi}(\Gamma_{\emptyset}\times\mathop{\rm
Frob}\nolimits_{{\cal X}_{\emptyset}}^{\mathop{\rm deg}(\xi )n}, {\cal
X}_{\emptyset})&=\mathop{\rm tr}(\mathop{\rm cl}(\Gamma)
\times\mathop{\rm Frob}\nolimits_{{\cal X}}^{\mathop{\rm deg}(\xi )n},
R\Gamma ({\cal X}_{\overline{\xi}},{\Bbb Q}_{\ell}))\cr
&\kern 4mm +\sum_{I\not=\emptyset}(-1)^{|I|}\mathop{\rm tr}
(\mathop{\rm cl}(\Gamma )_{I}^{(n_{0})}\times\mathop{\rm Frob}
\nolimits_{{\cal Y}_{I}}^{\mathop{\rm deg}(\xi )n},R\Gamma ({\cal
Y}_{I,\overline{\xi}}, {\Bbb Q}_{\ell})).\cr}
$$

On remarquera que, pour chaque $I\subset\{1,\ldots ,m\}$ non vide, il
existe une famille $(m_{I}(L))_{L}$ de scalaires dans $\overline{{\Bbb
Q}}_{\ell}$, index\'{e}e par un ensemble de repr\'{e}sentants des
syst\`{e}mes locaux $\ell$-adiques irr\'{e}ductibles sur
$(X-S_{g})^{2}$, \`{a} support fini, telle que
$$
\mathop{\rm tr}(\mathop{\rm cl}(\Gamma )_{I}^{(n_{0})}\times
\mathop{\rm Frob}\nolimits_{{\cal Y}_{I}}^{\mathop{\rm deg}(\xi )n},
R\Gamma ({\cal Y}_{I,\overline{\xi}},{\Bbb Q}_{\ell}))=
\sum_{L}m_{I}(L)\mathop{\rm Tr}(\mathop{\rm Frob}
\nolimits_{\xi}^{n},L),
$$
et que $m_{I}(L)=0$ pour tout $L$ qui n'est pas $r$-n\'{e}gligeable
d'apr\`{e}s 6.3.
\vskip 3mm

{\it Dans la suite de cette sous-section, on identifie
$\overline{{\Bbb Q}}_{\ell}$ \`{a} ${\Bbb C}$ par l'isomorphisme
fix\'{e} en} 1.3.
\vskip 3mm

Dans la formule de points fixes ci-dessus, on peut remplacer
$\overline{\xi}$ par $(\mathop{\rm Frob}\nolimits_{X}^{k}
\times\mathop{\rm Id}\nolimits_{X})(\overline{\xi})$ pour
$k=1,\ldots ,r!$ et faire la moyenne des $r!$ expressions
obtenues. Le premier membre de cette moyenne est \'{e}gal \`{a}
$$\eqalign{
&\sum_{{\scriptstyle\pi\in {\cal A}_{r}(K_{N})\atop\scriptstyle
\omega_{\pi}(a)=1}}\mathop{\rm Tr}\nolimits_{\pi}(1_{K_{N}gK_{N}})
q^{(r-1)\mathop{\rm deg}(\xi )n}S_{\infty}^{\left(-{\mathop{\rm
deg}(\xi )\over \mathop{\rm deg}(\infty )}n\right)}\!(\pi )
S_{o}^{\left({\mathop{\rm deg}(\xi )\over \mathop{\rm deg}(o)}n
\right)}\!(\pi )\cr
&+\kern -2mm\sum_{{\scriptstyle 1\leq r'<r\atop\scriptstyle 1\leq
r''<r}}\sum_{{\scriptstyle\pi'\in {\cal A}_{r'}(K_{N})\atop
\scriptstyle\pi''\in {\cal A}_{r''}(K_{N})}}\mathop{\rm Tr}
\nolimits_{\pi',\pi''}^{\leq_{\alpha}p}(1_{K_{N}gK_{N}},\mathop{\rm
deg}(\xi )n)S_{\infty}^{\left(-{\mathop{\rm deg}(\xi )\over
\mathop{\rm deg}(\infty )}n\right)}\! (\pi')S_{o}^{\left({\mathop{\rm
deg}(\xi )\over \mathop{\rm deg}(o)}n\right)}\!(\pi'')\cr}
$$
d'apr\`{e}s 5.2. Compte tenu de la proposition de la sous-section
pr\'{e}c\'{e}dente et de l'hypoth\`{e}se de r\'{e}currence ${\rm
P}_{1}(r')$ pour tout $r'<r$ (cf. 6.2), on a donc une expression pour
la diff\'{e}rence
$$\displaylines{
\qquad\mathop{\rm Tr}(1_{K_{N}gK_{N}}\times\mathop{\rm 
Frob}\nolimits_{\xi}^{n},H_{N,{\rm cusp}}^{r})
\hfill\cr\hfill
-q^{(r-1)\mathop{\rm deg}(\xi )n}\kern -2mm\sum_{{\scriptstyle\pi\in
{\cal A}_{r}(K_{N})\atop\scriptstyle \omega_{\pi}(a)=1}}\mathop{\rm
Tr}\nolimits_{\pi}(1_{K_{N}gK_{N}})S_{\infty}^{\left(-{\mathop{\rm
deg}(\xi )\over \mathop{\rm deg}(\infty )}n\right)}\!(\pi )
S_{o}^{\left({\mathop{\rm deg}(\xi )\over \mathop{\rm
deg}(o)}n\right)}\!(\pi )\qquad}
$$
de la forme
$$
\sum_{L}P_{L}(n\mathop{\rm deg}(\xi ))\mathop{\rm Tr}(\mathop{\rm
Frob}\nolimits_{\xi}^{n},L)
$$
o\`{u} $(\nu\mapsto P_{L}(\nu))_{L}$ est une famille \`{a} support
fini de fonctions complexes de l'entier $\nu$ de la forme
$$
P_{L}(\nu)=\sum_{z}P_{L,z}(\nu)z^{\nu}
$$
pour une famille \`{a} support fini de polyn\^{o}mes $P_{L,z}(u)\in
{\Bbb C}[u]$, $z\in {\Bbb C}^{\times}$, qui sont identiquement nuls 
quand $L$ n'est pas $r$-n\'{e}gligeable.

En d\'{e}veloppant les polyn\^{o}mes $P_{L,z}(u)$ et en utilisant
l'ind\'{e}pendance lin\'{e}aire des fonctions $\nu\mapsto
\nu^{k}z^{\nu}$ de l'entiers $\nu$ pour $k\in {\Bbb Z}$ et $z\in {\Bbb
C}^{\times}$, on voit tout d'abord que l'on ne perd rien en
repla\c{c}ant chaque $P_{L,z}(u)$ par son terme constant $P_{L,z}(0)$.

En utilisant la correspondance de Langlands et la conjecture de
Ramanujan-Petersson d\'{e}j\`{a} connues en les rangs $<r$ et un
dernier argument de fonctions $L$, Lafforgue conclut alors la preuve de
la proposition de 8.1 et la r\'{e}currence.

\vfill\eject

\newtoks\ref	\newtoks\auteur
\newtoks\titre	\newtoks\editeur
\newtoks\annee	\newtoks\revue
\newtoks\tome	\newtoks\pages
\newtoks\reste	\newtoks\autre

\def\bibitem#1{\parindent=20pt\itemitem{#1}\parindent=12pt}

\def\livre{\bibitem{[\the\ref]}%
\the\auteur ~-- {\sl\the\titre}, \the\editeur, ({\the\annee}).
\smallskip\smallskip\filbreak}

\def\article{\bibitem{[\the\ref]}%
\the\auteur ~-- \the\titre, {\sl\the\revue} {\the\tome},
({\the\annee}), \the\pages.\smallskip\filbreak}

\def\autre{\bibitem{[\the\ref]}%
\the\auteur ~-- \the\reste.\smallskip\filbreak}

\centerline{\twelvebf  BIBLIOGRAPHIE}
\vskip 10mm

{\bf Travaux de L. Lafforgue}
\vskip 5mm

\ref={La 1}
\auteur={L. {\pc LAFFORGUE}}
\titre={Chtoucas de Drinfeld et conjecture de Ramanujan-Petersson}
\editeur={Ast\'erisque, 243}
\annee={1997}
\livre

\ref={La 2}
\auteur={L. {\pc LAFFORGUE}}
\titre={Sur la conjecture de Ramanujan-Petersson pour les corps de
fonctions. I: \'{E}tude g\'{e}om\'{e}trique}
\revue={C. R. Acad. Sci. Paris}
\tome={322}
\annee={1996}
\pages={605-608}
\article

\ref={La 3}
\auteur={L. {\pc LAFFORGUE}}
\titre={Sur la conjecture de Ramanujan-Petersson pour les corps de
fonctions. II: \'{E}tude spectrale}
\revue={C. R. Acad. Sci. Paris}
\tome={322}
\annee={1996}
\pages={707-710}
\article

\ref={La 4}
\auteur={L. {\pc LAFFORGUE}}
\titre={Sur la d\'{e}g\'{e}n\'{e}rescence des chtoucas de Drinfeld}
\revue={C. R. Acad. Sci. Paris}
\tome={323}
\annee={1996}
\pages={491-494}
\article

\ref={La 5}
\auteur={L. {\pc LAFFORGUE}}
\titre={Compactification de l'isog\'{e}nie de Lang et
d\'{e}g\'{e}n\'{e}rescence des structures de niveau simple des
chtoucas de Drinfeld}
\revue={C. R. Acad. Sci. Paris}
\tome={325}
\annee={1997}
\pages={1309-1312}
\article

\ref={La 6}
\auteur={L. {\pc LAFFORGUE}}
\reste={Chtoucas de Drinfeld et applications, dans {\it Proceedings of
the International Congress of Mathematicians, Vol. II {\rm (}Berlin,
1998{\rm )}}, Doc. Math., (1998), 563-570}
\autre

\ref={La 7}
\auteur={L. {\pc LAFFORGUE}}
\titre={Une compactification des champs classifiant les chtoucas de
Drinfeld}
\revue={J. Amer. Math. Soc.}
\tome={11}
\annee={1998}
\pages={1001-1036}
\article

\ref={La 8}
\auteur={L. {\pc LAFFORGUE}}
\titre={Pavages des simplexes, sch\'{e}mas de graphes recoll\'{e}s et
compactification des $\mathop{\rm PGL}\nolimits_{r}^{n+1}/\mathop{\rm
PGL}\nolimits_{r}$}
\revue={Invent. Math.}
\tome={136}
\annee={1999}
\pages={233-271}
\article

\ref={La 9}
\auteur={L. {\pc LAFFORGUE}}
\reste={Compactification des $\mathop{\rm PGL}\nolimits_{r}^{n+1}/
\mathop{\rm PGL}\nolimits_{r}$ et restriction des scalaires \`{a} la
Weil, \`{a} para\^{\i}tre dans le volume du Cinquantenaire des Annales
de l'Institut Fourier}
\autre

\ref={La 10}
\auteur={L.  {\pc LAFFORGUE}}
\reste={La correspondance de Langlands sur les corps de fonctions,
manuscrit (en cours de frappe), juin 1999}
\autre

\vskip 5mm

{\bf Articles de V.G. Drinfeld sur les chtoucas}
\vskip 5mm

\ref={Dr 1}
\auteur={V. G. {\pc DRINFELD}}
\titre={Commutative subrings of certain noncommutative rings}
\revue={Funct. Anal. and its Appl.}
\tome={11}
\annee={1977}
\pages={9-12}
\article

\ref={Dr 2}
\auteur={V. G. {\pc DRINFELD}}
\titre={Proof of the global Langlands conjecture for $\mathop{\rm
GL}(2)$ over a function field}
\revue={Funct. Anal. and its Appl.}
\tome={11}
\annee={1977}
\pages={223-225}
\article

\ref={Dr 3}
\auteur={V. G. {\pc DRINFELD}}
\titre={A proof of Petersson's conjecture for function fields}
\revue={Uspehi Mat. Nauk}
\tome={32}
\annee={1977}
\pages={209-210}
\article

\ref={Dr 4}
\auteur={V. G. {\pc DRINFELD}}
\reste={Langlands' conjecture for $\mathop{\rm GL}(2)$ over functional fields,
dans {\it Proceedings of the International Congress of Mathematicians
{\rm (}Helsinki, 1978{\rm )}}, Acad. Sci. Fennica, (1980), 565-574}
\autre

\ref={Dr 5}
\auteur={V. G. {\pc DRINFELD}}
\titre={Moduli varieties of $F$-sheaves}
\revue={Funct. Anal. and its Appl.}
\tome={21}
\annee={1987}
\pages={107-122}
\article

\ref={Dr 6}
\auteur={V. G. {\pc DRINFELD}}
\titre={Proof of the Petersson conjecture for $\mathop{\rm GL}(2)$
over a global field of characteristic $p$}
\revue={Funct. Anal. and its Appl.}
\tome={22}
\annee={1988}
\pages={28-43}
\article

\ref={Dr 7}
\auteur={V. G. {\pc DRINFELD}}
\titre={Cohomology of compactified moduli varieties of $F$-sheaves of
rank $2$}
\revue={J. of Soviet Math.}
\tome={46}
\annee={1989}
\pages={1789-1821}
\article

\ref={Dr 8}
\auteur={V. G. {\pc DRINFELD}}
\titre={On my paper: ``Cohomology of compactified moduli varieties of
$F$-sheaves of rank $2$''}
\revue={Zap. Nauchn. Sem. Leningrad. Otdel. Mat. Inst. Steklov.
{\rm (}LOMI{\rm )}, Anal. Teor. Chisel i Teor. Funktsii. 9}
\tome={168}
\annee={1988}
\pages={45--47}
\article

\vskip 5mm

{\bf Autres publications sur la correspondance
de Langlands}
\vskip 5mm

\ref={Ca}
\auteur={H. {\pc CARAYOL}}
\reste={Vari\'{e}t\'{e}s de Drinfeld compactes, d'apr\`{e}s Laumon,
Rapoport et Stuhler, dans {\it S\'{e}minaire Bourbaki {\rm 1991/92}},
Ast\'{e}risque 206, (1992), 369-409}
\autre

\ref={D-H}
\auteur={P. {\pc DELIGNE}, D. {\pc HUSEMOLLER}}
\reste={Survey of Drinfeld modules, dans {\it Current trends in
arithmetical algebraic geometry {\rm (}Arcata, 1985{\rm )}}, Contemp.
Math. 67, (1987), 25-91}
\autre

\ref={Dr 9}
\auteur={V.G. {\pc DRINFELD}}
\titre={Elliptic modules}
\revue={Math. USSR Sbornik}
\tome={23}
\annee={1974}
\pages={561-592}
\article

\ref={Dr 10}
\auteur={V.G. {\pc DRINFELD}}
\titre={Elliptic modules II}
\revue={Math. USSR Sbornik}
\tome={31}
\annee={1977}
\pages={159-170}
\article

\ref={Dr 11}
\auteur={V. G. {\pc DRINFELD}}
\titre={Number of two-dimensional irreducible representations of
the fundamental group of a curve over a finite field}
\revue={Functional Anal. and its Appl.}
\tome={15}
\annee={1982}
\pages={294-295}
\article

\ref={F-K}
\auteur={Y. {\pc FLICKER}, D. {\pc KAZHDAN}}
\reste={Geometric Ramanujan conjecture and Drinfeld reciprocity
law, dans {\it Number theory, trace formulas and discrete groups},
Academic Press, (1989), 201-218}
\autre

\ref={H-K}
\auteur={G. {\pc HARDER}, D. {\pc KAZHDAN}}
\titre={Automorphic forms on $\mathop{\rm GL}\nolimits_{2}$ over
function fields (after V. G. Drinfeld)}
\revue={Proc. Sym. Pure Math.}
\tome={33, Part~2}
\annee={1979}
\pages={357-379}
\article

\ref={Ka}
\auteur={D. {\pc KAZHDAN}}
\titre={An introduction to Drinfeld's ``shtuka''}
\revue={Proc. Sym. Pure Math.}
\tome={33, Part~2}
\annee={1979}
\pages={347-356}
\article

\ref={Lau 1}
\auteur={G. {\pc LAUMON}}
\titre={Cohomology of Drinfeld modular varieties, Part I {\rm
(}Geometry, counting of points and local harmonic analysis{\rm )}}
\editeur={Cambridge University Press}
\annee={1995}
\livre

\ref={Lau 2}
\auteur={G. {\pc LAUMON}}
\titre={Cohomology of Drinfeld modular varieties, Part II {\rm
(}Automorphic forms, trace formulas and Langlands correspondence{\rm
)}}
\editeur={Cambridge University Press}
\annee={1996}
\livre

\ref={Lau 3}
\auteur={G. {\pc LAUMON}}
\reste={Drinfeld Shtukas, dans {\it CIME Session ``Vector bundles on
Curves. New Directions'' {\rm (}Cetraro, juin 1995{\rm )}}, Lecture
Notes in Mathematics 1649, (1996), 50-109}
\autre

\ref={Lau 4}
\auteur={G. {\pc LAUMON}}
\reste={The Langlands Correspondence for Function Fields following
Laurent Lafforgue, dans {\it Current Developments in Mathematics
Conference {\rm (}Harvard University, novembre 1999{\rm )}}}
\autre

\ref={L-R-S}
\auteur={G.  {\pc LAUMON}, M. {\pc RAPOPORT} et U. {\pc
STUHLER}}
\titre={D-elliptic sheaves and the Langlands correspondence}
\revue={Inv. Math.}
\tome={113}
\annee={1993}
\pages={217-338} 
\article

\vskip 5mm

{\bf Sur la correspondance de Langlands
g\'{e}om\'{e}trique}
\vskip 5mm

\ref={Dr 12}
\auteur={V. G. {\pc DRINFELD}}
\titre={Two-dimensional $l$-adic representations of the fundamental group
of a curve over a finite field and automorphic forms on $\mathop{\rm 
GL}(2)$}
\revue={Amer. J. Math.}
\tome={105}
\annee={1983}
\pages={85-114}
\article

\ref={Dr 13}
\auteur={V. G. {\pc DRINFELD}}
\titre={Two-dimensional $l$-adic representations of the Galois group of a
global field of characteristic $p$ and automorphic forms on $\mathop{\rm 
GL}(2)$}
\revue={J. of Soviet Math.}
\tome={36}
\annee={1987}
\pages={93-105}
\article

\ref={F-G-K-V}
\auteur={E. {\pc FRENKEL}, D. {\pc GAITSGORY}, D. {\pc KAZHDAN}, K.
{\pc VILONEN}}
\titre={Geometric realization of Whittaker
functions and the Langlands conjecture} 
\revue={J. Amer. Math. Soc.}
\tome={11}
\annee={1998}
\pages={451-484}
\article

\ref={Lau 5}
\auteur={G. {\pc LAUMON}}
\titre={Correspondance de Langlands g\'{e}om\'{e}trique pour les corps
de fonctions}
\revue={Duke Math. J.}
\tome={54}
\annee={1987}
\pages={309-359}
\article

\ref={Lau 6}
\auteur={G. {\pc LAUMON}}
\reste={Faisceaux automorphes pour $\mathop{\rm GL}\nolimits_{n}$: la
premi\`{e}re construction de Drinfeld, pr\'{e}publication
\'{e}lectronique, r\'{e}f. alg-geom/9511004}
\autre

\vskip 5mm

{\bf Autres r\'{e}f\'{e}rences de g\'{e}om\'{e}trie alg\'{e}brique}
\vskip 5mm

\ref={De}
\auteur={P. {\pc DELIGNE}}
\titre={La conjecture de Weil. II} 
\revue={Publ. Math. IHES}
\tome={52}
\annee={1980}
\pages={137-252}
\article

\ref={DC-P}
\auteur={C. {\pc DE} {\pc CONCINI}, C. {\pc PROCESI}}
\reste={Complete symmetric varieties, dans {\it CIME Session
``Invariant Theory'' {\rm (}Montecatini{\rm )}}, Lecture Notes in
Math. 996, Springer-Verlag, (1982), 1-44}
\autre

\ref={Fa}
\auteur={G.  {\pc FALTINGS}}
\titre={Explicit resolution of local singularities of moduli-spaces}
\revue={Journal f\"ur die reine und angewandte Mathematik}
\tome={483}
\annee={1997}
\pages={183-196}
\article

\ref={Ge}
\auteur={A. {\pc GENESTIER}}
\reste={A toroidal resolution for the bad reduction of some Shimura 
varieties, http://xxx.lanl.gov/abs/math/9912054, (1999)}
\autre

\ref={Gr}
\auteur={A. {\pc GROTHENDIECK}}
\reste={Formule de Lefschetz et rationalit\'{e} des fonctions $L$,
S\'{e}m. Bourbaki 1964/65, exp.~n$^{\circ}$~279, collection hors
s\'{e}rie Ast\'{e}risque {\bf 9} (1995), 41-55, et dans {\it Dix
expos\'{e}s sur la cohomologie des sch\'{e}mas}, North Holland,
(1968), 31-45}
\autre

\ref={H-N}
\auteur={G. {\pc HARDER}, M.S. {\pc NARASIMHAN}}
\titre={On the cohomology groups of moduli spaces of vector bundles on
curves}
\revue={Math. Ann.} 
\tome={212}
\annee={1975}
\pages={215-248} 
\article

\ref={H-T}
\auteur={M. {\pc HARRIS}, R. {\pc TAYLOR}}
\reste={On the geometry and cohomology of some simple Shimura
varieties, Institut de Math\'{e}matiques de Jussieu,
Pr\'{e}publication 227, (1999)}
\autre

\ref={Lak}
\auteur={D. {\pc LAKSOV}}
\reste={Completed quadrics and linear maps, dans {\it Algebraic
geometry, Bowdoin, 1985 {\rm (}Brunswick, Maine, 1985{\rm )}}, Proc.
Sympos. Pure Math., 46, Part 2, Amer. Math. Soc., (1987), 371-387}
\autre

\ref={Lau 7}
\auteur={G.  {\pc LAUMON}}
\titre={Transformation de Fourier, constantes d'\'{e}quations
fonctionnelles et conjecture de Weil}
\revue={Publ. Math. I.H.E.S.}
\tome={65}
\annee={1987}
\pages={131-210}
\article

\ref={L-M}
\auteur={G. {\pc LAUMON}, L. {\pc MORET}-{\pc BAILLY}}
\titre={Champs alg\'{e}briques}
\editeur={Springer-Verlag}
\annee={1999}
\livre

\ref={Pi}
\auteur={R. {\pc PINK}}
\titre={On the calculation of local terms in the Lefschetz-Verdier 
trace formula and its application
to a conjecture of Deligne}
\revue={Ann. Math.}
\tome={135}
\annee={1992}
\pages={483-525} 
\article

\ref={SD}
\auteur={B. {\pc SAINT}-{\pc DONAT}}
\reste={Affine Embeddings, dans {\it Toroidal Embeddings I}, Lecture Notes in 
Math. 339, Springer-Verlag, (1973), 1-40}
\autre

\ref={Se}
\auteur={J.-P. {\pc SERRE}}
\titre={Abelian $\ell$-Adic Representations and Elliptic Curves}
\editeur={Addison-Wesley}
\annee={1988}
\livre

\vskip 5mm

{\bf Autres r\'{e}f\'{e}rences automorphes}
\vskip 5mm

\ref={Ar}
\auteur={J. {\pc ARTHUR}}
\titre={A trace formula for reductive groups I: terms associated to 
classes in $G({\Bbb Q})$} 
\revue={Duke Math. J.}
\tome={45}
\annee={1978}
\pages={911-953}
\article

\ref={C-PS}
\auteur={J. W. {\pc COGDELL}, I. I. {\pc PIATETSKI}-{\pc SHAPIRO}}
\titre={Converse theorems for $\mathop{\rm GL}\nolimits_{n}$} 
\revue={Publ. Math. IHES}
\tome={79}
\annee={1994}
\pages={157-214}
\article

\ref={G-J}
\auteur={R. {\pc GODEMENT}, H. {\pc JACQUET}}
\reste={Zeta functions of simple algebras, Lecture Notes in Math. 260,
Springer-Verlag, (1972)}
\autre

\ref={J-PS-S}
\auteur={H. {\pc JACQUET}, I. I {\pc PIATETSKI}-{\pc SHAPIRO}, J.A. {\pc
SHALIKA}}
\titre={Rankin-Selberg convolutions} 
\revue={Amer. J. Math.}
\tome={105}
\annee={1983}
\pages={367-464}
\article

\ref={J-S 1}
\auteur={H. {\pc JACQUET}, J.A. {\pc SHALIKA}}
\titre={On Euler products and the classification of automorphic 
representations I} 
\revue={Amer. J. Math.}
\tome={103}
\annee={1981}
\pages={499-558}
\article

\ref={J-S 2}
\auteur={H. {\pc JACQUET}, J.A. {\pc SHALIKA}}
\titre={On Euler products and the classification of automorphic 
representations II} 
\revue={Amer. J. Math.}
\tome={103}
\annee={1981}
\pages={777-815}
\article

\ref={Lan 1}
\auteur={R.P. {\pc LANGLANDS}}
\reste={Problems in the theory of automorphic forms, dans {\it
Lectures in modern analysis and applications III}, Lecture Notes in
Math. 170, Springer-Verlag, (1970), 18-61}
\autre

\ref={Lan 2}
\auteur={R.P. {\pc LANGLANDS}}
\titre={On the functional equations satisfied by Eisenstein series}
\editeur={Lecture Notes in Math. 544, Springer-Verlag}
\annee={1976}
\livre

\ref={M-W}
\auteur={C. {\pc MOEGLIN}, J.-L. {\pc WALDSPURGER}}
\titre={D\'{e}composition spectrale et s\'{e}ries d'Ei\-sen\-stein, Une
paraphrase de l'\'{E}criture}
\editeur={Birkh\"{a}user}
\annee={1994} 
\livre

\ref={PS 1}
\auteur={I.I. {\pc PIATETSKI}-{\pc SHAPIRO}}
\reste={Zeta functions of $\mathop{\rm GL}(n)$, pr\'{e}publication de
L'uni\-ver\-si\-t\'{e} du Maryland, (1976)}
\autre

\ref={PS 2}
\auteur={I.I. {\pc PIATETSKI}-{\pc SHAPIRO}}
\titre={Multiplicity one theorem}
\revue={Proc. Sym. Pure Math.}
\tome={33, Part~1}
\annee={1979}
\pages={209-212}
\article

\ref={Sh}
\auteur={J.A. {\pc SHALIKA}}
\titre={The multiplicity one theorem for $\mathop{\rm
GL}\nolimits_{n}$}
\revue={Ann. of Math.}
\tome={100}
\annee={1974}
\pages={171-193} 
\article

\vskip 7mm

\line{\hfill\hfill\hbox{{\vtop{\tabalign G\'{e}rard Laumon\cr
\tabalign Universit\'{e} de Paris-Sud\cr
\tabalign et CNRS, UMR 8628\cr
\tabalign Math\'ematiques, B\^{a}t. 425\cr
\tabalign F-91405 Orsay Cedex (France)\cr
\tabalign Gerard.Laumon@math.u-psud.fr\cr}}}}
\bye